%% file: ms.tex
\documentclass[times,final]{elsarticle}
\pdfoutput=1

\usepackage{jcomp}
\usepackage{framed,multirow}

\usepackage{amssymb}
\usepackage{latexsym}
\usepackage{bm}
 \usepackage{indentfirst}
\usepackage{amsmath,amsthm,amssymb}

\usepackage{graphicx}
\usepackage{subcaption}
\usepackage{float}
\usepackage[framemethod=tikz]{mdframed}

\usepackage{import}

\import{\commandsdir}{cmdFluxReconstruction.tex}

\import{\commandsdir}{cmdPartialDifferentialEquations.tex}
\import{\commandsdir}{cmdPhysicalQuantities.tex}

\import{\commandsdir}{cmdPunctuation.tex}

\import{\commandsdir}{cmdTurbulenceModeling.tex}
\import{\commandsdir}{cmdImmersedBoundary.tex}

\usepackage{url}
\usepackage{xcolor}
\definecolor{newcolor}{rgb}{.8,.349,.1}

\journal{Journal of Computational Physics}

\begin{document}

\verso{Jiaqing Kou \textit{etal}}

\begin{frontmatter}

\title{Eigensolution analysis of immersed boundary method based on volume penalization: applications to high-order schemes}

\author[1,2]{Jiaqing \snm{Kou}\corref{cor1}}
\cortext[cor1]{Corresponding author}
\ead{jiaqingkou@gmail.com}
\author[1,2]{Aurelio \snm{Hurtado-de-Mendoza}}
\author[1,2]{Saumitra \snm{Joshi}}
\author[1]{Soledad \snm{Le Clainche}}
\author[1,3]{Esteban \snm{Ferrer}}

\address[1]{ETSIAE-UPM-School of Aeronautics, Universidad Politécnica de Madrid, Plaza Cardenal Cisneros 3, E-28040 Madrid, Spain}
\address[2]{NUMECA International S.A., Chaussee de la Hulpe 187, Brussels, B-1170, Belgium}
\address[3]{Center for Computational Simulation, Universidad Politécnica de Madrid, Campus de Montegancedo, Boadilla del Monte, 28660 Madrid, Spain}

\received{1 May 2013}
\finalform{10 May 2013}
\accepted{13 May 2013}
\availableonline{15 May 2013}
\communicated{xxx}

\begin{abstract}
This paper presents eigensolution and non-modal analyses for immersed boundary methods (IBMs) based on volume penalization for the linear advection equation. This approach is used to analyze the behavior of flux reconstruction (FR) discretization, including the influence of polynomial order and penalization parameter on numerical errors and stability. Through a semi-discrete analysis, we find that the inclusion of IBM adds additional dissipation without changing significantly the dispersion of the original numerical discretization. This agrees with the physical intuition that in this type of approach, the solid wall is modelled as a porous medium with vanishing viscosity. From a stability point view, the selection of penalty parameter can be analyzed based on a fully-discrete analysis, which leads to practical guidelines on the selection of penalization parameter. Numerical experiments indicate that the penalization term needs to be increased to damp oscillations inside the solid (i.e. porous region), which leads to undesirable time step restrictions. As an alternative, we propose to include a second-order term in the solid for the no-slip wall boundary condition. Results show that by adding a second-order term we improve the overall accuracy with relaxed time step restriction. This indicates that the optimal value of the penalization parameter and the second-order damping can be carefully chosen to obtain a more accurate scheme. Finally, the approximated relationship between these two parameters is obtained and used as a guideline to select the optimum penalty terms in a Navier-Stokes solver, to simulate flow past a cylinder.
\end{abstract}


\end{frontmatter}

\tableofcontents
\section{Introduction}
Dispersion and dissipation behavior is of paramount importance for all numerical schemes that solve Partial Differential Equations (PDEs). This behavior can be characterised using eigensolution analysis, which quantifies how the amplitude and frequency of a wave-like solution evolve with time, providing useful information on the numerical errors. Eigensolution analysis, also known as von Neumann or Fourier analysis, has been widely applied to different numerical schemes for space discretization, including finite difference \citep{mahesh1998family,hirsch2007numerical}, finite volume \citep{leveque2002finite}, finite element \citep{hughes2012finite}, Lattice Boltzmann \citep{chavez2018improving,chavez2020optimizing}, as well as high-order methods \citep{hu1999analysis,van2007dispersion,gassner2011comparison,vincent2011insights,moura2015linear,alhawwary2018fourier,manzanero2018dispersion}. Typically, eigendecomposition is applied to the global discretization matrix, where dispersion and dissipation characteristics reflect in the corresponding eigenvalues. The dispersion error, which is the error in wave advection, is represented by the modified frequencies of the solution. The dissipation error, which corresponds to the nonphysical damping or amplification effect, is represented by the modified amplitudes of the solution. This dispersion-dissipation behavior is crucial to understand and quantity the numerical error of a numerical scheme, and can be used to evaluate its stability and robustness. In addition, insights from the analysis can be extended to design better numerical schemes for turbulence simulation \citep{moura2016eigensolution,manzanero2020design,solan2021application}. Eigensolution analysis in multi-dimensions can also be used to evaluate the effect of mesh quality for high-order schemes \citep{trojak2020effect}. Most analyses belong to temporal eigensolution analysis, since the focus of such analysis is on the temporal evolution of the solution by considering periodic boundary conditions. To investigate the evolution in space, spatial eigensolution analysis \citep{mengaldo2018spatial,mengaldo2018spatial2} has been proposed recently for high-order schemes with inflow-outflow boundary conditions. Both spatial and temporal eigensolution analyses have been extended in a non-intrusive and data-driven manner to evaluate dispersion-dissipation behavior only from simulation data \citep{kou2021STKD}. Another alternative perspective for the dissipation error is to look at the non-modal behavior characterized by the short-term dissipation of a numerical scheme \citep{fernandez2019non}, which has been shown to obtain useful insights for the hybridized Discontinuous Galerkin (DG) scheme. When the time integration scheme is incorporated into the analysis, a fully-discrete analysis \citep{yang2013dispersion,vermeire2017behaviour} is performed and takes into account the numerical error both in space and time. 

So far, eigensolution analyses are still limited to body-fitted grids. The immersed boundary method (IBM) \citep{mittal2005immersed,sotiropoulos2014immersed,griffith2020immersed} provides an alternative to body-fitted meshes. This technique has gained attention during recent years and has the potential to handle complex geometry and moving bodies on simple Cartesian grids, thus reducing the cost in mesh generation. IBM originates from the work of Peskin \citep{peskin1972IBM} to introduce a singular source term to the background flow grid in the vicinity of the solid body. So far, IBM been extensively studied in various applications including benchmark geometries \citep{taira2007immersed,wu2009implicit}, turbulent flows \citep{iaccarino2003immersed,tamaki2017near}, fluid-structure interaction (FSI) \citep{huang2007simulation,yang2012simple,tian2014fluid}, acoustics \citep{seo2011high,sun2012immersed}, multiphase flows \citep{wang2017immersed,o2018volume}, etc. To understand the numerical behaviors and better implement IBM treatment in practical numerical simulation, dispersion-dissipation analysis can be adapted to IBM, which is a subject that has been deeply explored in the present study.

In general, IBM can be achieved in multiple ways, which are broadly divided into two strategies, including cutting the cell by the solid boundary \citep{ye1999accurate,udaykumar2001sharp} or introducing additional forcing terms  \citep{goldstein1993modeling,angot1999penalization,fadlun2000combined,luo2012numerical} to mimic the effect of solid objects. As a typical approach of the second kind, the volume penalization method \citep{angot1999penalization,brown2014CBVP,abgrall2014IBM,schneider2015immersed} has been widely adopted for boundary treatment due to its good robustness, simplicity and rigorous theoretical background. It is based on the physical intuition that the solid wall can be modelled as a porous medium with vanishing diffusivity \citep{kadoch2012volume}. A source or penalization term is introduced inside the solid domain, whose value depends on the diffusivity coefficient (the penalization parameter) that needs to be decided by the user. Applications of volume penalization include flapping wings~\cite{kolomenskiy2009fourier}, two-phase flow~\cite{horgue2014penalization}, FSI~\cite{engels2015numerical} and thermal flows~\cite{cui2018coupled}. Recently, volume penalization has been applied to high-order flux reconstruction to evaluate the accuracy of IBM in high-order numerical schemes \citep{kou2021IBMFR1,kou2021IBMFR2}. This type of boundary treatment is the focus of this study. 

There have been several attempts to investigate the the numerical properties of IBM. For classical IBM where the singular force is described by delta functions, linear stability analysis is performed to study the stability and stiffness of IBM coupled with time-stepping methods \citep{tu1992stability,stockie1999analysis,gong2008stability}. It was found that the explicit time stepping scheme has the worst stability and can be unstable even though the underlying physical system is stable due to the stiffness of the IBM source term. Spectral analysis for Laplace and Stokes operators based on volume penalization and FD discretization has been studied in \citep{kolomenskiy2014approximation,kolomenskiy2015analysis} for either Dirichlet or Neumann boundary conditions, where it is suggested that if one wants to reach higher accuracy, the underlying discretization should behave well in the presence of discontinuous solutions. Linear stability analysis for the second-order Adams–Bashforth time integration of the penalization term in the incompressible Navier–Stokes equations was performed in \citep{kolomenskiy2009fourier}, where the stability condition $\Delta t < \IBMparam$ was proposed. The eigenvalue problem for elliptic problems based on immersed finite element method is studied in \citep{lee2017immersed}. Recently, eigenvalue and error analysis of the immersed boundary method based on direct forcing is studied by Zhou and Balachandar \citep{zhou2021analysis}. All of these efforts provide valuable insights to deepen the understanding of IBM, while some remaining issues still exist. For example, the IBM treatment based on high-order discretization has not been discussed, as well as some guidelines to improve the implementation in these settings. 

In addition to commonly used low-order methods based on FD and FV, there is an emerging interest in developing high-order methods for computational fluid dynamics (CFD) due to their potential in providing higher accuracy with relatively low cost compared to low-order methods \citep{wang2013high}. Currently, different high-order methods have been developed, including DG \citep{hesthaven2007nodal,karniadakis2013spectral}, spectral difference (SD) \citep{kopriva1996conservative,liu2006SD}, flux reconstruction (FR) \citep{huynh2007FR}, and correction procedure via reconstruction (CPR) \citep{wang2009unifying}. Despite vast efforts, applying high-order methods on unstructured grids remains a challenge due to the difficulty of mesh generation for complex geometries. This makes the development of IBM under high-order frameworks a very appealing alternative. As a representative high-order method, FR is chosen in the present study, which unifies nodal DG and SD schemes in certain conditions \citep{vincent2011new} (e.g. pure advection, as considered in this work). Therefore our analysis covers both FR and DG formulations. This scheme will be combined with volume penalization to simulate the advection equation with a solid wall in the middle of computational domain, thus providing insights into the numerical behavior of IBM on high-order schemes.

The remainder of this paper is organized as follows. In Section 2, the FR method for space discretization and volume penalization method are introduced. Next, Section 3 introduces the eigensolution and non-modal analyses, as well as the novel analysis approach considering the IBM treatment. Additionally, both semi-discrete and fully-discrete analysis are considered to study the behaviors of various factors, including the effects of polynomial order and selection of penalization parameter. Furthermore, numerical experiments are performed in Section 4 to validate and better understand the results from the analysis. A novel approach to include artificial viscosity (second-order term) inside the solid is also proposed and evaluated. Finally, conclusions are drawn in Section 5.

\section{Methodology}
\label{sec1}
\subsection{The governing equation}
We consider a one-dimensional advection equation defined in space $x$ and time $t$: 
\begin{equation}
    \frac{\partial \velx}{\partial t} + \frac{\partial f}{\partial x} = 0,
\end{equation}
where $\velx(x,t)$ is the transported unknown solution and $f = \advcoef \velx$ is the flux, with advection speed $\advcoef$. To apply eigensolution analysis to this equation, we first obtain the analytical solution for a harmonic wave under the periodic boundary conditions and the initial condition $u(x,0) = \text{exp}(ikx)$: 
\begin{equation}
    u(x,t) = \text{exp} [i(kx -\omega t)],
\end{equation}
where $k$ is the wavenumber, $\omega = \advcoef k$ is the (angular) frequency and $i = \sqrt{-1}$ is the imaginary unit. In the present work, this equation is discretized in space based on the high-order flux reconstruction method \citep{huynh2007FR,vincent2011new} (which is equivalent to DG). Details in space and time discretization of the advection equation are given in Appendix \ref{FRdist}. In classic eigensolution analyses, one only needs to look at one element rather than the global domain since it represents the global behavior of the system under periodic or inflow-outflow boundary conditions. However, when IBM is considered, the global matrices for all elements will be taken into account, since some of the solution points should be penalized by additional source terms to impose the IBM conditions. This will be detailed in the following sections.

\subsection{Immersed boundary method based on volume penalization}
Volume penalization imposes boundary conditions by introducing source terms inside the solid region, which is assumed to be a porous medium whose permeability tends to zero. As shown in previous studies \citep{angot1999penalization,kolomenskiy2009fourier,kadoch2012volume}, this method is easy to implement with rigorous theoretical foundations. A mask function $\mask(x,t)$ that distinguishes between the fluid region $\Omega _{f}$ and solid region $\Omega _{s}$ is defined first:
\begin{equation}
\mask  (x,t) = \left\{\begin{matrix}1, \, \, \text{if}\, \, x \in \Omega _{s}
\\0,\, \, \text{otherwise}
\end{matrix}\right. .
\end{equation}

For the one-dimensional advection equation with penalized Dirichlet boundary conditions, this problem is formulated as follows:
\begin{equation}
    \frac{\partial \velx}{\partial t} + \advcoef \frac{\partial \velx}{\partial x} + \frac{\mask}{\IBMparam}(\velx - \velx_s) = 0,
\end{equation}
where $\velx_s$ refers to the penalized solution in the solid, which is fixed to $0$ in this study to approximate the no-slip wall boundary condition. The additional source term is usually a pointwise operation imposed to each solution point. The extension to Neumann boundary conditions is discussed in \citep{kadoch2012volume,sakurai2019volume}. The penalization parameter is defined as $\IBMparam$, which drives the solution to $\velx_s$ as $\IBMparam \rightarrow 0$. For solution inside the solid, this results in the following governing equation similar with convection-diffusion-reaction equations \citep{sengupta2020global}
\begin{equation}
    \frac{\partial \velx}{\partial t} + \advcoef \frac{\partial \velx}{\partial x} + \frac{1}{\IBMparam}\velx = 0.
\end{equation}

The analytical solution of this equation is
\begin{equation}
    \velx(x,t) = \text{exp} [i(kx - (\advcoef k - \frac{i}{\IBMparam}) t)],
\label{eq:ana_solid}
\end{equation}
indicating that the penalization term introduces additional damping to the solution, controlled by the penalization parameter $\IBMparam$. This indicates that volume penalization mimics the porous media with a very low permeability which dissipate the wave-like solution entering the solid and drive it to the expected solution $\velx_s$. From the following eigensolution analysis, this penalization term will have an impact on both the primary and secondary modes by proving additional dissipation.

There are a few works that have discussed the convergence and error estimation of volume penalization, where rigorous proofs have been given~\citep{angot1999penalization,carbou2003boundary}. From Angot et al.~\cite{angot1999penalization} and Carbou and Fabrie~\cite{carbou2003boundary}, it was proven that, as the penalization parameter $\eta$ approaches $0$, the solution of the penalized Navier-Stokes equations will converge to the solution of the Navier-Stokes equations with no-slip boundary conditions. This is one of the advantages of volume penalization over other IBM approaches since the numerical error introduced from the penalization term can be controlled a-priori \cite{brown2014CBVP}. The total error of the solution of penalized equation compared with the body-fitted simulation (non-penalized equation) includes two parts \cite{engels2015numerical}

\begin{equation}
\left \|  \velx^{exact}  - \velx_{\eta} ^ {N} \right \| \leq \left \|  \velx^{exact}  - \velx_{\eta}^{exact} \right \| + \left \|  \velx_{\eta}^{exact} - \velx_{\eta} ^ {N} \right \|,
\end{equation}
where the first part is the penalization error and the second part is the discretization error. $u^{exact} $, $u_{\eta} $ and $u_{\eta} ^ {N}$ are the exact analytical solution of the governing equations, the exact and numerical solution of the penalized equations, respectively. The error is quantified by the norm $\left \|  \cdot \right \|$ between solutions. The penalization error depends on the penalization parameter \cite{schneider2015immersed}:
\begin{equation}
\left \|  \velx^{exact}  - \velx_{\eta}^{exact} \right \| \propto \eta ^ {\alpha},
\end{equation}

This indicates that as the penalization parameter approaches zero, the error between exact solutions of penalized and original equation will converge to zero, i.e., $\lim_{\eta \rightarrow 0} \left \|  u^{exact}  - u_{\eta}^{exact} \right \| \rightarrow 0$. It has been proved that the volume penalization gives $\alpha = \frac{1}{2}$ for the Dirichlet boundary conditions, indicating the penalization error has a decay rate of $\mathcal{O}(\sqrt{\eta})$. For Neumann boundary condition, $\mathcal{O}(\eta)$ can be obtained \cite{kadoch2012volume,kolomenskiy2015analysis}. The discretization error refers to the error between the exact solution and the numerical solution of the penalized equations. With consistent discretization and a stable numerical scheme, the discretization error usually follows ($\beta > 0$):
\begin{equation}
	 \left \|  u_{\eta}^{exact} - u_{\eta} ^ {N} \right \| \propto N ^ {-\beta}.
\end{equation}

However, as pointed out by Schneider et al. \cite{kadoch2012volume,schneider2015immersed}, the convergence for the discretization error is not only determined by the numerical scheme, but also limited by the regularity of the solution, which refers to the continuity (smoothness) of the exact solution $u_{\eta}$ at the boundary of the penalized equation. This regularity can be improved by designing proper $\velx_s$ inside the solid \citep{gautier2014dns,stein2016immersed}. However, in this study we will discuss the classical volume penalization method with $\velx_s = 0$ for a no-slip wall. This analysis also indicates the difficulty of IBM to achieve high order convergence near the boundary, with any type of space discretization. From these theories, it is suggested to use a penalization parameter $\eta$ which is small enough to ensure low penalization error (high penalization since $\eta$ is in the denominator). In the meantime, $\eta$ should also be limited since the penalization source term can become very stiff, leading to numerical stabilities. When explicit time integration is considered, it is suggested to use $ \Delta t < \eta $ \cite{kolomenskiy2009fourier} or $ \Delta t \approx \eta $ \cite{jause2012simulation} to ensure stability. An interpretation for this choice is that the penalty term acts as a strong damping term with order $\eta$ on the velocity, which has to be resolved by the time discretization scheme \cite{engels2015numerical}. We will provide further insights into the selection of the penalty parameter in what follows.

\section{Modal and non-modal analysis}
\subsection{Eigensolution analysis}
In the present temporal  eigensolution analysis, we assume periodic boundary condition and consider solution $u(x,t) = \text{exp} [i(kx -\omega t)]$ with real wavenumber $k$. Therefore, for each element with index $n$, we have $\boldsymbol{u}_{n-1} = e^{-ikh}\boldsymbol{u}_{n}$ and $\boldsymbol{u}_{n+1} = e^{+ikh}\boldsymbol{u}_{n}$. The space discretization results in 
\begin{equation}
    \frac{d\boldsymbol{u}_n}{dt} = \left ( \boldsymbol{L} e^{-ikh} +  \boldsymbol{C}  +  \boldsymbol{R} e^{+ikh} \right ) \boldsymbol{u}_{n} = \boldsymbol{M} \boldsymbol{u}_{n}.
\end{equation}

This semi-discrete matrix formulation is the basis for the eigensolution analysis and the non-modal analysis. It can be transformed into the compact form of an eigenvalue problem
\begin{equation}
    - i \omega^* \boldsymbol{u}_n = \boldsymbol{M} \boldsymbol{u}_{n},
\end{equation}
where $\omega^* = \advcoef k^*$ becomes a complex value due to the dispersion and dissipation errors of the space discretization. This wavenumber $\omega^*$ relates to the eigendecomposition of the coefficient matrix $\boldsymbol{M}$ (with $P+1$ solutions)
\begin{equation}
\label{eq:LCR}
    - i \omega^*_m = \lambda_m \\, \ \lambda_m \boldsymbol{v}_m = \boldsymbol{M} \boldsymbol{v}_m,
\end{equation}
where $\lambda_m$ and $\boldsymbol{v}_m$ are the $m$th eigenvalues and eigenvectors of matrix $\boldsymbol{M}$, respectively. By defining the modified wavenumber $k^*_m$ for each $\omega^*_m$ and $\lambda_m$, we have:
\begin{equation}
    \text{Real}(k^*_m) = \frac{\text{Real}(\omega^*_m)}{\advcoef} = -\frac{\text{Imag}(\lambda_m)}{\advcoef}\\, \ \text{Imag}(k^*_m) = \frac{\text{Imag}(\omega^*_m)}{\advcoef} = \frac{\text{Real}(\lambda_m)}{\advcoef}.
\end{equation}

For the advection equation, the difference between $\text{Real}(k^*)$ and $k$ is due to the dispersion error, indicating the change in wavenumber of the solution. The difference between $\text{Imag}(k^*)$ and $0$ corresponds to the dissipation (diffusion) error, where $\text{Imag}(k^*) \leq 0$ holds for a stable scheme. Note that for the analytical solution we have $\text{Real}(k^*) = k$ and $\text{Imag}(k^*) = 0$. For high-order methods, the wavenumber is normalized by the $h/(P+1)$ which represents the smallest length scale that can be captured by the scheme. In addition, $\boldsymbol{M}$ will have $P+1$ eigenmodes, where the one that recovers $k^* = k$ as $k \rightarrow 0$ will be identified as the physical mode, while others are secondary modes which are indeed translated replicas of the primary mode \citep{moura2015linear}. 

When the IBM is included in the analysis, one has to consider the whole computational domain rather than one element, while penalizing the solution points inside the solid domain. To account for multiple (equispaced) elements, we consider the direct sum of $N$-element contributions as shown in \citep{manzanero2018dispersion}, with solution vector defined as $\boldsymbol{U} = [\boldsymbol{u}_1,\boldsymbol{u}_2,...,\boldsymbol{u}_N]^T$. We focus on the computational domain $[-T, T]$ with periodic boundary condition for the advection equation. In this case, the effect of the physical boundary conditions need to be introduced through two ghost elements for external states, $\boldsymbol{u}_0$ and $\boldsymbol{u}_{N+1}$, derived as $\boldsymbol{u}_0 = \text{exp}[-2ikT] \boldsymbol{u}_N$ and $\boldsymbol{u}_{N+1} = \text{exp}[2ikT] \boldsymbol{u}_1$. For cases without IBM, this represents the continuation of a hypothetical infinite mesh. Since volume penalization imposes pointwise source terms, the penalization term is only added to the diagonal of the global matrix. Following \citep{manzanero2018dispersion}, we can reach the global matrix form of semi-discrete formulation considering IBM as follows:
\begin{equation}
    \frac{d}{dt} \begin{pmatrix} \boldsymbol{\velx}_{1}\\  \boldsymbol{\velx}_{2}\\ \boldsymbol{\velx}_{3}\\ \vdots\\ \boldsymbol{\velx}_{N} \end{pmatrix} = \begin{pmatrix} \boldsymbol{C}-\frac{\chi}{\eta}\boldsymbol{I}
    & \boldsymbol{R} & \boldsymbol{0} & \cdots & \text{exp}[-2ikT] \boldsymbol{L}\\ \boldsymbol{L} & \boldsymbol {C}-\frac{\chi}{\eta} \boldsymbol{I} & \boldsymbol{R} & \cdots & \boldsymbol{0} &\\ \boldsymbol{0} & \boldsymbol{L} &\boldsymbol{C}-\frac{\chi}{\eta}\boldsymbol{I} & \cdots & \boldsymbol{0}\\ \vdots & \vdots & \vdots & \ddots &\vdots\\ \text{exp}[2ikT] \boldsymbol{R} & \boldsymbol{0} & \boldsymbol{0} & \cdots &\boldsymbol{C}-\frac{\chi}{\eta}\boldsymbol{I} \end{pmatrix}
    \begin{pmatrix} \boldsymbol{\velx}_{1}\\  \boldsymbol{\velx}_{2}\\ \boldsymbol{\velx}_{3}\\ \vdots\\ \boldsymbol{\velx}_{N} \end{pmatrix},
\end{equation}
where $I$ denotes the identity matrix with dimension $P+1$. The volume penalization term is imposed when the corresponding solution point lies inside the solid ($\mask = 1$). For example, assuming we have one solid element in the middle, the global matrix becomes
\begin{equation}
\label{eq:IBMglobal}
    \frac{d}{dt} \begin{pmatrix} \boldsymbol{\velx}_{1}\\  \boldsymbol{\velx}_{2} \\ \vdots\\ \boldsymbol{\velx}_{IBM} \\ \vdots \\ \boldsymbol{\velx}_{N} \end{pmatrix} = \begin{pmatrix} \boldsymbol{C}
    & \boldsymbol{R} & \boldsymbol{0} & \cdots & \cdots & \cdots & \text{exp}[-2ikT] \boldsymbol{L}\\ \boldsymbol{L} & \boldsymbol{C} &\boldsymbol{R} & \cdots & \cdots & \cdots & \boldsymbol{0} &\\ \vdots & \vdots & \vdots & \ddots & \ddots & \ddots & \vdots\\ \vdots & \vdots & \vdots & \boldsymbol{L} & \boldsymbol{C}-\frac{1}{\eta}\boldsymbol{I} & \boldsymbol{R}& \vdots \\ \vdots & \vdots & \vdots & \ddots & \ddots & \ddots & \vdots \\
    \text{exp}[2ikT] \boldsymbol{R} & \boldsymbol{0} & \boldsymbol{0} & \cdots & \cdots & \cdots &\boldsymbol{C} \end{pmatrix}
    \begin{pmatrix} \boldsymbol{\velx}_{1}\\  \boldsymbol{\velx}_{2} \\ \vdots\\ \boldsymbol{\velx}_{IBM} \\ \vdots \\ \boldsymbol{\velx}_{N} \end{pmatrix}.
\end{equation}

This is equivalent to an advection problem with a solid wall in the middle of a periodic computational domain. By investigating this global matrix, the dispersion-dissipation behavior of IBM for high-order schemes can be analyzed. Note that this matrix will contain $N (P+1)$ eigenmodes, but as pointed out in \citep{manzanero2018dispersion}, only one of the modes tracks the physical propagation speed and damping, which is referred to as the primary mode or physical mode (again, the one that recovers $k^* = k$ as $k \rightarrow 0$ \citep{moura2015linear}). All remaining modes are referred to as secondary modes.

\subsection{Non-modal analysis}
Non-modal analysis \citep{fernandez2019non} is an approach to analyze the numerical diffusion, which is more relevant in under-resolved turbulence simulations. A short-term diffusion is defined that describes the change of magnitude of the numerical solution evolving right after time $t=0$. It differs from classic eigensolution analysis since the contribution of all eigenmodes is taken into account, which leads to bigger difference in the high wavenumber range. The short-term diffusion parameter $\Tilde{\omega}^*$ is defined as \citep{fernandez2019non}
\begin{equation}
    \Tilde{\omega}^* :=\left(\frac{d\log ||\boldsymbol{\velx}_n||}{d\tau^*} \right)_{\tau^*=0},
\end{equation}
where $|| \cdot ||$ is the $L_2$ norm and $\tau ^ * = \tau (P+1) = t \advcoef (P+1)/h $ is a non-dimensional time based on the non-dimensional convection time per degree of freedom. For $\Tilde{\omega}^* \leq 0$ the scheme is stable. Considering the orthogonality of solution polynomials and the wave-like behavior of the numerical solution, the above equation can be rewritten as the form of a Rayleigh quotient:
\begin{equation}\label{stermdiff}
\Tilde{\omega}^* = \frac{1}{P+1}\text{Real} \left [ {\frac{ (\boldsymbol{\velx}_{n,0})^\dagger\, \boldsymbol{M}\,\boldsymbol{\velx}_{n,0} }{ (\boldsymbol{\velx}_{n,0})^\dagger \boldsymbol{\velx}_{n,0}}} \right ],
\end{equation}
where $\boldsymbol{\velx}_{n,0}$ is the initial condition, and $\dagger$ denotes the complex conjugate. $\boldsymbol{M}$ is the discretization operator in Equation \ref{eq:LCR} and Equation \ref{eq:IBMglobal}. The short-term diffusion parameter can be understood as the combined contribution of all eigenmodes of the discretization operator identified in a von Neumann analysis \citep{fernandez2019non}. Therefore, non-modal analysis is consistent with the von Neumann analysis whenever only one eigenmode exists (at $P=0$) or when the initial solution is an eigenvector of the discretization operator $\boldsymbol{M}$. 

\subsection{Semi-discrete analysis}
In this section, a semi-discrete eigensolution analysis of the FR scheme with volume penalization for the advection equation is performed. Since semi-discrete analysis only looks at the space discretization, it performs eigensolution analysis and non-modal analysis on the space discretization operator matrix. The analysis of the standard advection equation focuses on Equation \ref{eq:LCR}, while the analysis of IBM focuses on Equation \ref{eq:IBMglobal}. In the following analysis, we consider the advection equation with $\advcoef = 1$ and a computational domain defined in $x \in [-1,1]$ discretized by equispaced elements. All the cases considered here are based on the fully upwind flux with $\lambda = 1$. The solid region lies in the middle, starting from $x=0$, whose width is defined as $\Delta$:
\begin{equation}
\mask  (x,t) = \left\{\begin{matrix}1, \, \, \text{if}\, \, 0 < x < \Delta 
\\0,\, \, \text{otherwise}
\end{matrix}\right. .
\end{equation}

In the present study, we will consider $\Delta$ as integer multiples of element size $h$ ($\Delta = Zh$, where $Z$ is an integer), thus the solid boundaries will appear exactly at the interface between elements. This allows us to define the solid ratio $r = Z/N$ as the ratio between the solid region and the computational domain. For a wavelike initial condition, the solution is expected to approach $0$ when $x>0$ and $x < t$ (consider a short period). As time evolves, the global solution will eventually become $0$ in the whole computation domain, due to the solid and the periodic boundary condition. Moreover, when the solid wall exists, it is expected to have $\advcoef = 0$ for all solution points inside the wall. This indicates that the actual fluid domain is shorter than the whole computational domain, therefore the exact wavenumber in the fluid region is higher than the initial wavenumber. To obtain the real wavelength in the fluid region, the initial wavenumber $k$ should be re-scaled by the solid ratio $r$
\begin{equation}
\label{eq:rescaling}
    \hat{k} = k / (1-r).
\end{equation}

\begin{figure*}[htbp]
		\centering
		\includegraphics[width=300pt]{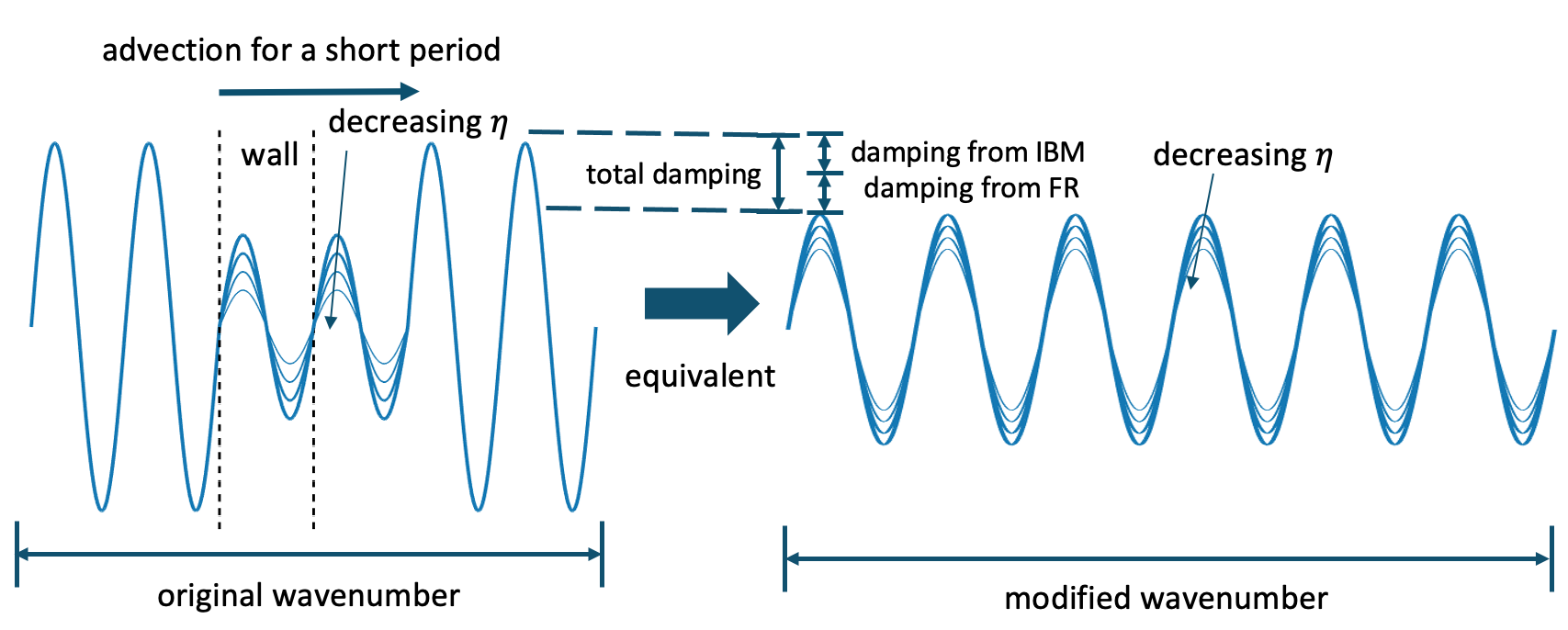}
		\caption{Schematic illustration of the advection problem with IBM.}
		\label{fig:cartoon}
\end{figure*}

This definition is consistent with the standard advection equation where we have $r = 0$ thus $\hat{k} = k$. Note that this re-scaling only applies to the initial wavenumber $k$ to match the modified wavenumber obtained from the eigensolution analysis $k^*$. It should also be noted that the inclusion of IBM sometimes requires shifting the wavenumber by $\pi / N$ to accurately trace the physical mode in eigensolution analysis. In summary, a schematic illustration of the considered problem is shown in Figure \ref{fig:cartoon}, where we consider an initial wave for a short period convection. The IBM wall is introduced in the middle, therefore the initial wave in the middle starts to be damped and will move towards right as the solution evolves. This problem is different from the assumption of classical eigensolution analysis where the solution is decomposed into the combination of harmonics. Therefore to better understand the present problem settings, we can think of eigensolution analysis as transforming the present problem into an equivalent problem in the right of Figure \ref{fig:cartoon}, where the solution can be represented as a harmonic wave with modified wavelength $\hat{k}$ and equivalent amplitude (damping). The wavelength is re-scaled due to the existence of the solid region, as defined in Equation \ref{eq:rescaling}. Since the analysed equation is linear, the equivalent damping is an average dissipation effect coming of both the local volume penalization term (averaged over the whole computational domain) and the FR scheme itself. By looking as this equivalent problem, we have the following questions: 1) how does the penalization parameter influences the dispersion-dissipation behavior? 2) will IBM change the spectral behaviors of the original FR scheme? Both of them will be answered in the following analysis and numerical tests.

\begin{figure*}[htbp]
\begin{subfigure}{.24\textwidth}
		\centering
		\includegraphics[width=110pt]{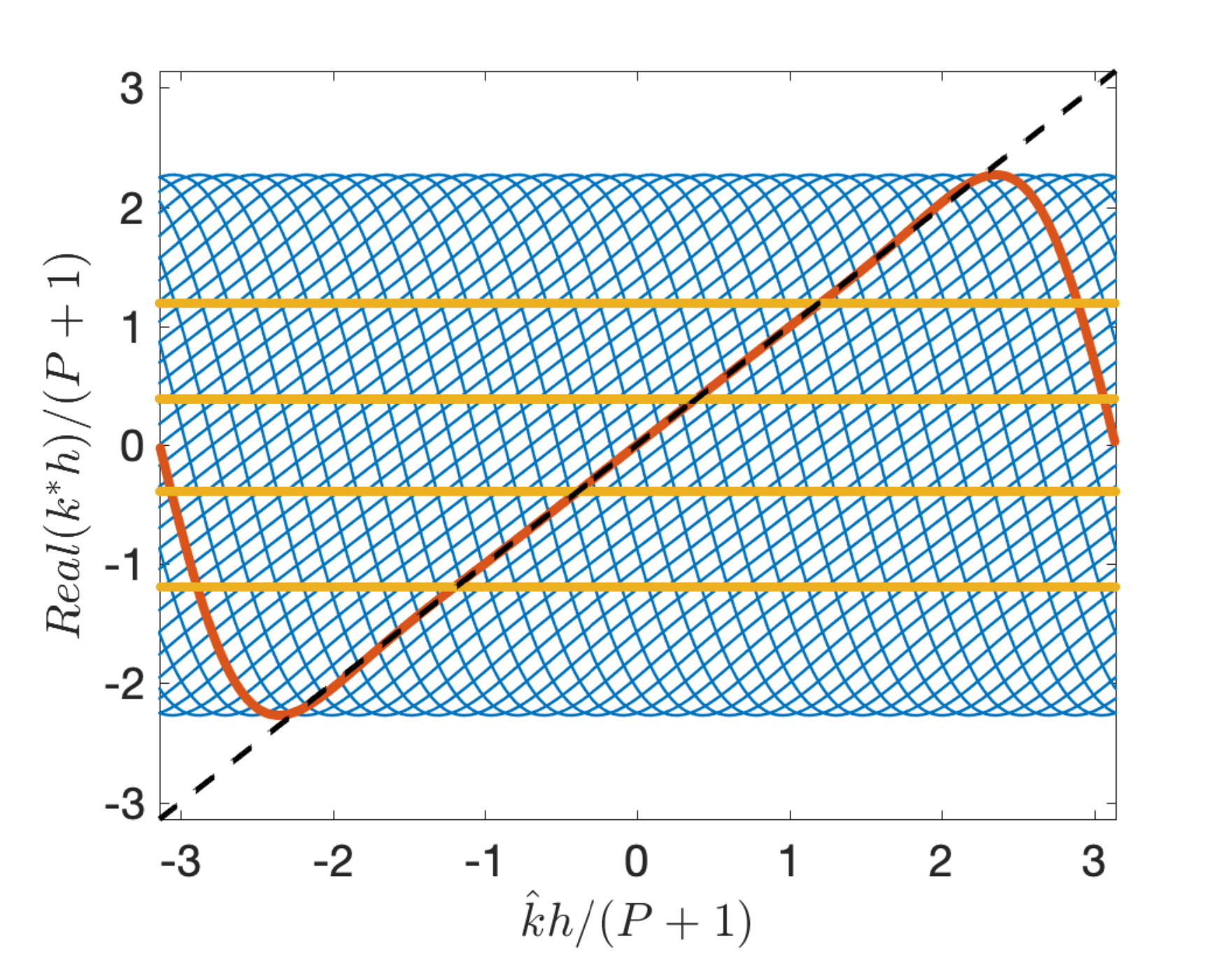}
		\caption{}
	\end{subfigure}
	\begin{subfigure}{.24\textwidth}
		\centering
		\includegraphics[width=110pt]{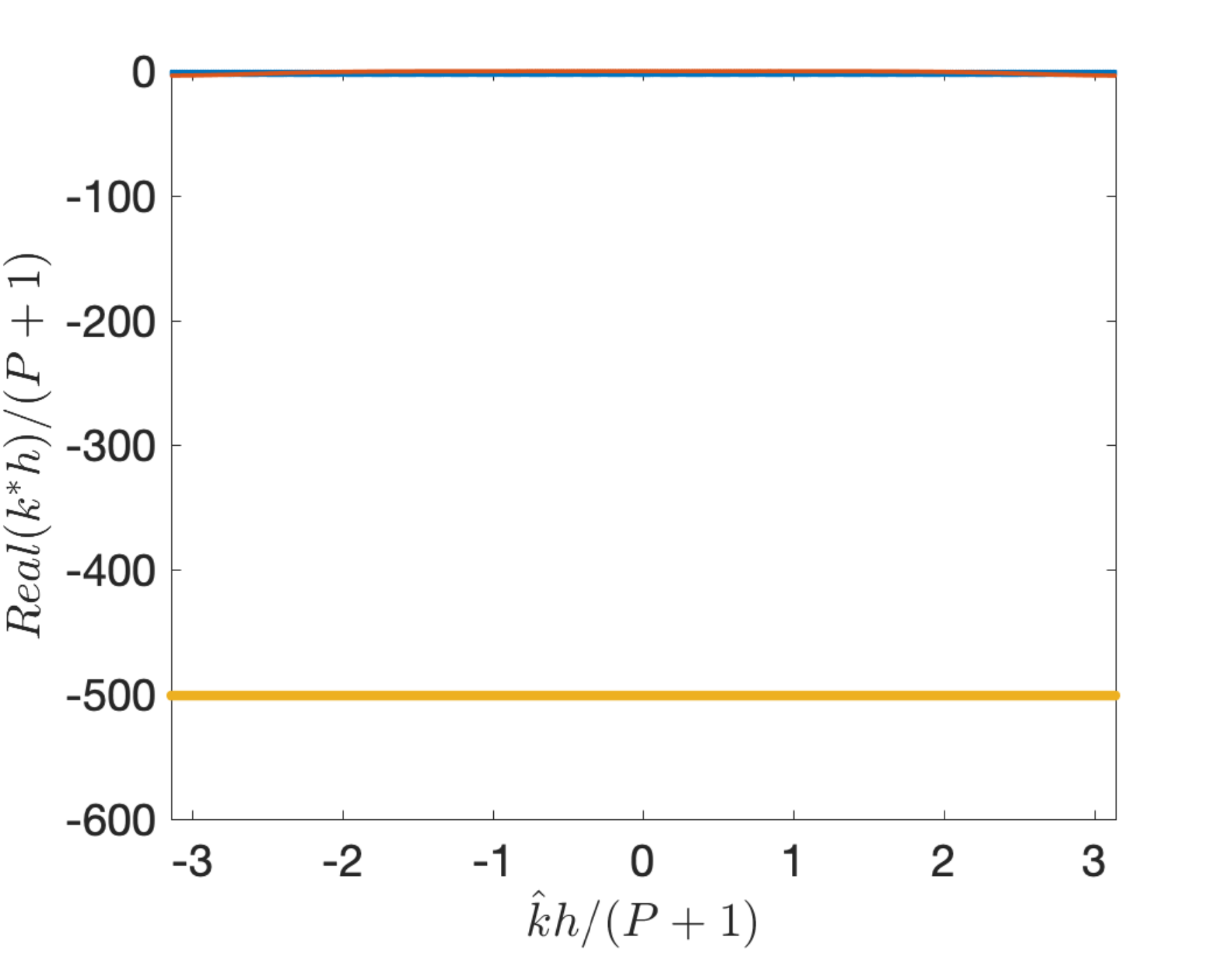}
		\caption{}
	\end{subfigure}
	\begin{subfigure}{.24\textwidth}
		\centering
		\includegraphics[width=110pt]{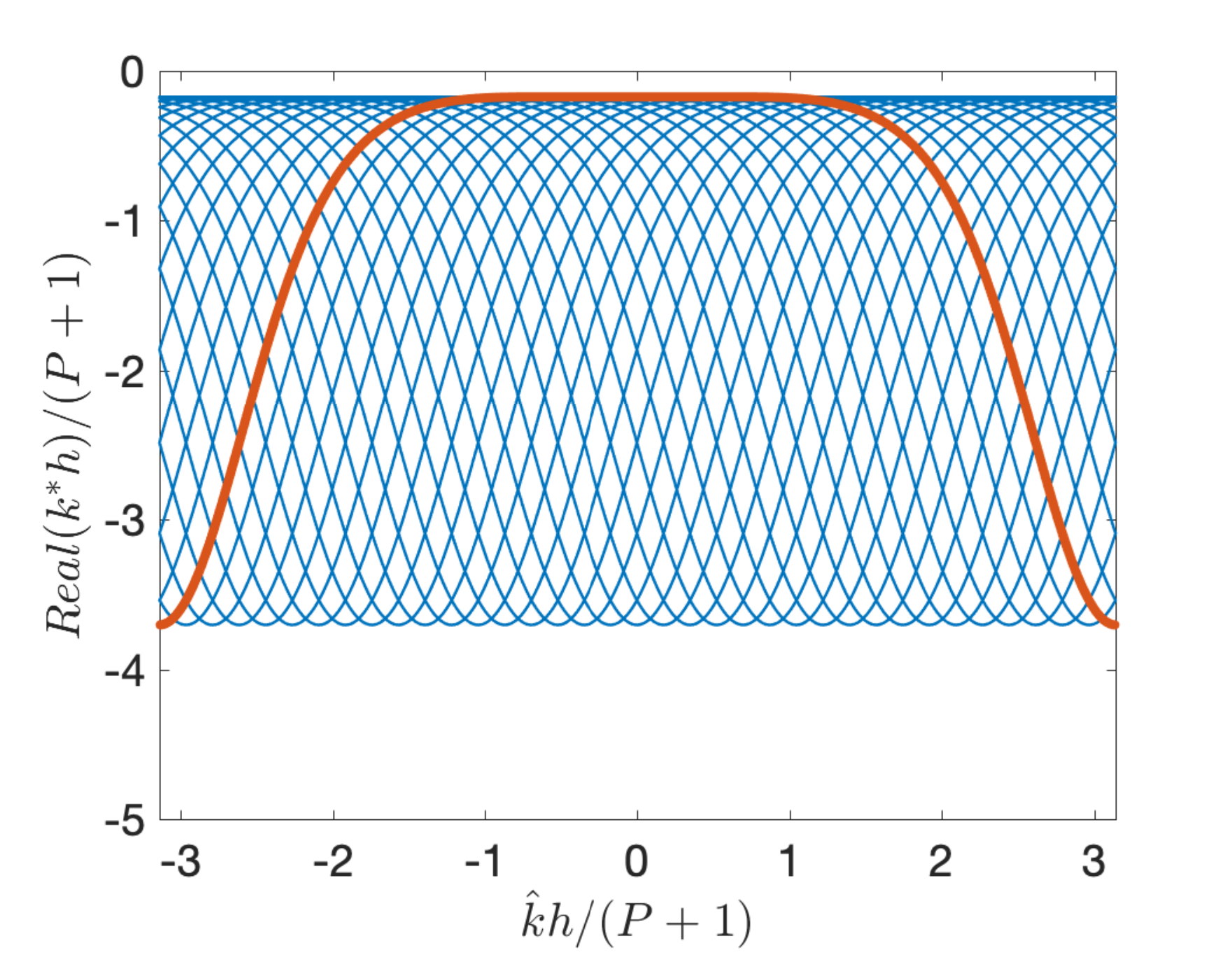}
		\caption{}
	\end{subfigure}
	\begin{subfigure}{.24\textwidth}
		\centering
		\includegraphics[width=110pt]{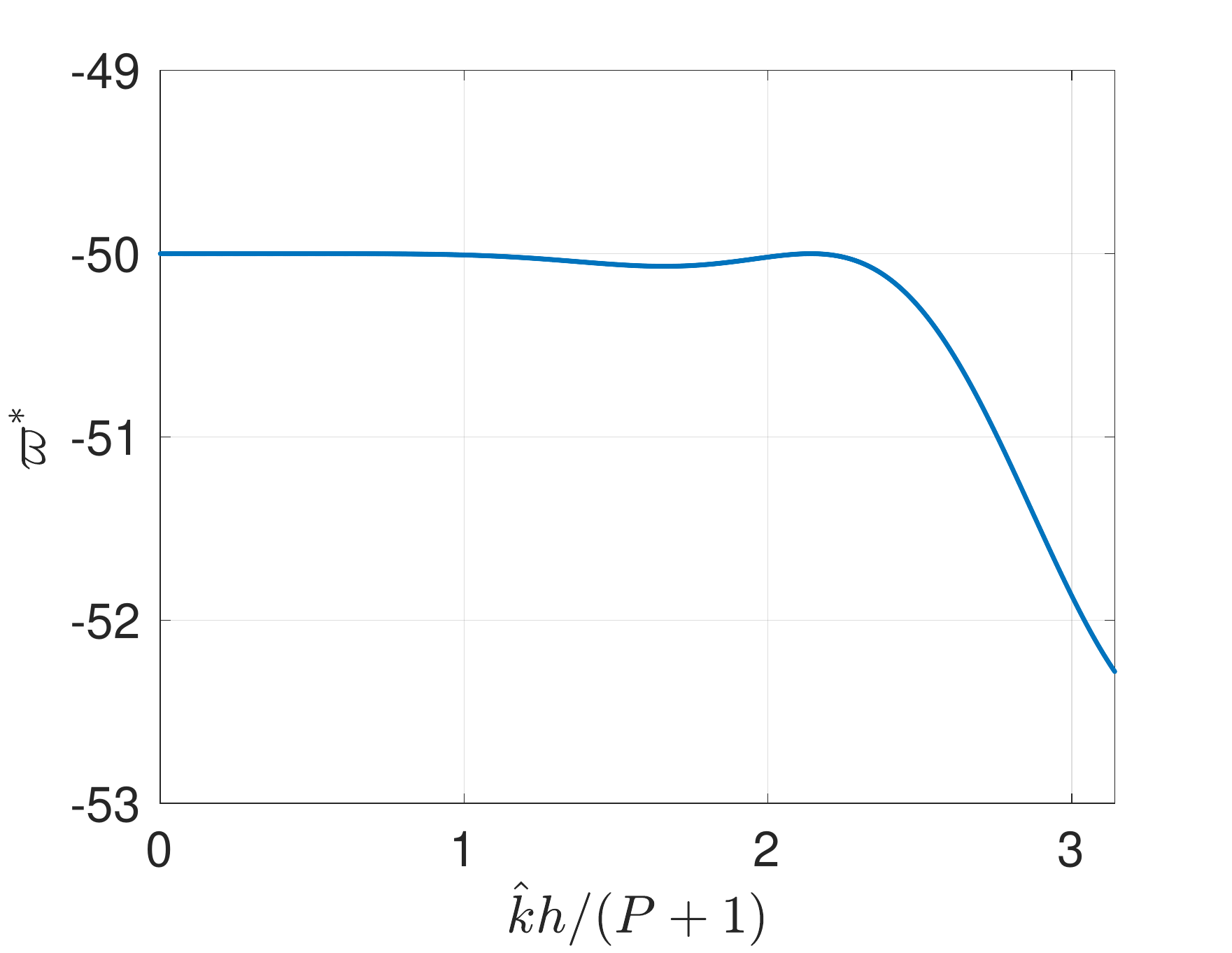}
		\caption{}
	\end{subfigure}
	\centering
	\caption{Eigensolution and non-modal analysis of advection equation with IBM ($N = 10$, $r = 1/10$, $P=3$ and $\IBMparam = 1 \times 10^{-4}$). a) Dispersion. b) Dissipation. c) Zoom-in dissipation. d) Short-term dissipation. Red and yellow curves represent the physical mode and the mode induced by the IBM source term. Deshed line is the refernce curve $Real(k^*) = k$, $Imag(k^*) = 0$ and $\omega^* = 0$.}
	\label{fig:example}
\end{figure*}

As the first example, we consider a coarse discretization, where the computational domain is discretized into $N=10$ elements, with the solid width $\Delta = h = 0.2$ and $r = 1/10$. The polynomial order is set to $P=3$, and the penalization parameter is set to $\IBMparam = 1 \times 10^{-4}$. Results from eigensolution and non-modal analysis are shown in Figure \ref{fig:example}, where all eigenmodes are shown. In general, the dispersion-dissipation behavior agrees well with those of the standard advection equation without IBM, as reported in \citep{vincent2011insights,moura2015linear,manzanero2018dispersion}, where one physical / primary mode is identified while other modes are mostly replications of the primary mode. One distinction of the present case is that the volume penalization term will introduce several solid modes (depending on the number of degrees of freedom inside the solid) with constant dispersion and dissipation across all wavenumbers, as shown in Figure \ref{fig:example}a and Figure \ref{fig:example}b. They are new secondary modes when IBM is considered, which also agrees with the analysis of Equation \ref{eq:ana_solid} that constant dissipation is introduced inside the wall. To understand the spectral behavior, the physical mode is more interesting. It can be seen that IBM retains reasonable dispersion behavior at low and medium wavenumbers, and introduces additional dissipation for the primary mode across the wavenumber range, as shown in Figure \ref{fig:example}c. This is reasonable because when any solution (including the constant solution $k = 0$) passes through the IBM wall, it will be damped in order to satisfy the penalized solution $\velx_s$. Therefore, the dissipation at $k = 0$ refers to the damping only induced by the IBM condition, which is the shift of a constant solution. When we look at the effect of IBM on the FR scheme, this constant dissipation needs to be subtracted from the total dissipation, as illustrated in Figure \ref{fig:cartoon}. Finally, the total dissipation effect for all fluid and solid modes can be seen in the short-term dissipation curve in Figure \ref{fig:example}d, reflecting the combined damping effect in the whole computational domain. It can be concluded that the short-term dissipation is dominated by the IBM source term, resulting in very large negative dissipation depending on $\IBMparam$. This indicates that most of the damping is due to the effect of IBM, resulting very large short-term dissipation at $k=0$. Again, the effect of IBM on the FR scheme should be obtained by subtracting this initial dissipation. The magnitude of this short-term damping will be smaller than that of the single solid mode in the dissipation curve in Figure \ref{fig:example}b, since the short-term damping is an averaged damping in the computational domain, but the solid mode is only a characteristic of the solid body. It should be mentioned that the magnitude of these damping factors rely heavily on $\IBMparam$ which determines the decay of the solution inside the solid. Following this analysis, we can define the IBM-induced dissipation $\gamma_{IBM}$ for the primary mode as
\begin{equation}
    \gamma_{IBM} = Imag[k^*(0)]h / (P+1),
\end{equation}
where $k^*(0)$ refers to the modified wavenumber $k^*$ at $k = 0$. Therefore, the total dissipation $\gamma$ of primary mode is a combination of IBM (constant dissipation across all wave numbers) and the numerical scheme (modified behavior due to the existence of IBM)
\begin{equation}
    \gamma = \gamma_{IBM} + \gamma_{FR}^*,
\end{equation}
which also applies to short-term dissipation based on the non-modal analysis. By comparing $\gamma_{FR}^*$ with the dissipation of FR for standard advection equation, the influence of IBM on the dissipation of the FR scheme can be studied.
In the following investigation, for eigensolution analysis, we will concentrate on the primary mode in the semi-discrete analysis, and the primary and IBM modes in the fully-discrete analysis, since the latter is the main cause of instability. Moreover, to compare the dissipation behavior more reasonably, the IBM-induced dissipation $\gamma_{IBM}$ of the physical mode and the short-term dissipation should be subtracted, in order to only look at the dissipation effect of IBM on FR itself ($\gamma_{FR}^*$). 

\begin{figure*}[htbp]
\begin{subfigure}{.24\textwidth}
		\centering
		\includegraphics[width=110pt]{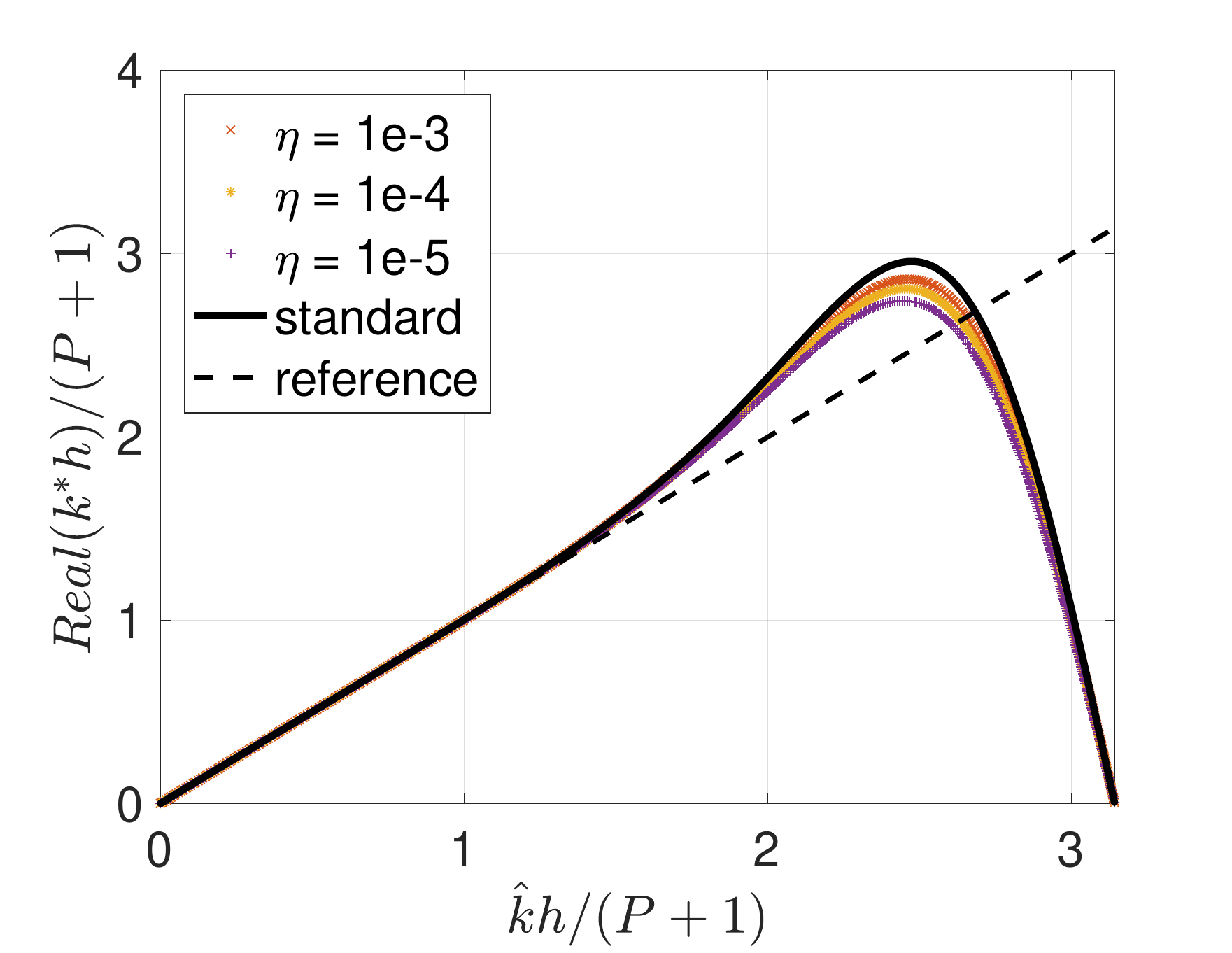}
		\caption{}
	\end{subfigure}
	\begin{subfigure}{.24\textwidth}
		\centering
		\includegraphics[width=110pt]{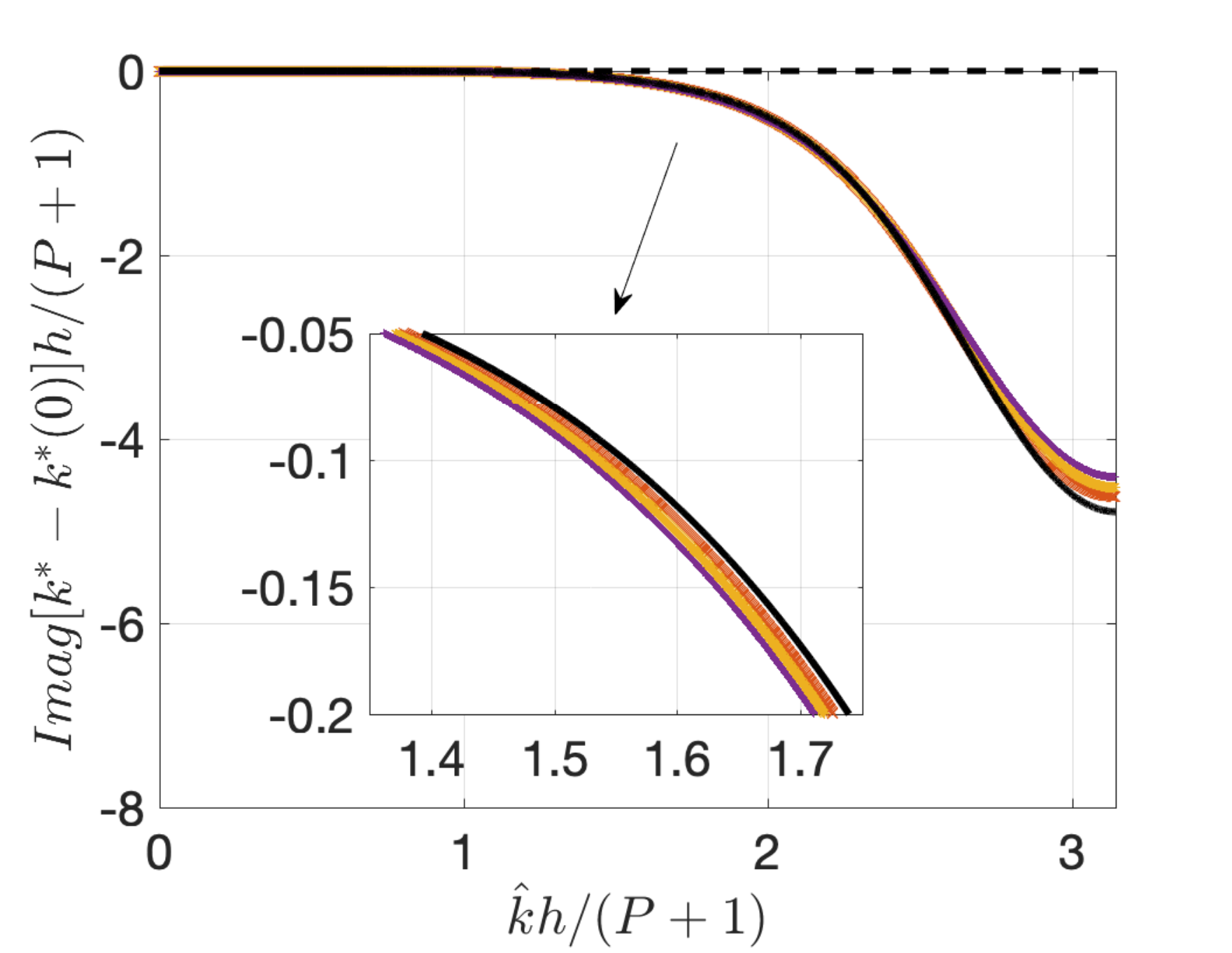}
		\caption{}
	\end{subfigure}
	\begin{subfigure}{.24\textwidth}
		\centering
		\includegraphics[width=110pt]{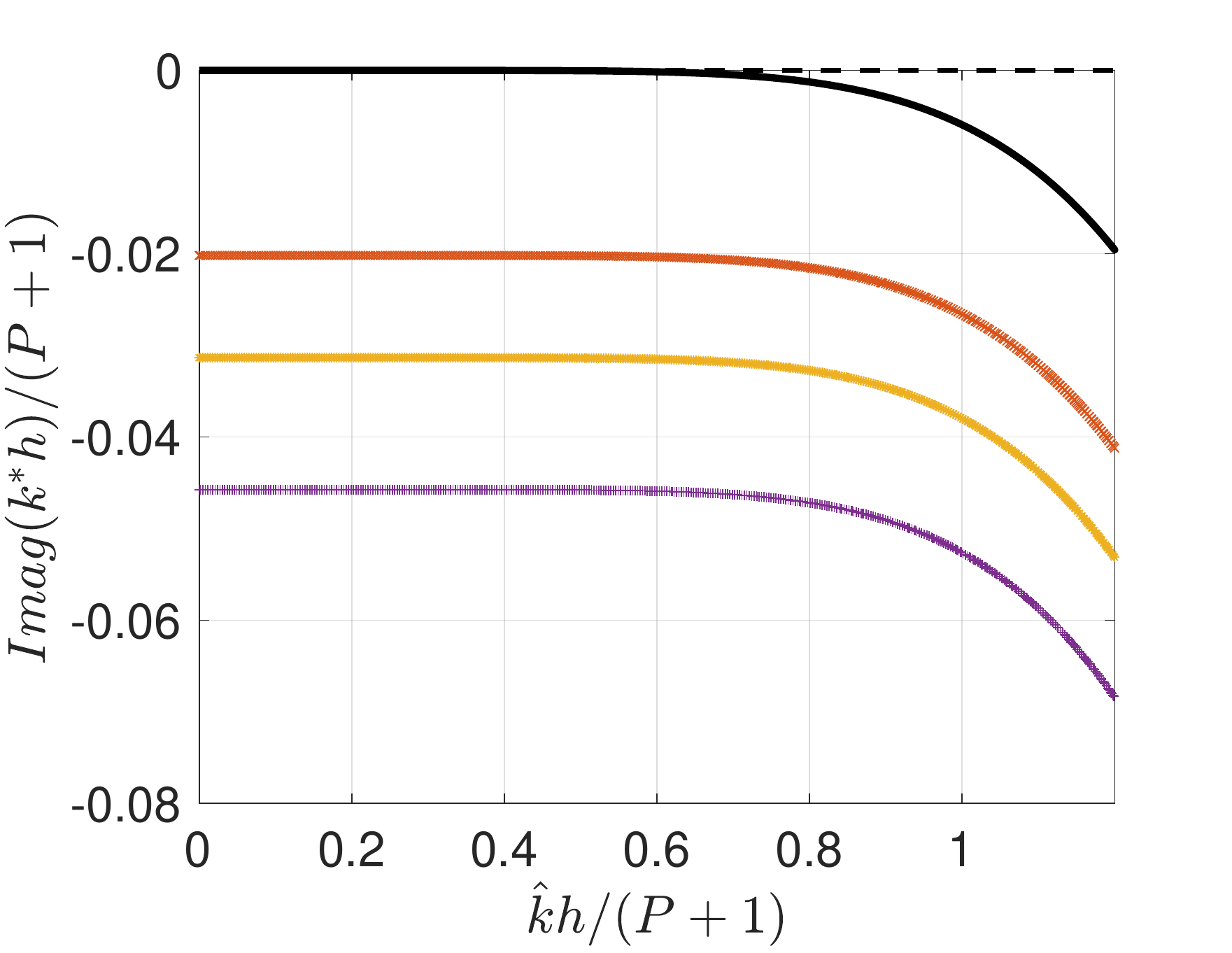}
		\caption{}
	\end{subfigure}
	\begin{subfigure}{.24\textwidth}
		\centering
		\includegraphics[width=110pt]{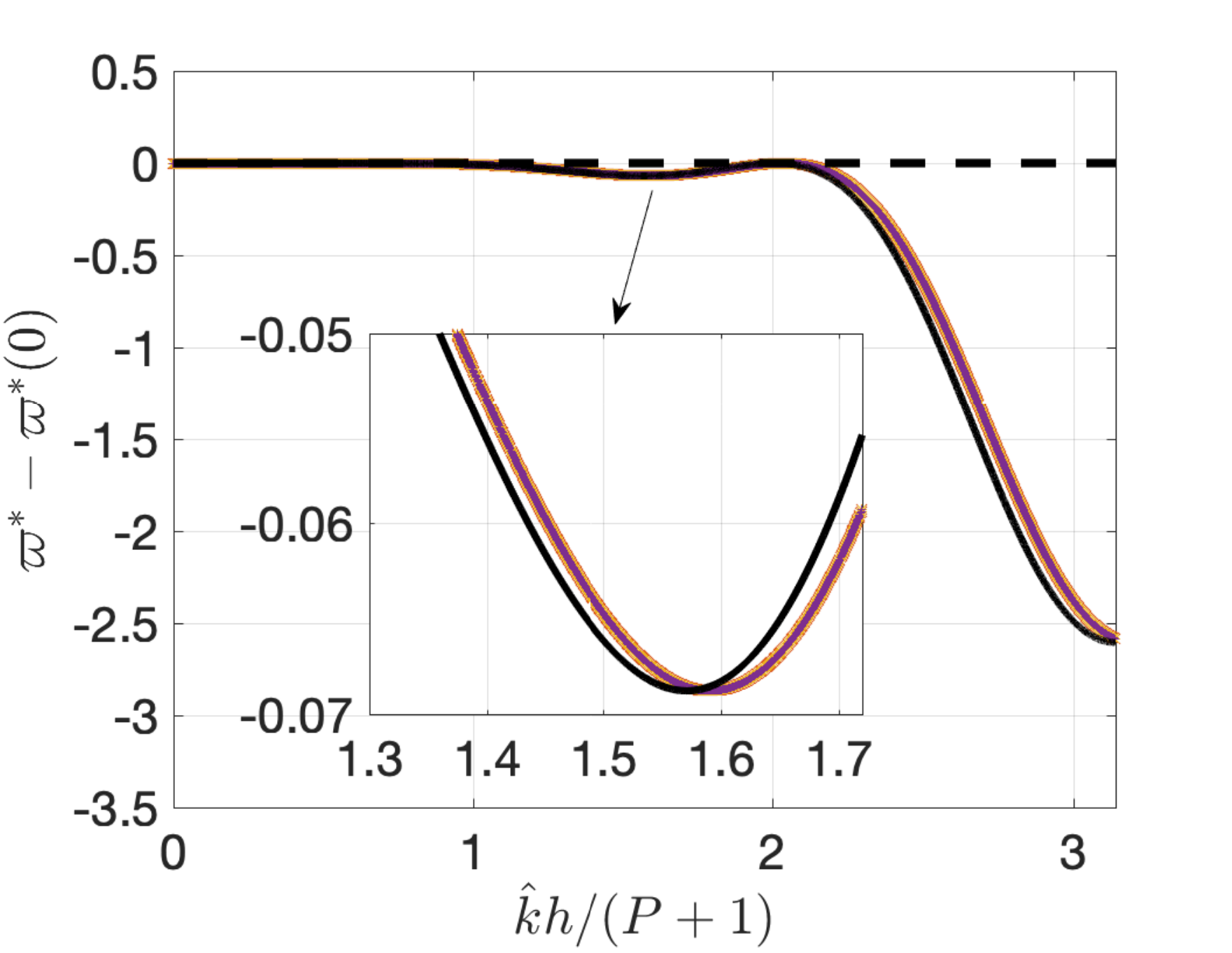}
		\caption{}
	\end{subfigure}
	\caption{Eigensolution and non-modal analysis of advection equation with IBM ($N = 40$, $r = 1/40$, $P=3$) and different penalization parameter $\eta$. a) Dispersion. b) Dissipation subtracted by the IBM-induced dissipation $\gamma_{IBM}$ (dissipation at $k=0$). c) Zoom-in dissipation. d) Short-term dissipation subtracted by the dissipation at $k=0$. The reason for subtracting the initial value is because this value depends on $\eta$ and is only an effect of IBM source term.}
	\label{fig:semi-penalty}
\end{figure*}

The influence of the penalization parameter is studied first. We consider a finer discretization with $N = 40$ elements in the computational domain ($h=0.05$). We set $\Delta = h$ with $r = 1/40$ and $P = 3$. Three penalization parameters $1 \times 10^{-3}$, $1 \times 10^{-4}$ and $1 \times 10^{-5}$ are selected, which are relatively small because we need $\IBMparam$ to be small enough to guarantee good accuracy \citep{angot1999penalization}. The results are compared against the standard analysis obtained for the standard advection equation, which are shown in Figure \ref{fig:semi-penalty}. As explained, the initial dissipation ($\gamma_{IBM}$) is subtracted for comparison purpose, in order to only look at the dissipation effect of the numerical scheme, as shown in Figure \ref{fig:semi-penalty}b. It can be seen that in well resolved regions (small to medium wavenumber), the dispersion-dissipation behaviors similar with the standard advection scheme are recovered, while a slight difference is seen in the high wavenumber region. This observation applies for both dispersion-dissipation and non-modal analyses, as shown in Figure \ref{fig:semi-penalty}a, \ref{fig:semi-penalty}b, and \ref{fig:semi-penalty}d. This is expected because it indicates that IBM does not affect too much the underlying FR scheme in the resolved wavenumber range. Another information is about total dissipation (both the IBM-induced dissipation and the dissipation of numerical scheme) for the physical mode, as shown in Figure \ref{fig:semi-penalty}c. It is clear that as $\eta \rightarrow 0$, more damping is added to the primary mode, indicating that the solution will be further damped and the boundary condition will be imposed more strongly. This analysis can justify the use of volume penalization, since this method affects the spectral behavior of the underlying numerical scheme marginally, but only adds constant dissipation inside the wall to mimic the no-slip wall boundary condition.

\begin{figure*}[htbp]
\begin{subfigure}{.3\textwidth}
		\includegraphics[width=120pt]{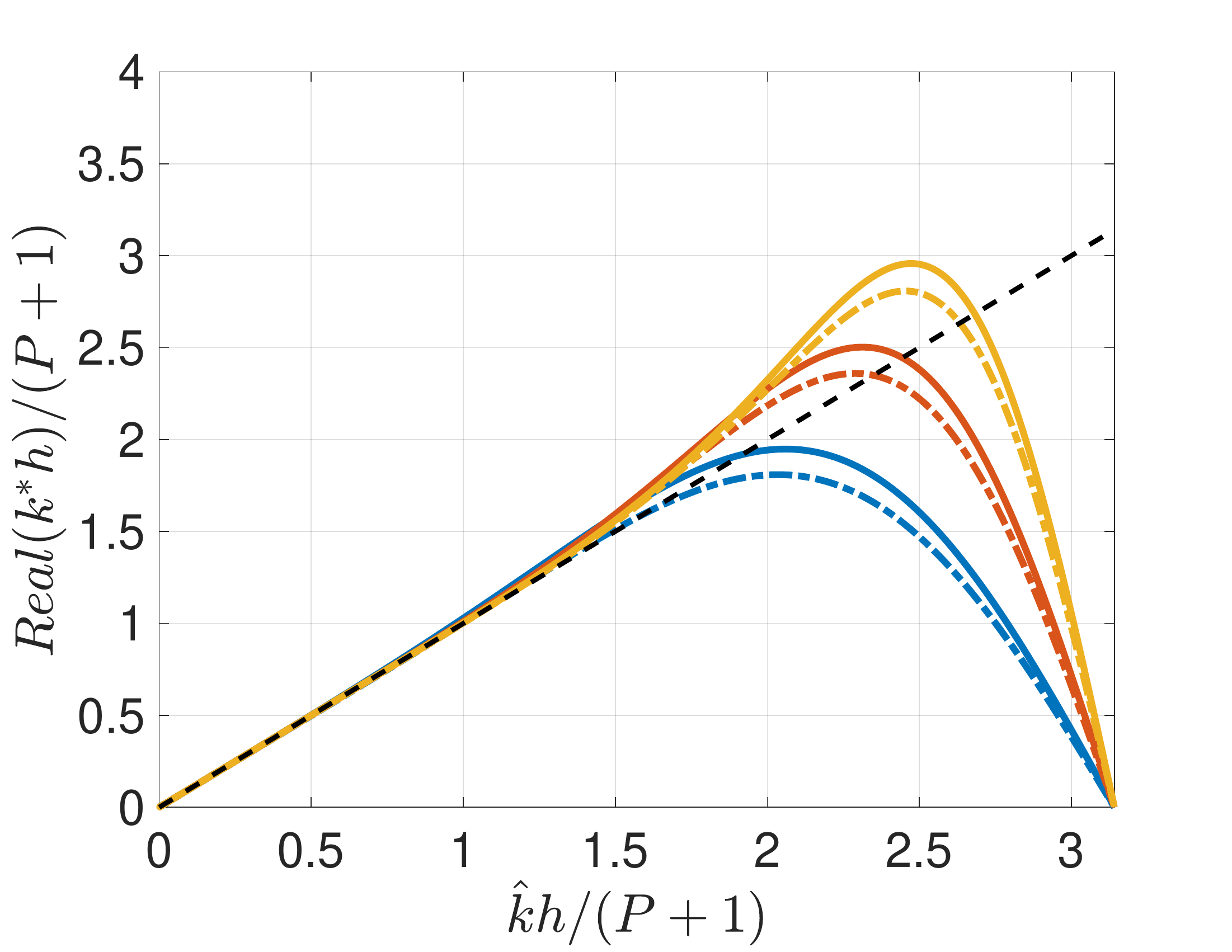}
		\caption{}
	\end{subfigure}
	\begin{subfigure}{.3\textwidth}
		\includegraphics[width=120pt]{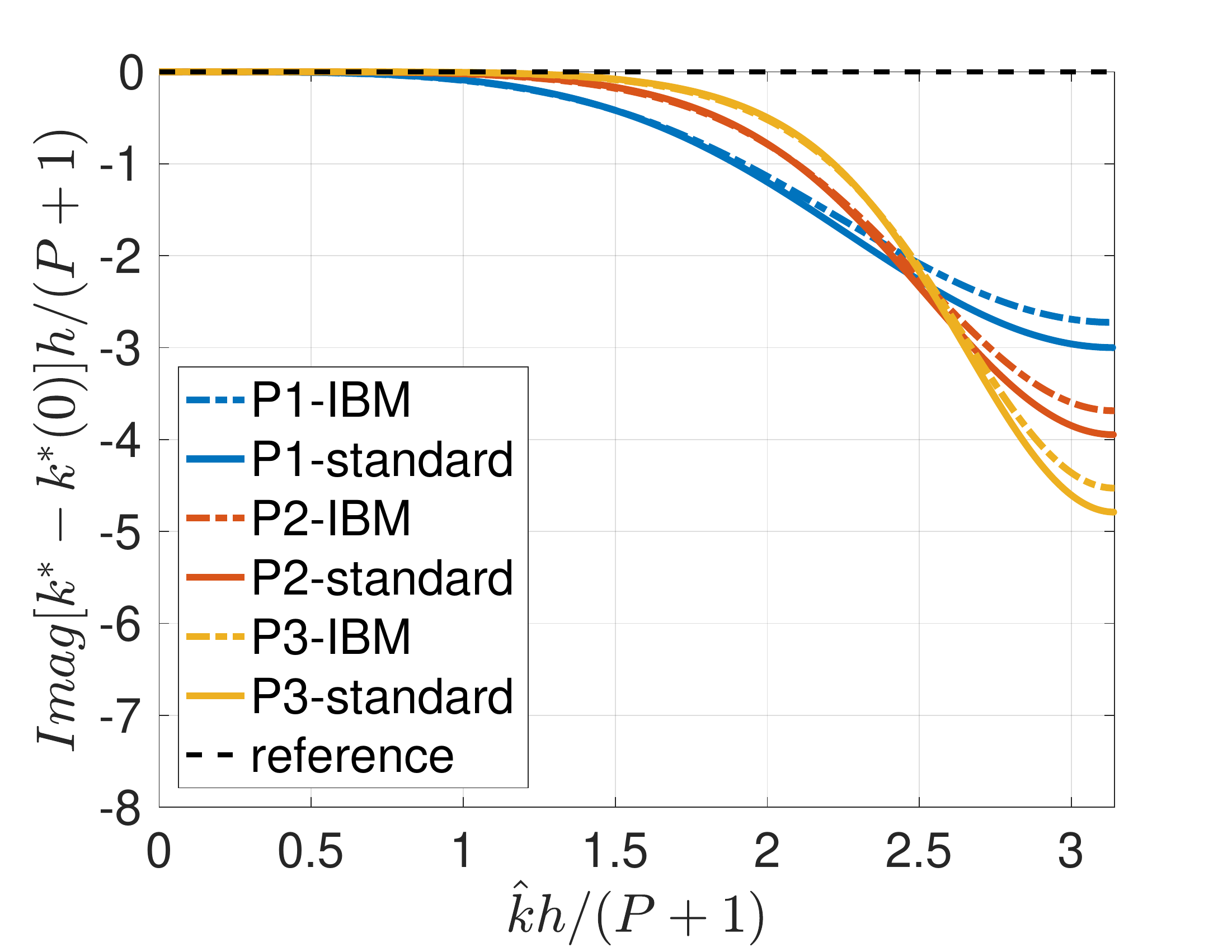}
		\caption{}
	\end{subfigure}
	\begin{subfigure}{.3\textwidth}
		\includegraphics[width=120pt]{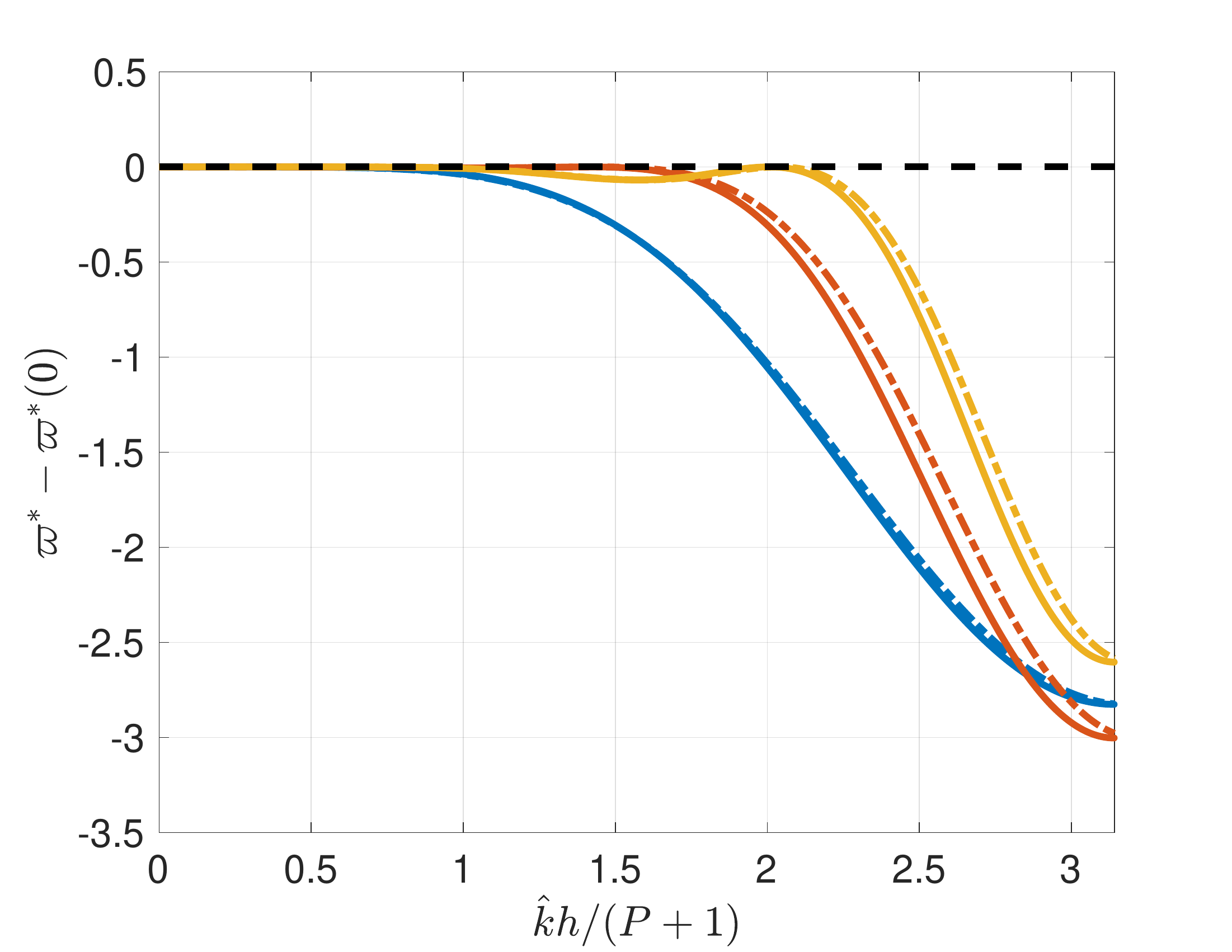}
		\caption{}
	\end{subfigure}
	\begin{subfigure}{.3\textwidth}
		\includegraphics[width=120pt]{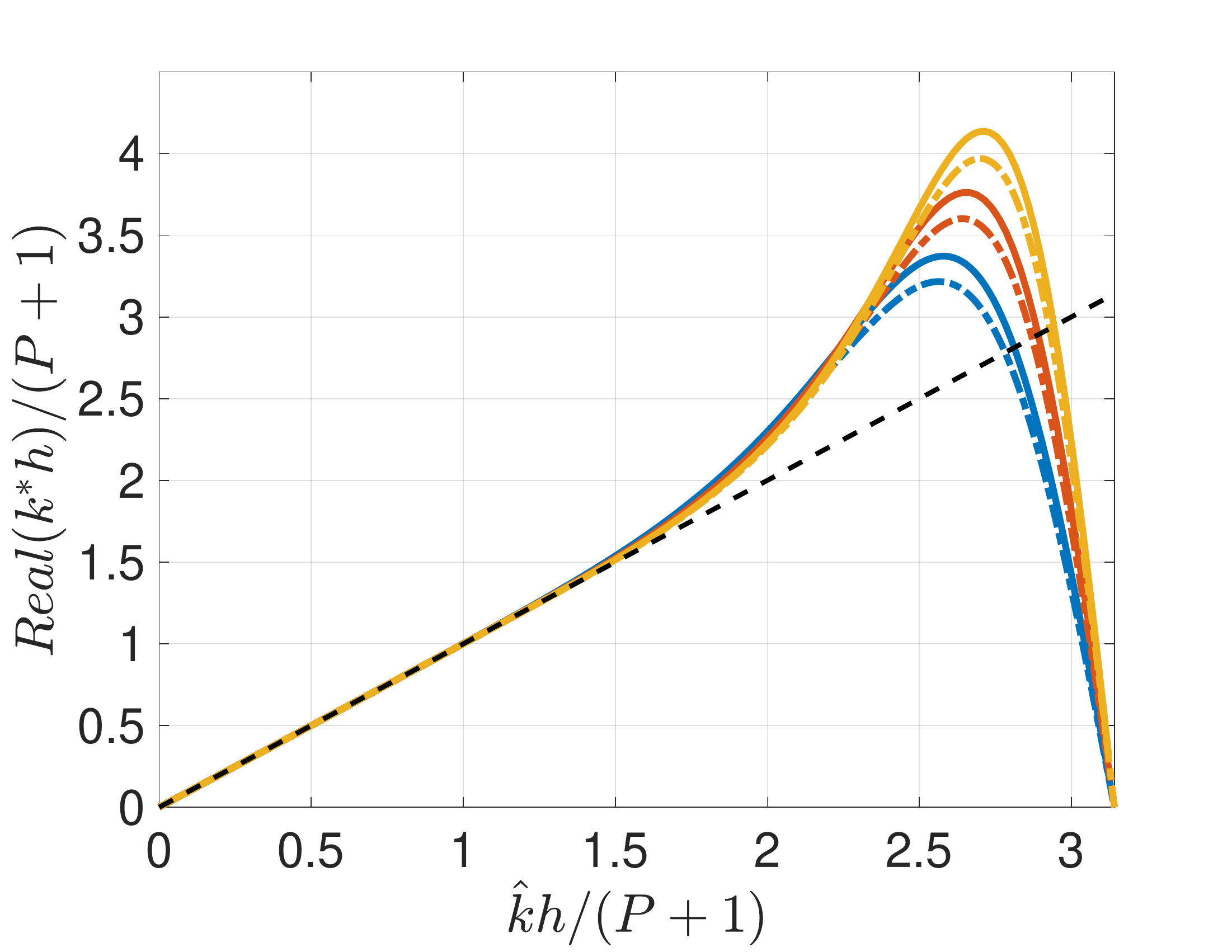}
		\caption{}
	\end{subfigure}
	\begin{subfigure}{.3\textwidth}
		\includegraphics[width=120pt]{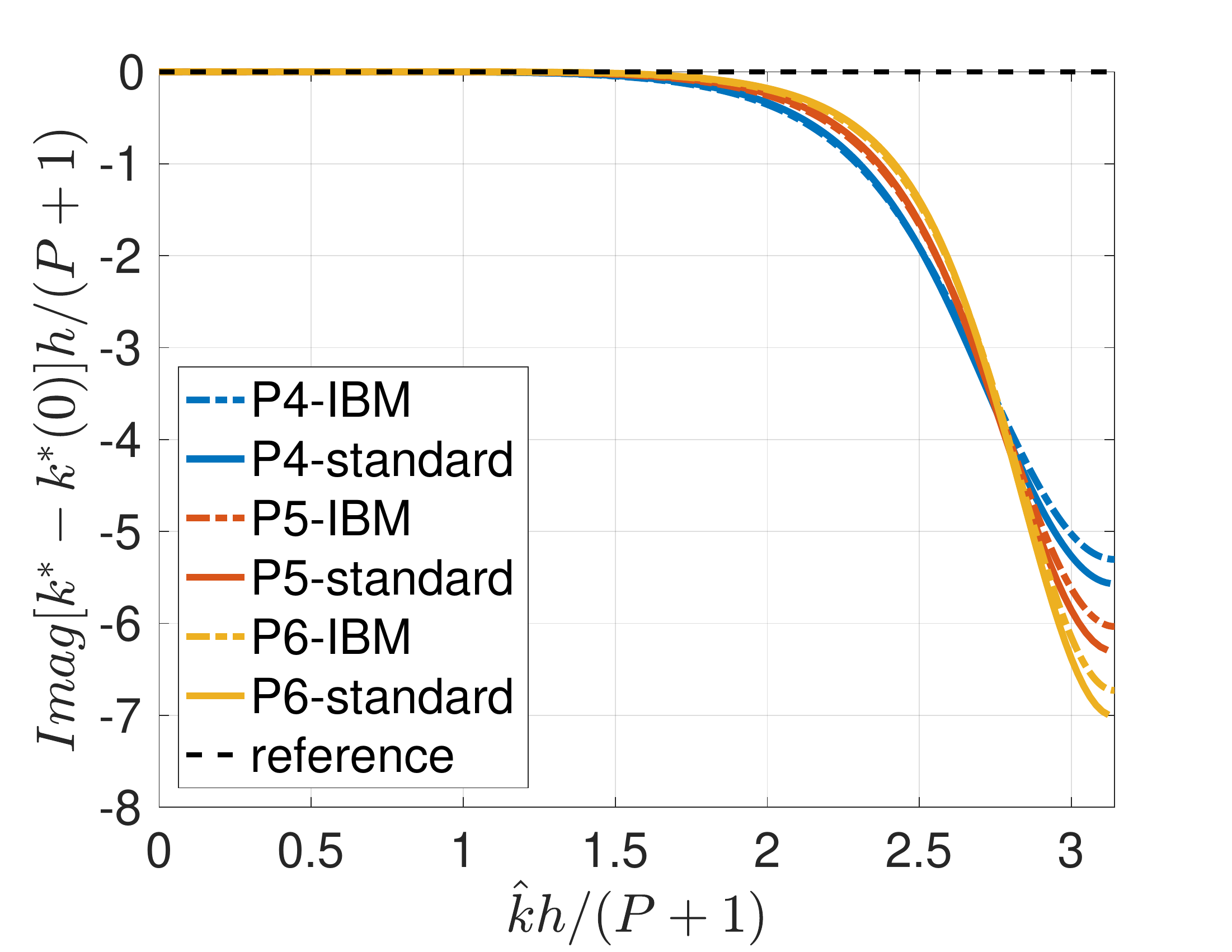}
		\caption{}
	\end{subfigure}
	\begin{subfigure}{.3\textwidth}
		\includegraphics[width=120pt]{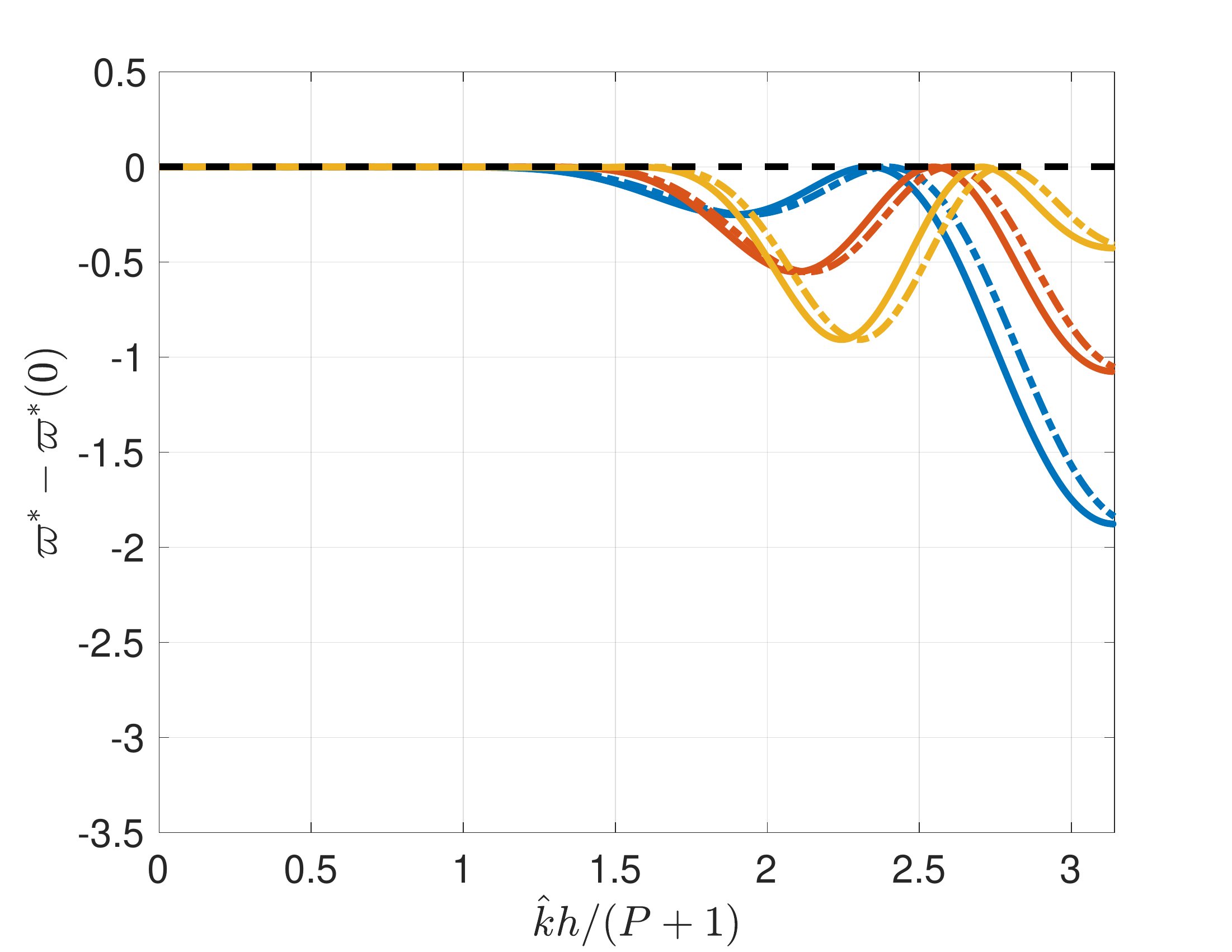}
		\caption{}
	\end{subfigure}
	\centering
	\caption{Eigensolution and non-modal analysis of advection equation with IBM ($N = 40$, $r = 1/40$, $\IBMparam = 1 \times 10^{-4}$) and comparison with those of standard advection equation. a) Dispersion for $P = 1,2,3$. b) Dissipation subtracted by the IBM-induced dissipation $\gamma_{IBM}$ (dissipation at $k=0$) for $P = 1,2,3$. c) Short-term dissipation subtracted by the value at $k=0$ for $P = 1,2,3$. d) Dispersion for $P = 4,5,6$. e) Dissipation subtracted by the IBM-induced dissipation $\gamma_{IBM}$ (dissipation at $k=0$) for $P = 4,5,6$. f) Short-term dissipation subtracted by the value at $k=0$ for $P = 4,5,6$.}
	\label{fig:semi-poly}
\end{figure*}

We then look at the dispersion-dissipation and non-modal behaviors changing with polynomial order $P$, in order to show the consistency with the standard FR scheme. We fix the penalization parameter to $\IBMparam = 1 \times 10^{-4}$, and set $N=40$ with $r=1/40$. The results are shown in Figure \ref{fig:semi-poly}, where the IBM-induced dissipation is subtracted to show $\gamma_{FR}^*$. Good agreement with the standard advection equation in the resolved wavenumber range is shown. This range is also increasing as we go to higher $P$. These results indicate that volume penalization guarantees a good spectral behavior for high-order methods in the resolved wavenumber range, thus can be a good candidate for high-order schemes involving IBM treatment.

\subsection{Fully-discrete analysis}
Fully-discrete analysis considers both time and space discretization and is able to evaluate the stability of the space-time system. To march the solution in time, we choose a standard explicit three-stage, third-order Runge-Kutta scheme for time integration, which leads to the following fully-discrete operator matrix \citep{vermeire2017behaviour,he2020dispersion}
\begin{equation}
    \boldsymbol{A} = \boldsymbol{I} + \Delta t \boldsymbol{M} + \frac{1}{2} (\Delta t \boldsymbol{M})^2 + \frac{1}{6} (\Delta t \boldsymbol{M})^3.
\end{equation}

Note here we take the standard advection equation $M$ as an example, while it can be easily generalized to penalized advection equation by changing the discretization operator matrix to the one in Equation \ref{eq:IBMglobal}. In addition, since the short-term dissipation is only defined when $t \rightarrow 0$, non-modal analysis will not be considered. For a given $k$ the eigenvalues related to the modified wavenumber are given by $\lambda_m = e^{-i \omega^*_m \Delta t}$, therefore we have the numerical dispersion and dissipation defined as
\begin{equation}
    \text{Real}(k^*_m) = \frac{\text{Real}(\omega^*_m)}{\advcoef} = \frac{\text{Real}[i \ln(\lambda_m)]}{\advcoef \Delta t}\\, \ \text{Imag}(k^*_m) = \frac{\text{Imag}(\omega^*_m)}{\advcoef} = \frac{\text{Imag}[i \ln(\lambda_m)]}{\advcoef \Delta t}.
\end{equation}

The fully-discrete dispersion and dissipation curves are again normalized by the smallest length scale $h/(P+1)$ and the initial wavenumber is again re-scaled based on $r$. Since the dissipation curve obtained from fully-discrete analysis indicates the stability of the space-time scheme, it can help to investigate the stability criterion for the selection of penalization parameter. The scheme becomes unstable when any one of the eigenmodes shows positive dissipation. Again, we consider the computational domain discretized into $N = 40$ elements with one solid element $r = 1/40$. The polynomial order $P = 3$ is used which is a typical value used for high-order schemes. To set the penalization parameter we firstly consider the guideline $\IBMparam = \Delta t$ \cite{kolomenskiy2009fourier,engels2015numerical}. This case is compared with results from \citep{vermeire2017behaviour}, where fully-discrete analysis of the advection equation based FR discretization is considered. When time discretization is considered, the effect of time step $\Delta t$ should be taken into account, which is determined by the Courant–Friedrichs–Lewy (CFL) condition. The CFL number is defined as \citep{cockburn1989tvb,chalmers2014relaxing}

\begin{equation}
    CFL = \advcoef  \Delta t  (2P+1) / h.
\end{equation}

\begin{figure*}[htbp]
\begin{subfigure}{.3\textwidth}
		\includegraphics[width=120pt]{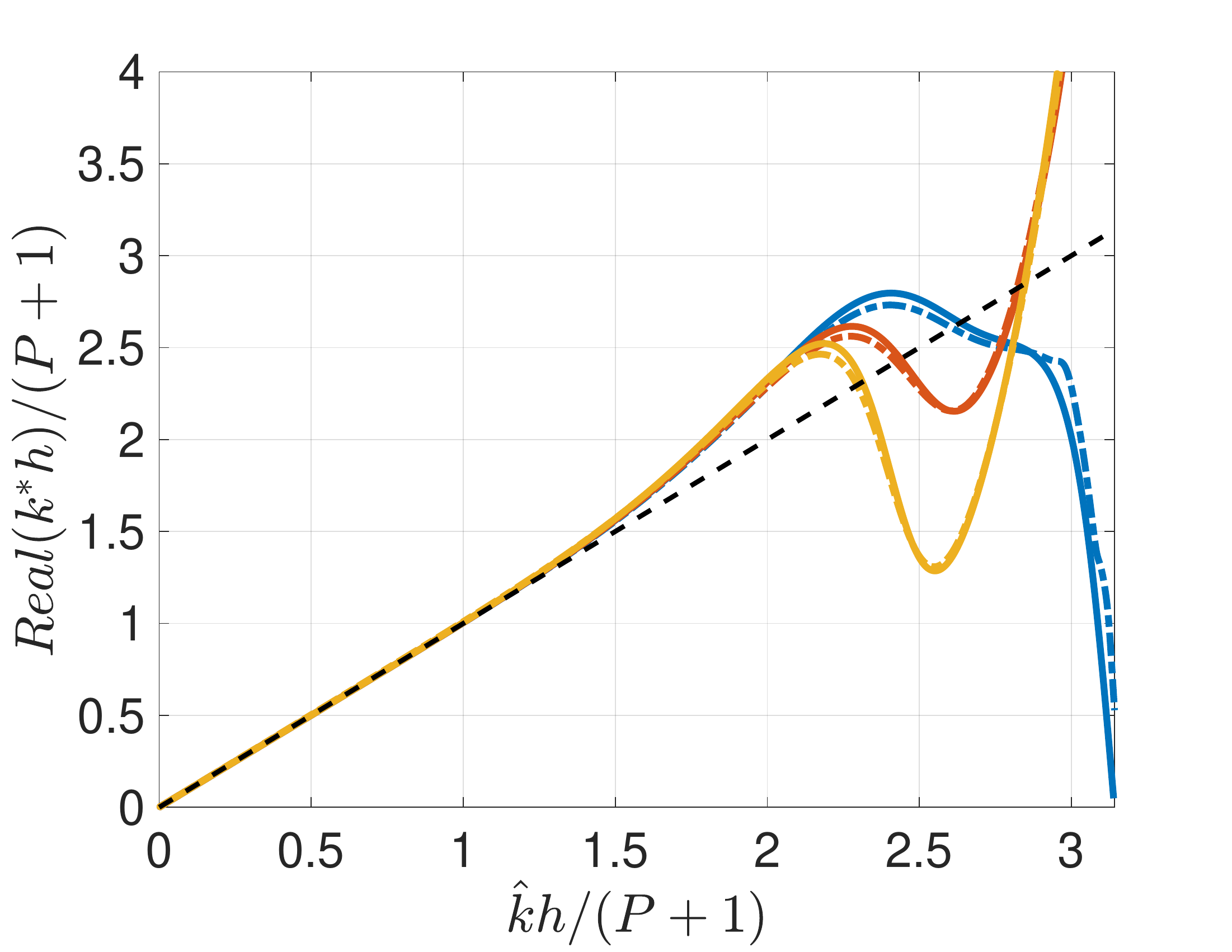}
		\caption{}
	\end{subfigure}
	\begin{subfigure}{.3\textwidth}
		\includegraphics[width=120pt]{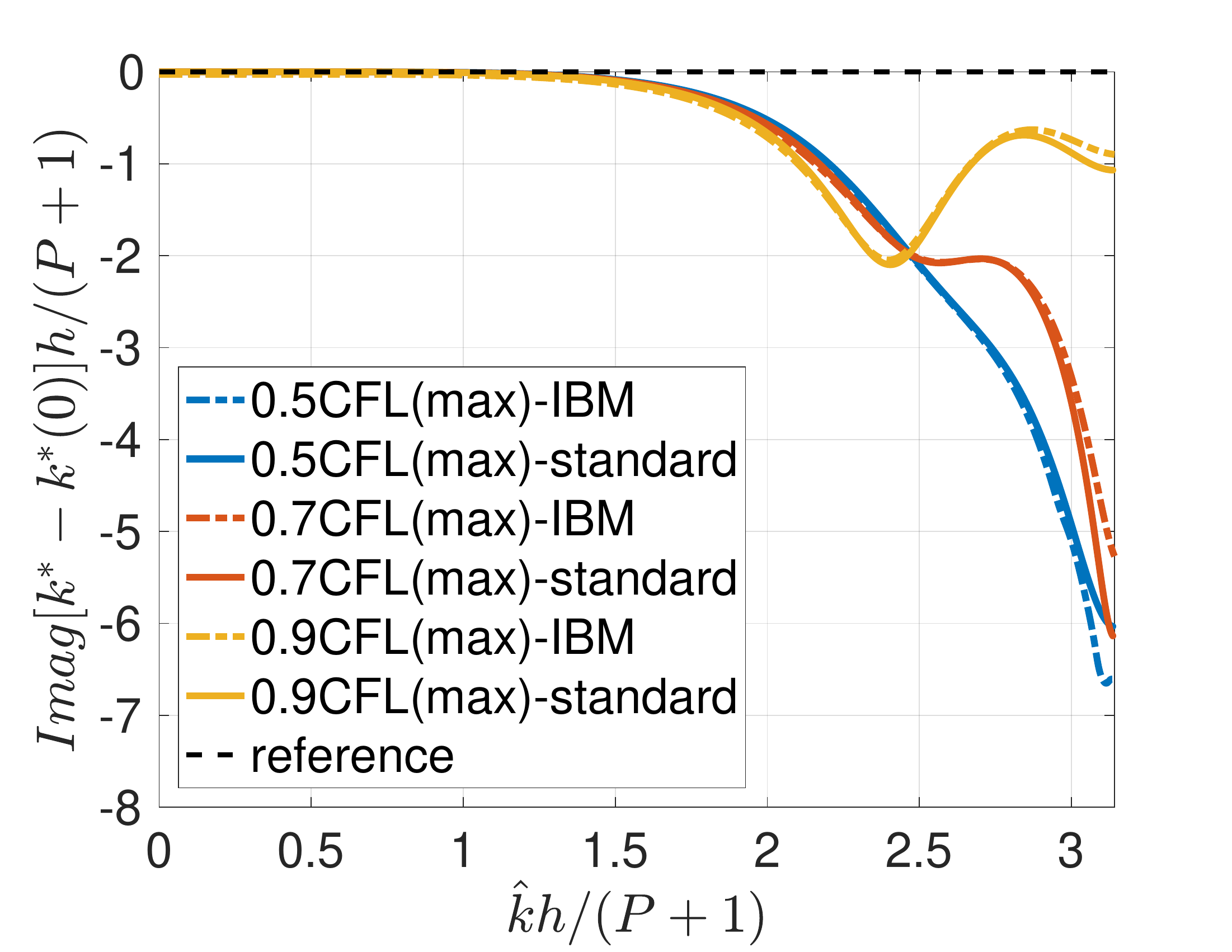}
		\caption{}
	\end{subfigure}
	\centering
	\caption{Fully-discrete analysis of advection equation with IBM ($N = 40$, $r = 1/40$, $P = 3$, $\IBMparam = 1 \times 10^{-4}$) and comparison with those of standard advection equation, under different CFL numbers. Third-order Runge-Kutta method is considered for time integration. a) Dispersion. b) Dissipation subtracted by the IBM-induced dissipation $\gamma_{IBM}$ (dissipation at $k=0$).}
	\label{fig:full-CFL}
\end{figure*}

Following \citep{vermeire2017behaviour}, the maximum time-step for stability is firstly determined and given by CFL(max). Then we select three typical CFL numbers through multiplying this CFL number by $0.5$, $0.7$ and $0.9$. The physical modes are compared in Figure \ref{fig:full-CFL}, where the dispersion-dissipation behavior between IBM and standard advection equation shows good agreement. Very slight difference is seen when $k > 2$, which again belongs to the under-resolved wavenumber region. These results indicate that IBM also has very slight impact on the spectral behavior when both space and time discretizations are considered, and the penalization parameter is chosen within the stability limit. 

\begin{figure*}[htbp]
\begin{subfigure}{.3\textwidth}
		\includegraphics[width=120pt]{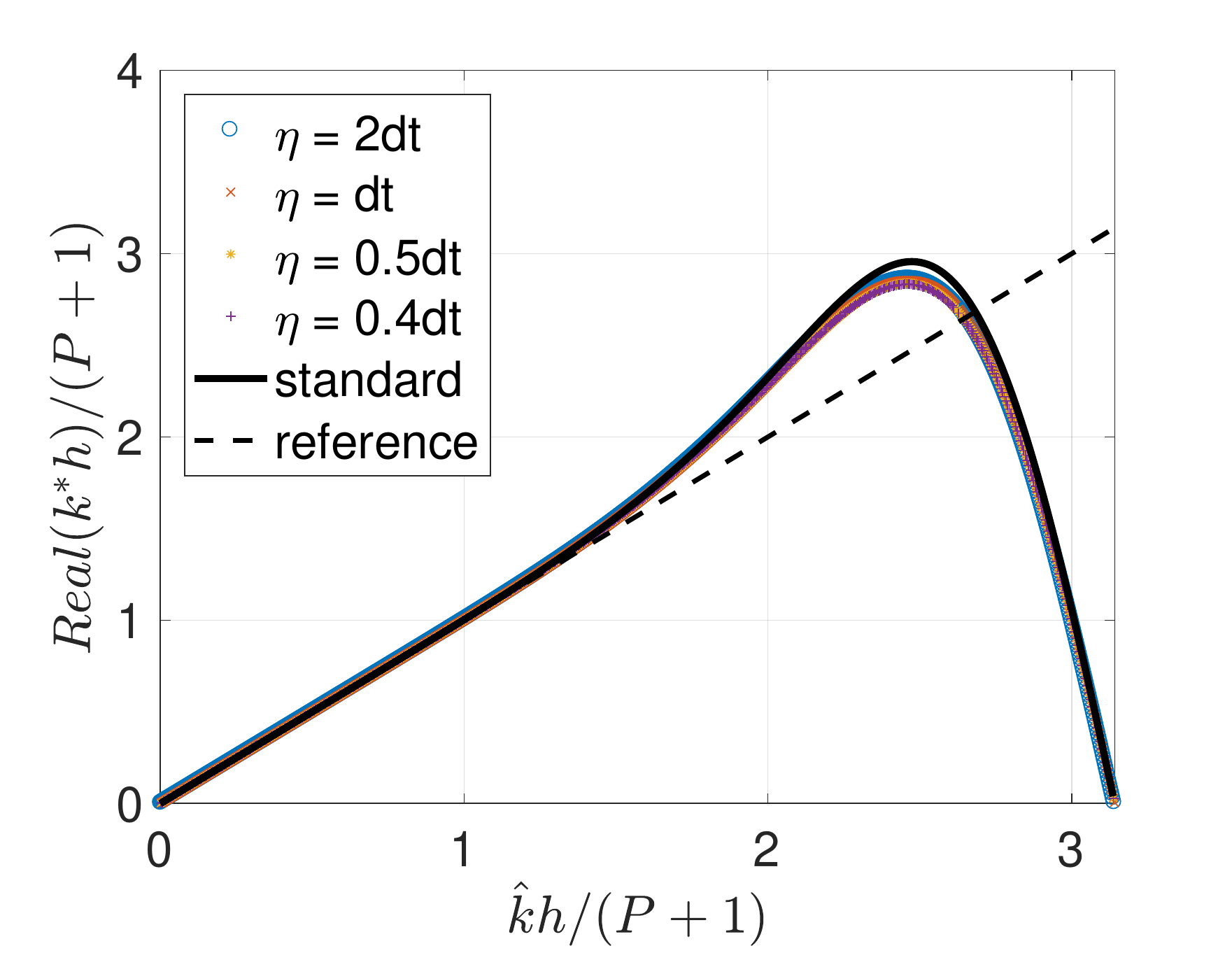}
		\caption{}
	\end{subfigure}
	\begin{subfigure}{.3\textwidth}
		\includegraphics[width=120pt]{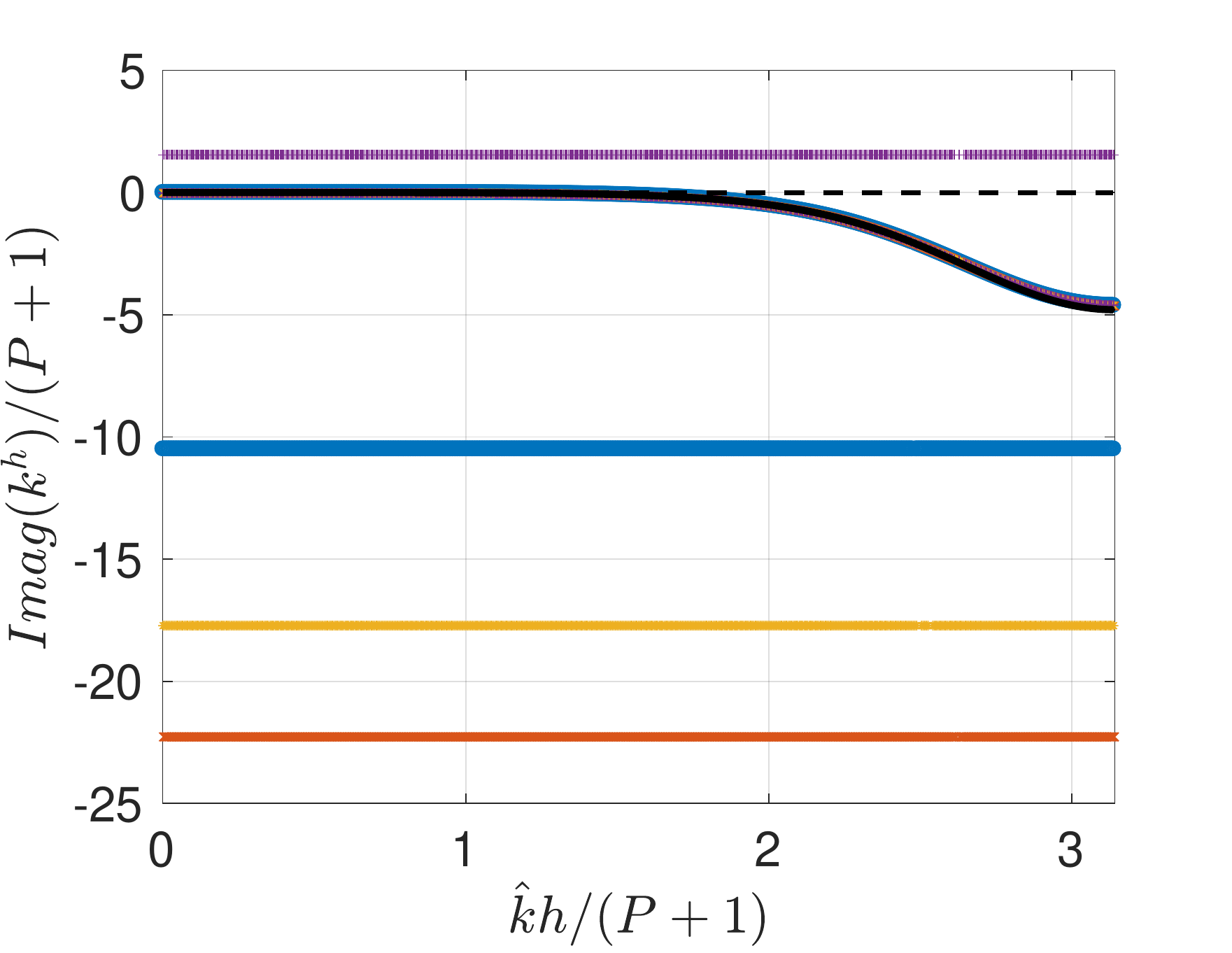}
		\caption{}
	\end{subfigure}
	\centering
	\caption{Fully-discrete analysis of advection equation with IBM ($N = 40$, $r = 1/40$, $P = 3$, $\IBMparam = 1 \times 10^{-4}$, CFL $= 0.1$ CFL(max)) and comparison with those of standard advection equation. Third-order Runge-Kutta method is considered for time integration. a) Dispersion. b) Dissipation of both primary mode and the most unstable solid mode.}
	\label{fig:full-cfl01}
\end{figure*}

\begin{figure*}[htbp]
\begin{subfigure}{.3\textwidth}
		\includegraphics[width=120pt]{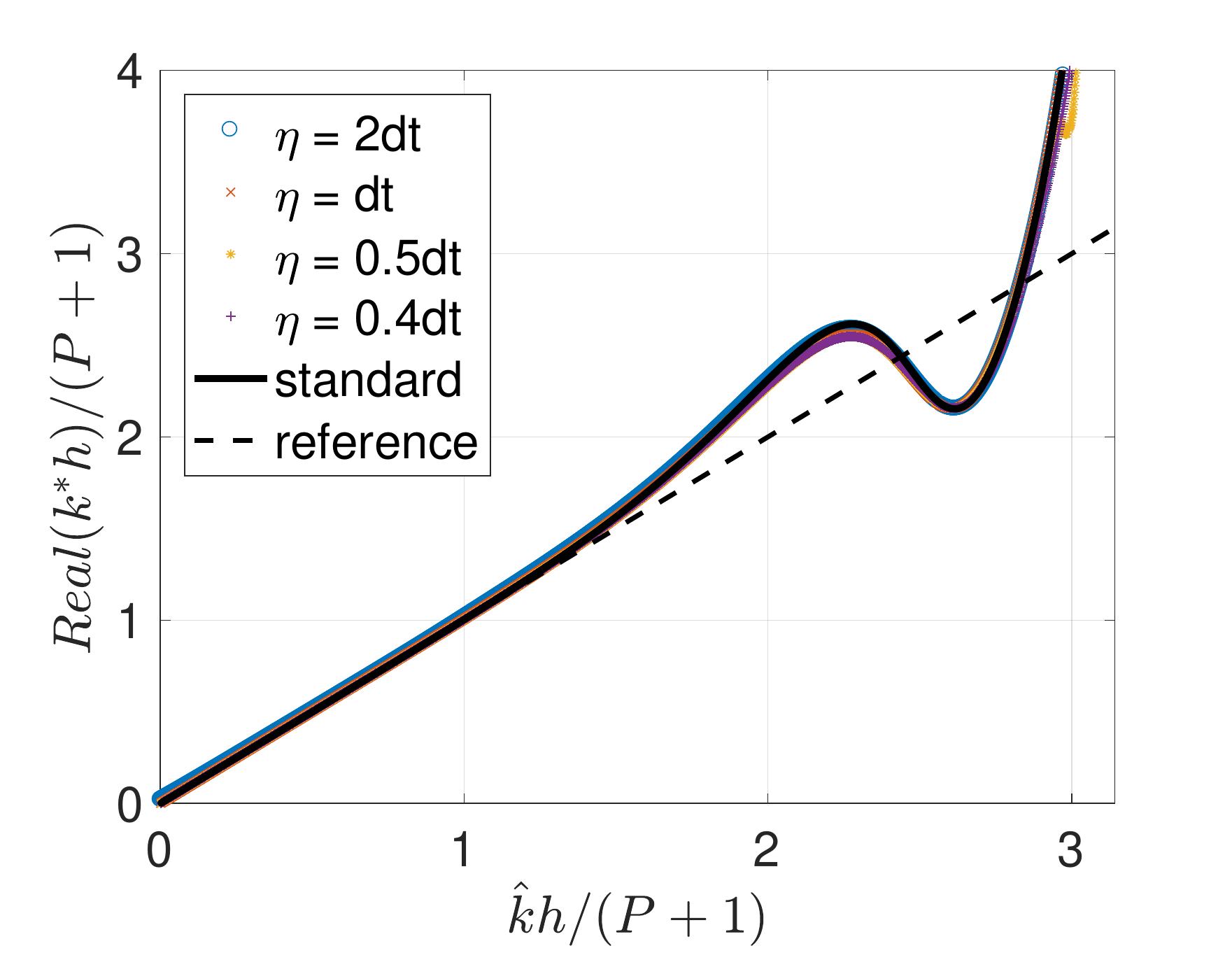}
		\caption{}
	\end{subfigure}
	\begin{subfigure}{.3\textwidth}
		\includegraphics[width=120pt]{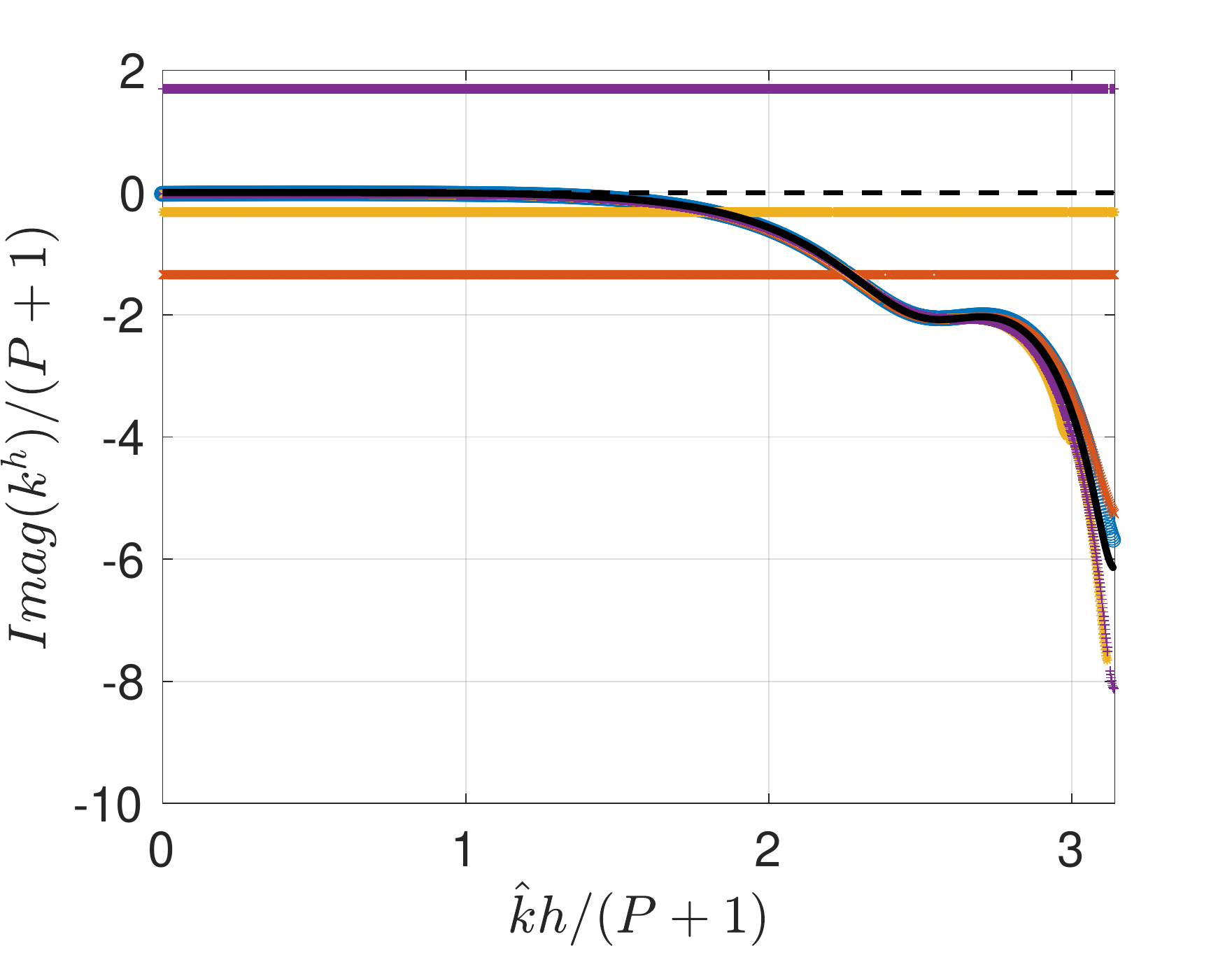}
		\caption{}
	\end{subfigure}
	\centering
	\caption{Fully-discrete analysis of advection equation with IBM ($r = 1/40$, $N = 40$, $P = 3$, $\IBMparam = 1 \times 10^{-4}$, CFL $= 0.7$ CFL(max)) and comparison with those of standard advection equation. Third-order Runge-Kutta method is considered for time integration. a) Dispersion. b) Dissipation of both primary mode and the most unstable solid mode.}
	\label{fig:full-cfl07}
\end{figure*}

Furthermore it has been shown that even though the time step is chosen below the critical CFL number, instability can still occur when $\IBMparam$ is too small since the problem becomes very stiff in the limit $\IBMparam \rightarrow 0$ \citep{schneider2015immersed}. Therefore, the critical value for penalization parameter, with respect to the time step, needs to be carefully studied to guide the selection of penalization parameter. From previous studies, we introduced some variations to $\IBMparam$ around $\IBMparam = \Delta t$ and look at the instability behavior based on the fully-discrete dissipation curve. As shown in Figure \ref{fig:example}, there are several solid modes among secondary modes induced by volume penalization, with constant dispersion and dissipation across all wavenumbers. When the penalization parameter varies, these solid modes will move and once $\IBMparam$ becomes smaller than the critical value $\IBMparam_{critical}$, these modes will show positive dissipation, indicating the numerical instability occurs. In order to investigate the range of $\IBMparam_{critical}$ and obtain more accurate estimation of this parameter, the variation of both penalization parameter and the CFL number is studied. The dispersion behavior of the physical mode, and the dissipation behavior of the physical and the most unstable solid mode at two CFL numbers are shown in Figure \ref{fig:full-cfl01} and Figure \ref{fig:full-cfl07}. This solid mode is identified as the most unstable eigenmode with maximal dissipation among all eigenmodes that have a constant value across the wavenumber range. It should be noted that this mode can be merged into other secondary modes when a relatively large $\IBMparam$ is selected with very large negative dissipation, and it will appear when $\IBMparam$ is approaching $\Delta t$. As shown in the figures, when $\IBMparam$ decreases from $0.5 \Delta t$ to $0.4 \Delta t$, the dissipation of this mode becomes positive, indicating that instability occurs. Therefore we can conclude that the instability of the numerical scheme is mainly driven by the instability of the solid mode introduced by IBM. This is consistent with the basic understanding of IBM, where the instability is induced by the stiffness of the IBM source term. Moreover, it also agrees with previous observations \citep{manzanero2018dispersion} that the secondary modes are responsible of the instability, rather than the primary mode. 

\begin{figure*}[htbp]
		\centering
		\includegraphics[width=200pt]{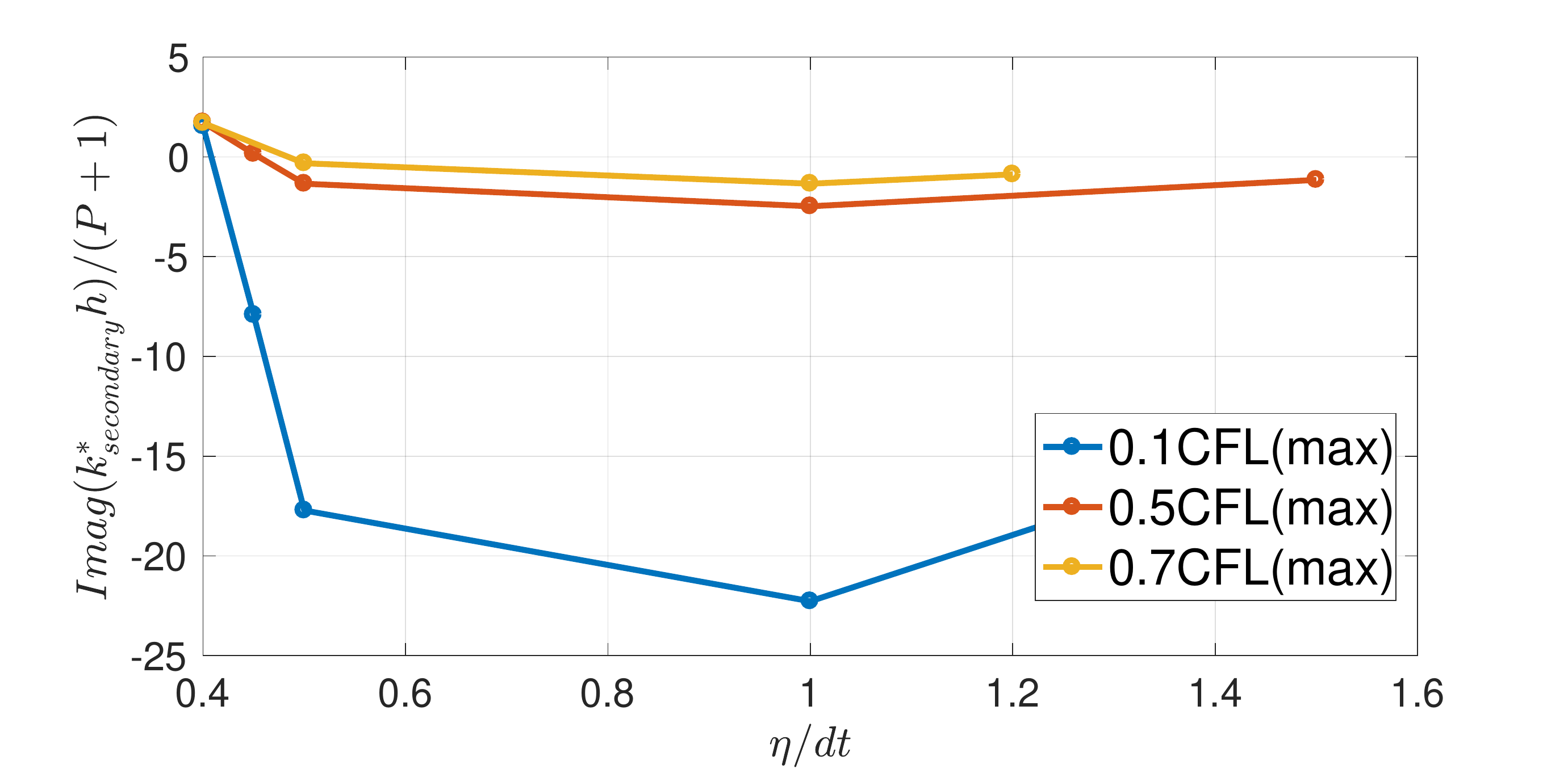}
		\caption{Comparison of the dissipation of the most unstable solid mode under different $\eta$ and CFL number.}
		\label{fig:full-diss}
\end{figure*}

Finally, we compute the dissipation of the solid mode under different CFL numbers and compare the dissipation in Figure \ref{fig:full-diss}. It is found that for different CFL numbers, the critical $\IBMparam$ for stability based on the present time integration is always about $0.4 \Delta t < \IBMparam_{critical} < 0.5 \Delta t$. As the CFL number increases, this value will slightly increase within this range. Compared with previous guidance of $\IBMparam = \Delta t$, this gives a more relaxed guidance to select the penalization parameter $\IBMparam$. 

\section{Numerical experiments}
\subsection{Advection equation with volume penalization}
\begin{figure*}[htbp]
\begin{subfigure}{.3\textwidth}
		\includegraphics[width=120pt]{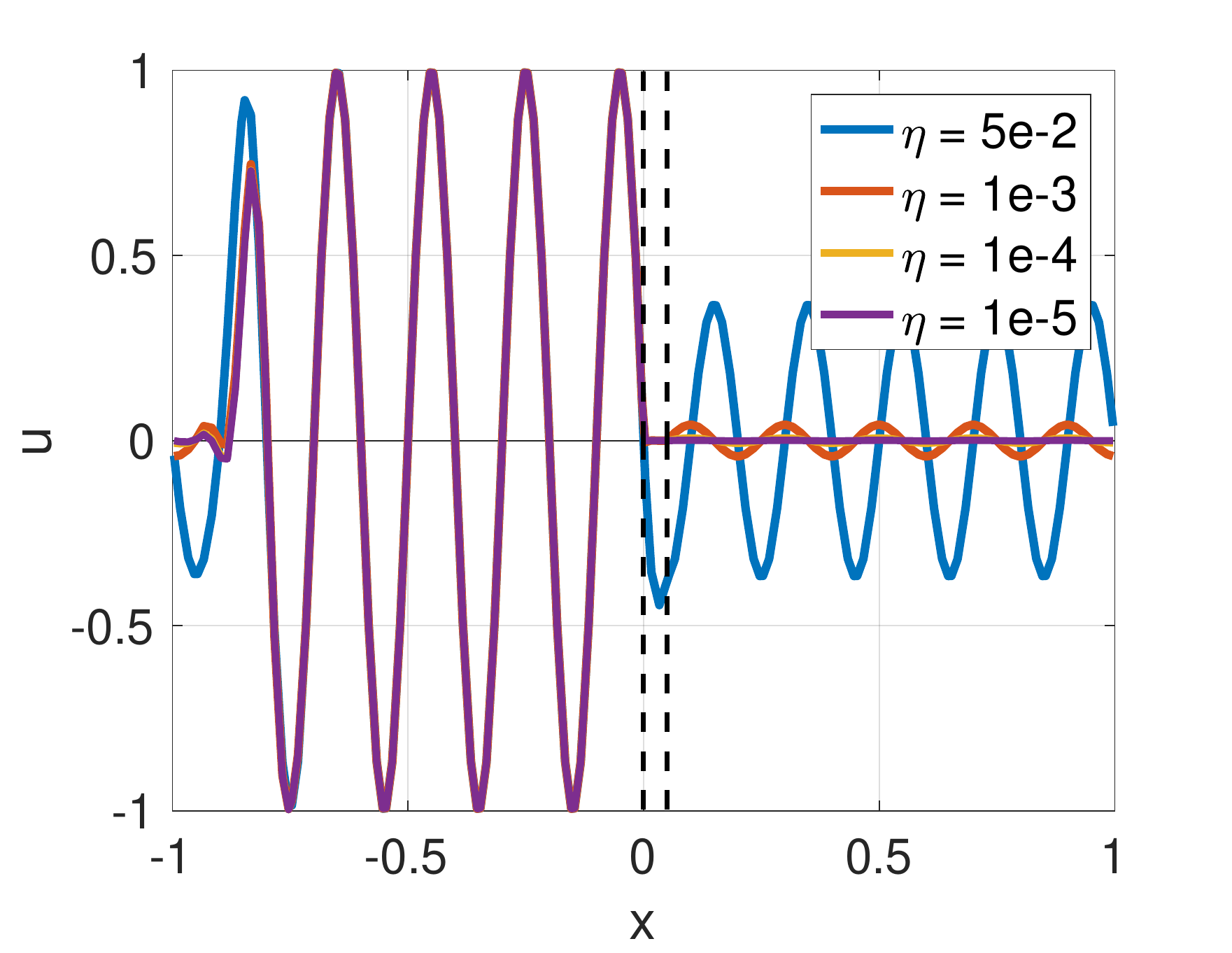}
		\caption{}
	\end{subfigure}
	\begin{subfigure}{.3\textwidth}
		\includegraphics[width=120pt]{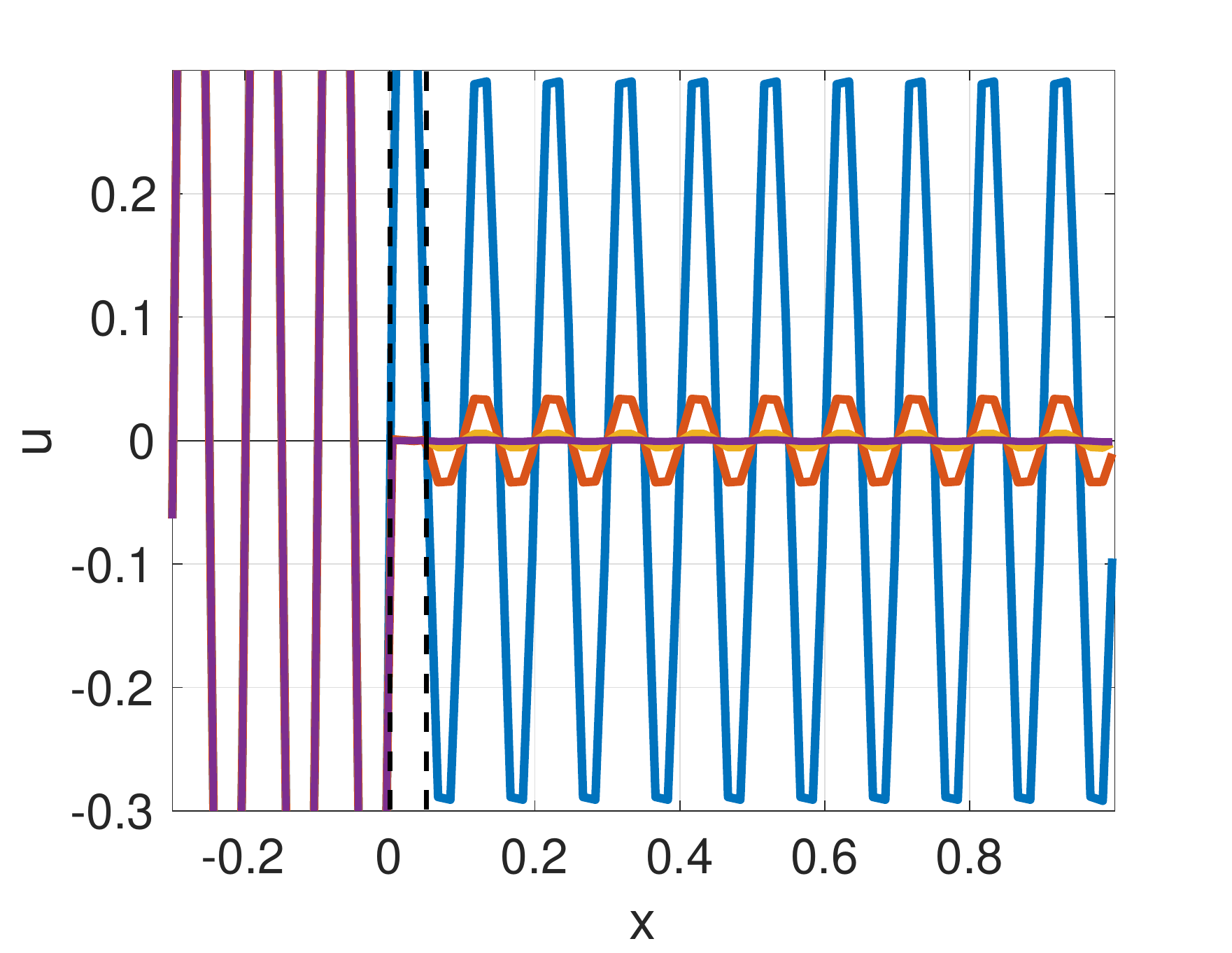}
		\caption{}
	\end{subfigure}
	\begin{subfigure}{.3\textwidth}
		\includegraphics[width=120pt]{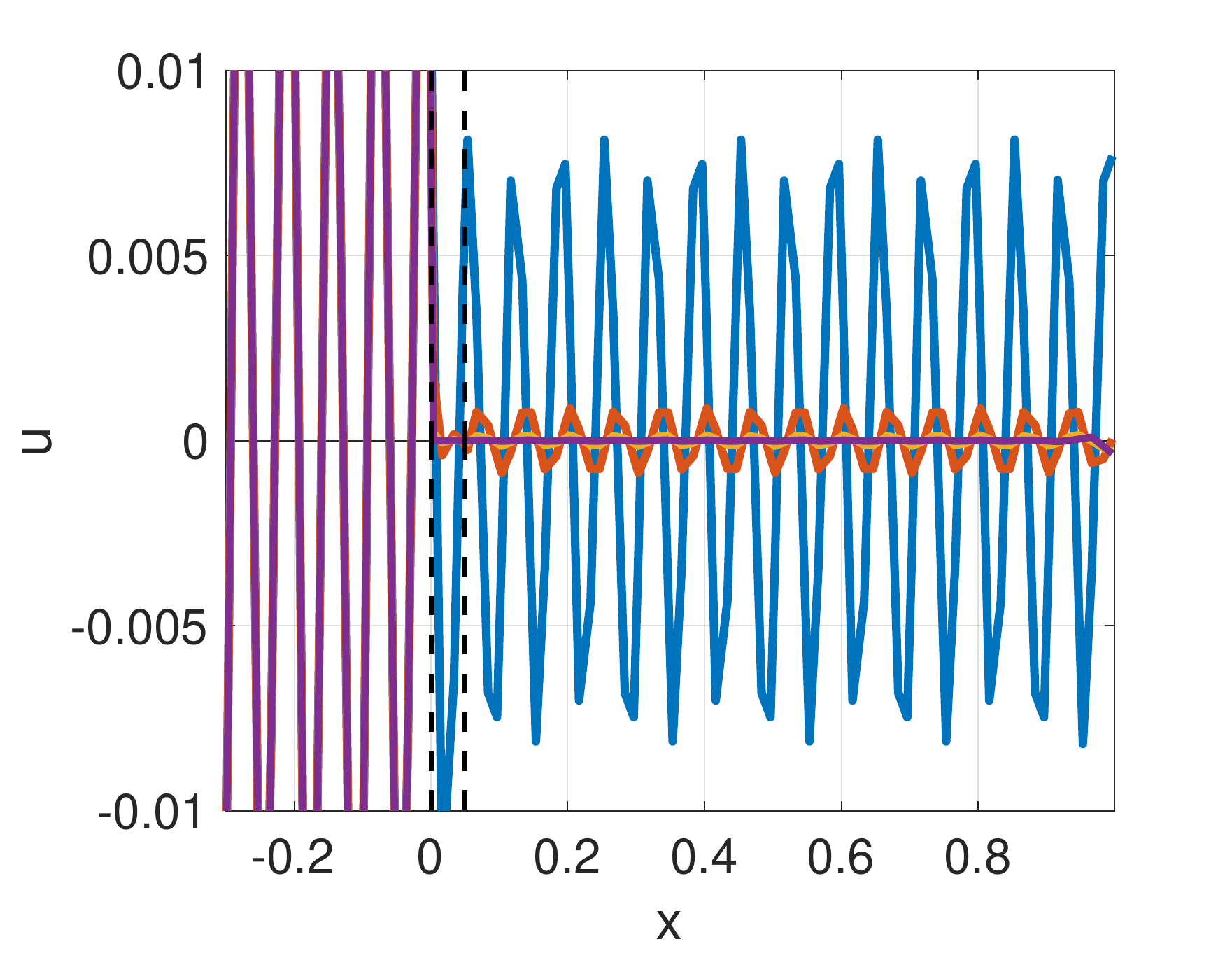}
		\caption{}
	\end{subfigure}
	\centering
	\caption{Simulation of advection equation with IBM under different penalization parameters ($N = 40$, $r = 1/40$, $P = 3$). a) $kh/(P+1) = 0.3927$. b) $kh/(P+1) = 0.7854$. c) $kh/(P+1) = 1.9635$.}
	\label{fig:exp1-sim}
\end{figure*}

In this subsection, the numerical simulation of the one-dimensional advection equation is presented, in order to validate the conclusions drawn in previous sections. We select the same case as before, where $N = 40$, $\Delta = h$ ($r = 1/40$) and $P = 3$. The final simulation time is set to $1.1$, to guarantee that the solution in the computational domain $x \in [\Delta, 1]$ is sufficiently penalized. Since we consider the no-slip wall boundary condition $\velx_s = 0$, the expected solution in this domain should be $0$. The time integration is based on the third-order Runge-Kutta scheme used for fully-discrete analysis. To reduce the temporal error, a sufficiently small time step is set as $\Delta t = 1 \times 10^{-5}$. Three initial conditions ranging from low to medium wavenumbers, including $kh/(P+1) = 0.3927$, $kh/(P+1) = 0.7854$ and $kh/(P+1) = 1.9635$, are considered as test cases. The simulation results with different penalization parameters $\IBMparam$ are shown in Figure \ref{fig:exp1-sim}. It can be seen that as $\IBMparam$ decreases, more damping inside the solid is imposed, therefore the solution shows a smaller amplitude. This agrees with the observation in Figure \ref{fig:semi-penalty}c where more damping is seen on the primary mode when $\IBMparam$ is decreasing. In addition, the influence on frequency is not evident among these cases, indicating the fact that volume penalization only provides additional damping but does not influence significantly the dispersion behavior. It should also be noted that the resulting amplitude among these cases is due to the combined dissipation effects of IBM and the FR scheme, where the pure damping effect from IBM can be extracted by considering a constant initial condition. As $\IBMparam \rightarrow 0$, the solution will converge to zero, indicating that the incoming wave is sufficiently damped by the wall. However, due to the stiffness of the source term, as indicated from the fully-discrete analysis, it is not feasible to make $\IBMparam = 0$ but to use a small $\IBMparam$ instead that balances the accuracy and stability. Moreover, as we go from low to high wavenumber, the amplitude of the solution for a given $\IBMparam$ decreases. This is due to the dissipation effect of the FR scheme, since IBM only has a constant damping across all wavenumbers. 


\begin{figure*}[htbp]
\begin{subfigure}{.3\textwidth}
		\includegraphics[width=120pt]{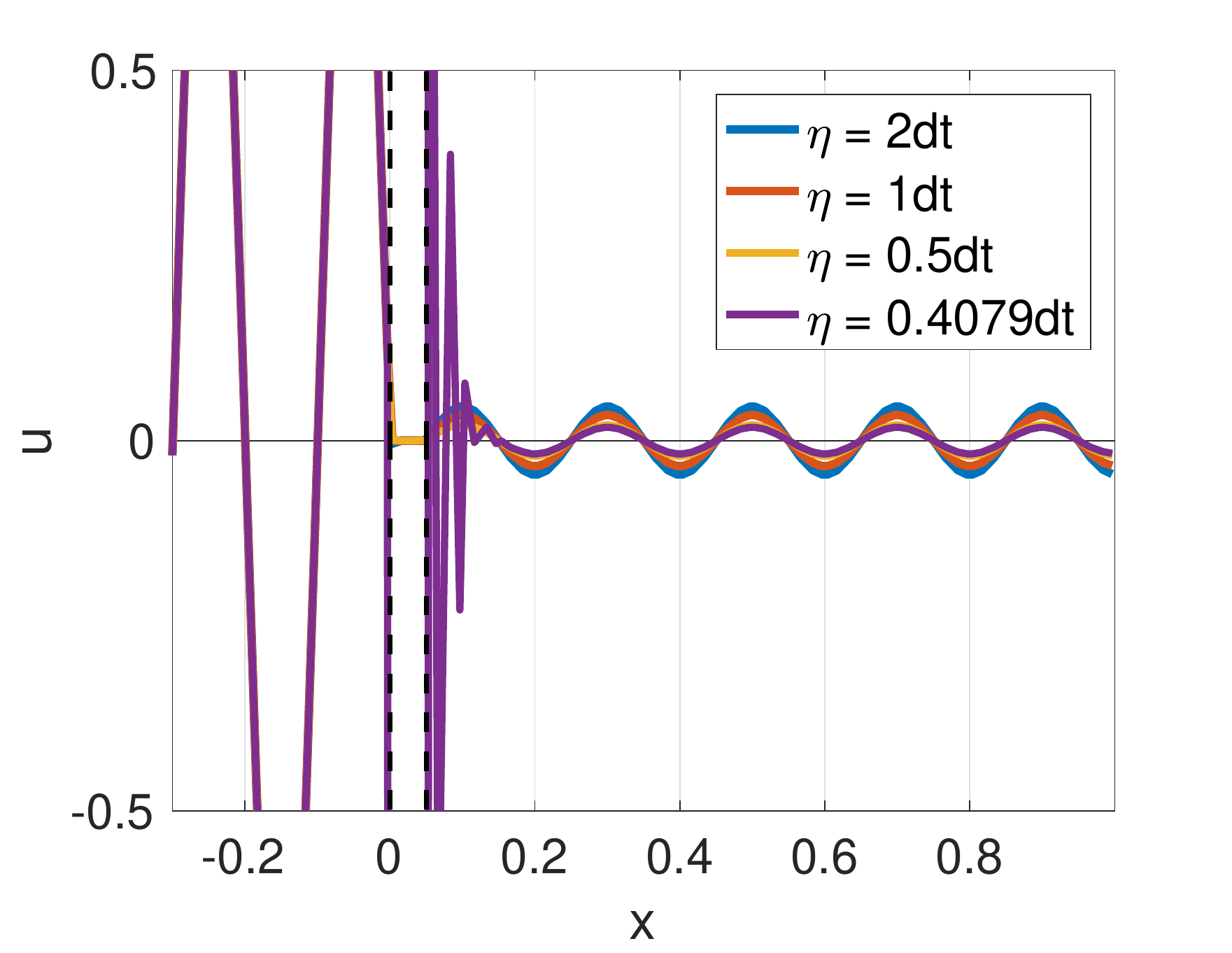}
		\caption{}
	\end{subfigure}
	\begin{subfigure}{.3\textwidth}
		\includegraphics[width=120pt]{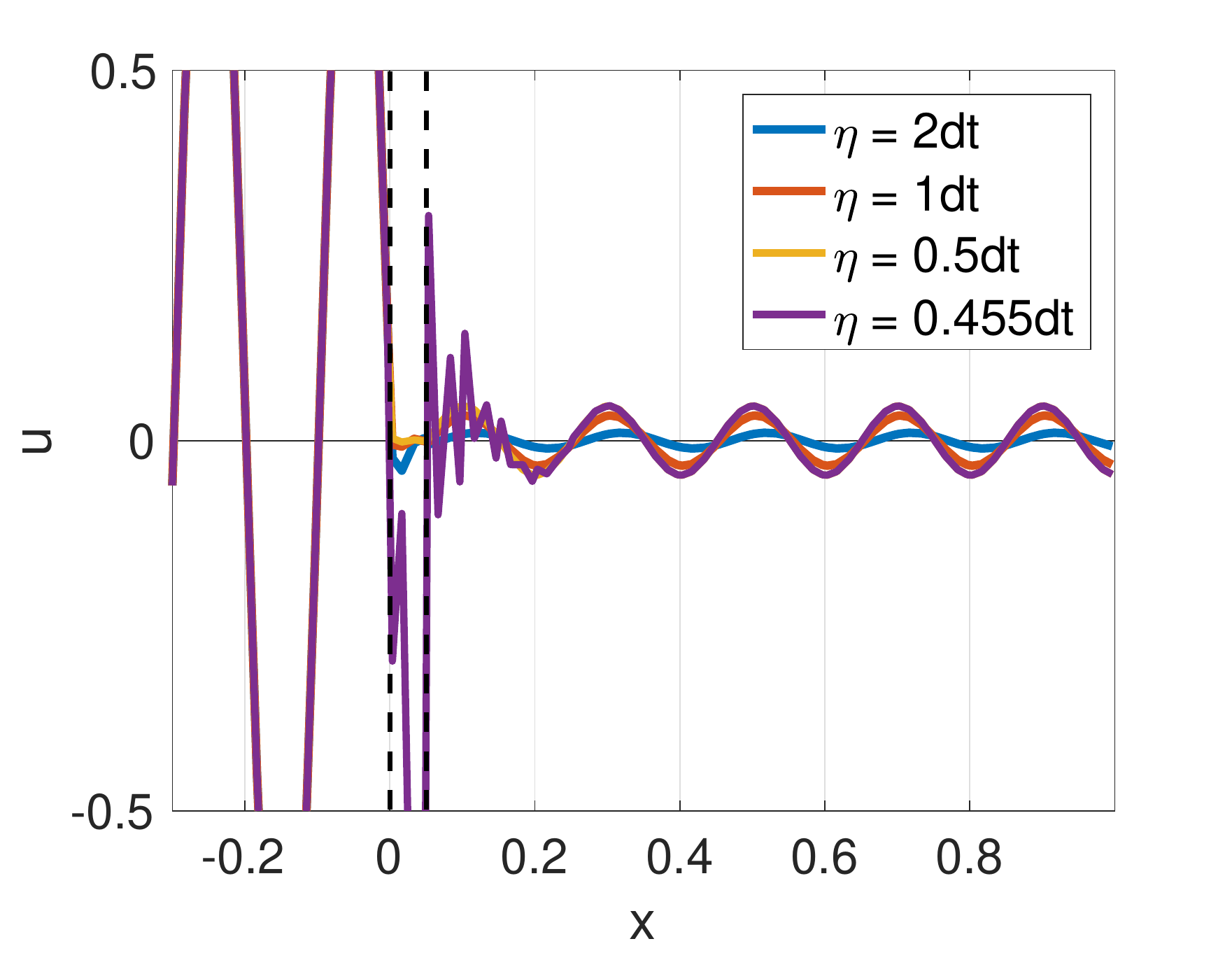}
		\caption{}
	\end{subfigure}
	\begin{subfigure}{.3\textwidth}
		\includegraphics[width=120pt]{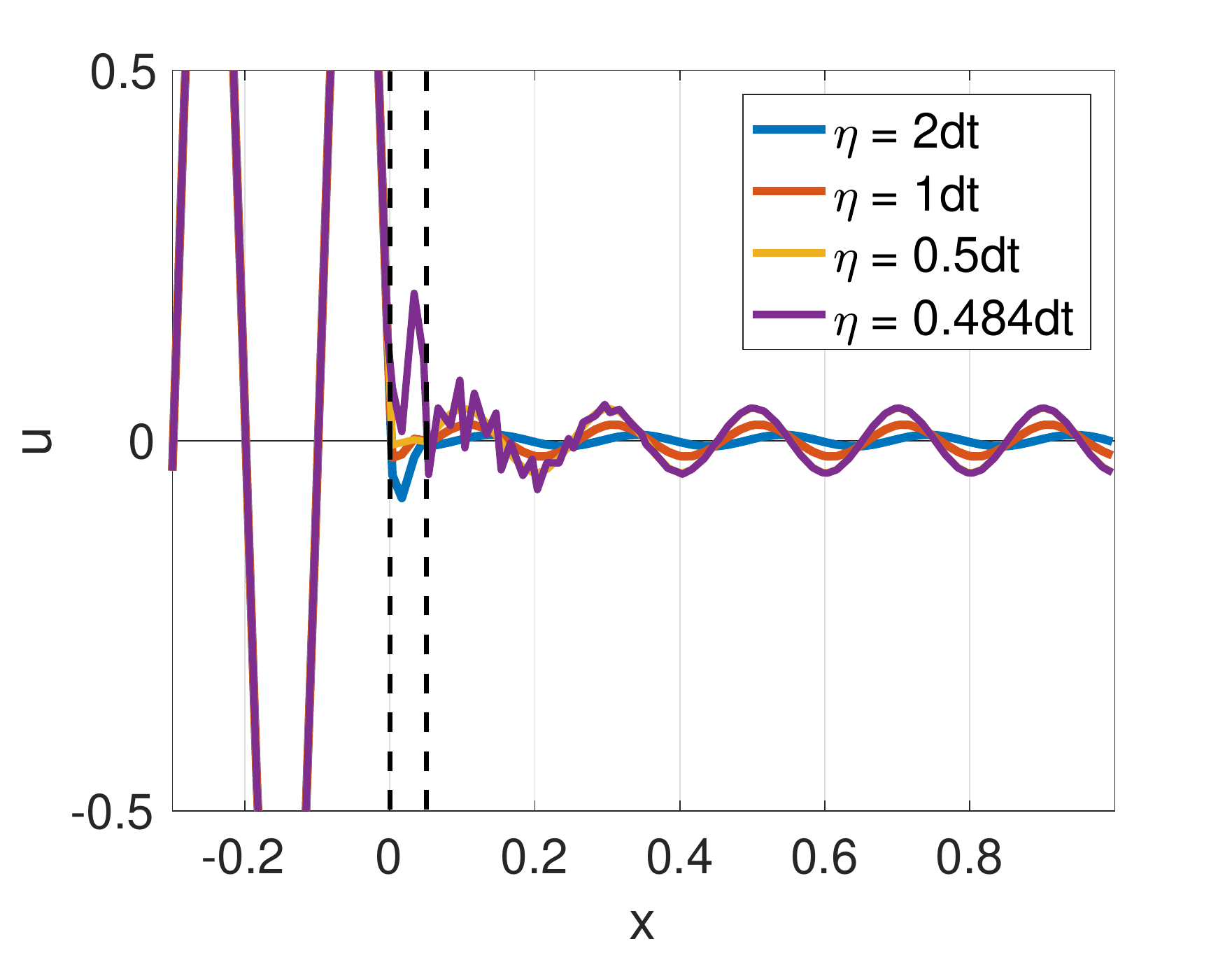}
		\caption{}
	\end{subfigure}
	\centering
	\caption{Simulation of advection equation with IBM under different CFL numbers and penalization parameters ($N = 40$, $r = 1/40$, $P = 3$, $kh/(P+1) = 0.3927$). a) CFL = 0.1CFL(max). b) CFL = 0.5CFL(max). c) CFL = 0.7CFL(max).}
	\label{fig:exp1-cfl}
\end{figure*}

The second set of test cases helps to validate the conclusions obtained from the fully-discrete analysis. We take the same test case as in Figure \ref{fig:exp1-sim}a, where the initial wavenumber is set to $kh/(P+1) = 0.3927$. Three CFL numbers are chosen for investigation, including $0.1$, $0.5$, and $0.7$, which correspond to time step $\Delta t = 6.5 \times 10^{-4}$, $0.00325$ and $0.00455$, respectively. The penalization parameter is gradually reduced from $2\Delta t$ to $0.4\Delta t$, with some typical solutions shown in Figure \ref{fig:exp1-cfl}. From the figure it can be seen that the numerical method remains stable until $\IBMparam = 0.5 \Delta$ for all CFL numbers, which agree with the fully-discrete analysis. Instabilities are observed when $\IBMparam$ approaches $0.4 \Delta t$, whose critical value increases with increasing CFL number. Another interesting observation is that for large CFL numbers, see Figure \ref{fig:exp1-cfl}b and Figure \ref{fig:exp1-cfl}c, reducing $\IBMparam$ does not lead to a reduction of the error, which means for the present problem that there is no need to set $\IBMparam$ too small when the time step is large. 

\subsection{Improvement of the volume penalization through adding second-order terms inside the solid}
\begin{figure*}[htbp]
\begin{subfigure}{.24\textwidth}
		\includegraphics[width=120pt]{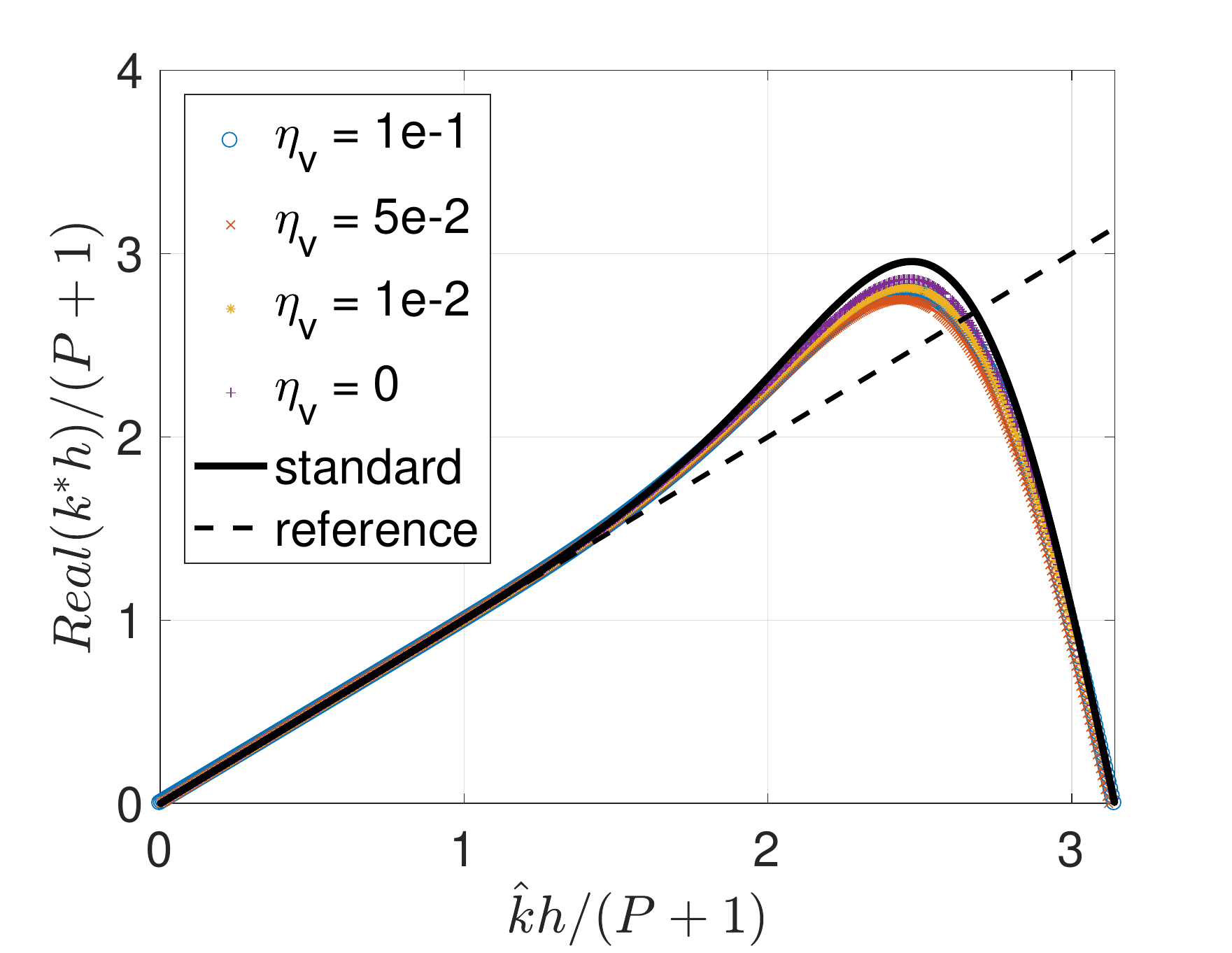}
		\caption{}
	\end{subfigure}
	\begin{subfigure}{.24\textwidth}
		\includegraphics[width=120pt]{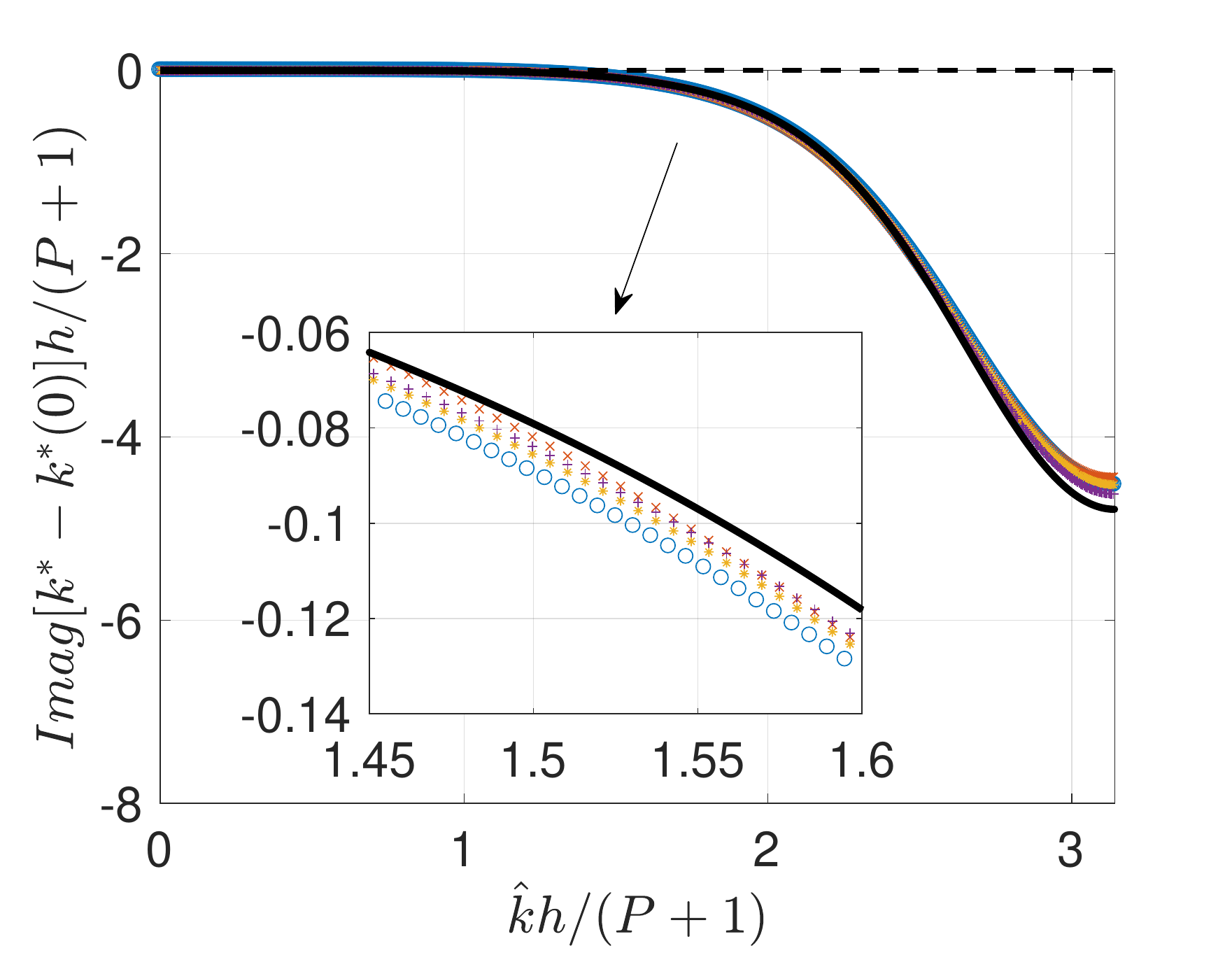}
		\caption{}
	\end{subfigure}
	\begin{subfigure}{.24\textwidth}
		\includegraphics[width=120pt]{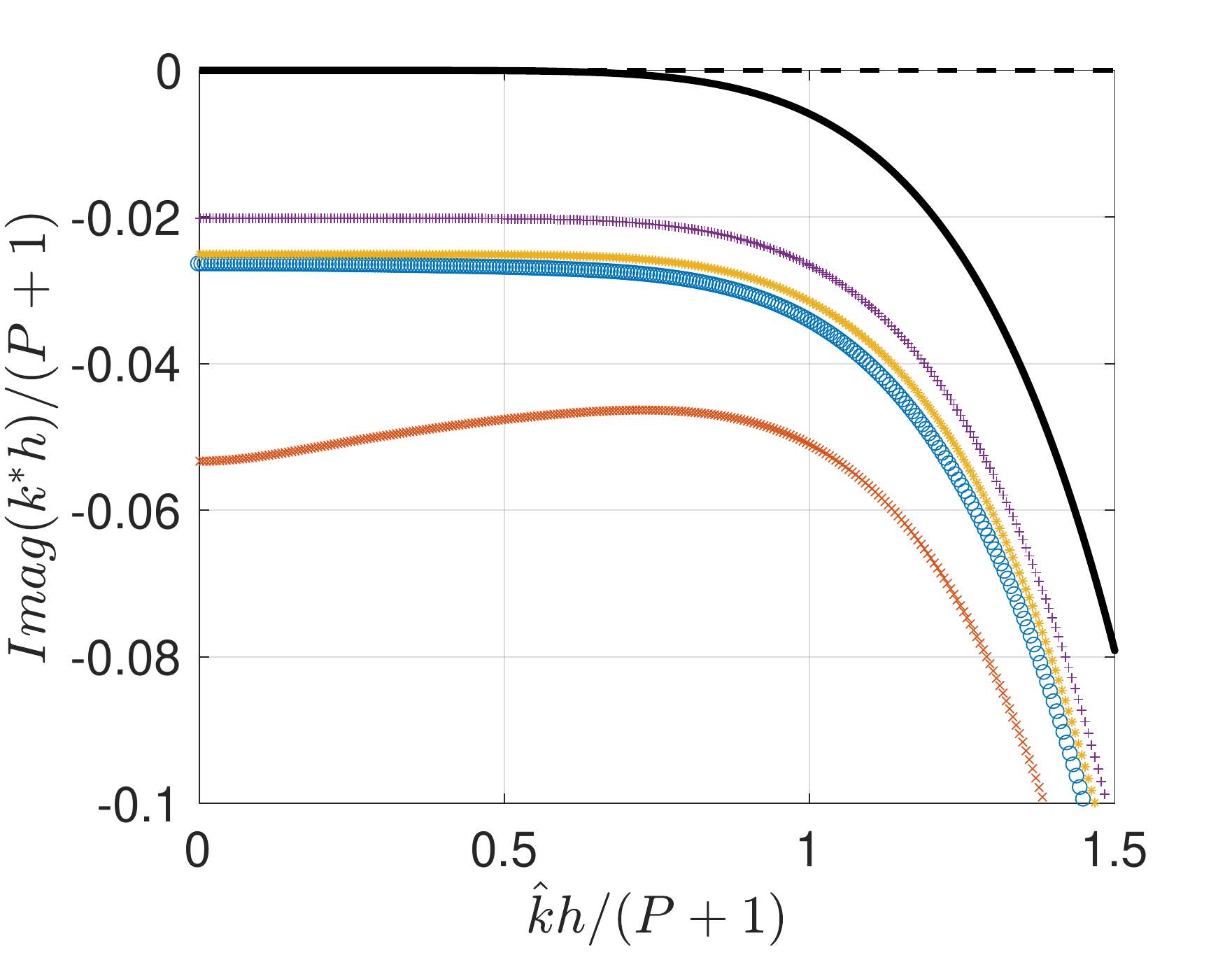}
		\caption{}
	\end{subfigure}
	\begin{subfigure}{.24\textwidth}
		\includegraphics[width=120pt]{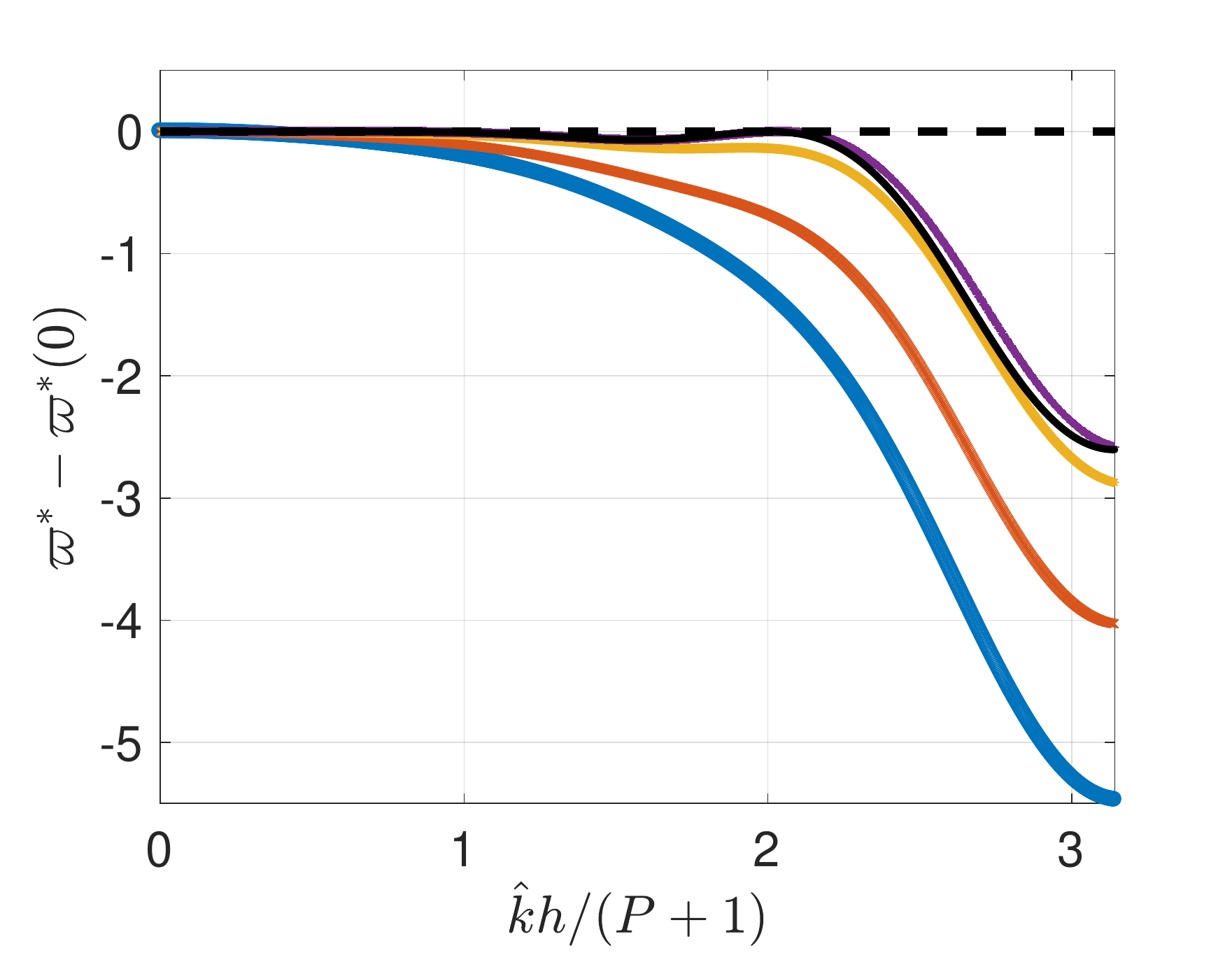}
		\caption{}
	\end{subfigure}
	\centering
	\caption{Eigensolution and non-modal analysis of advection equation with IBM ($N = 40$, $r = 1/40$, $P=3$, $\IBMparam = 1 \times 10^{-3}$) and different diffusivity $\IBMparam_v$. a) Dispersion. b) Dissipation subtracted by the IBM-induced dissipation $\gamma_{IBM}$ (dissipation at $k=0$). c) Zoom-in dissipation. d) Short-term dissipation subtracted by the dissipation at $k=0$.}
	\label{fig:exp2-analysis}
\end{figure*}

From the eigensolution analysis and the numerical experiments, it can be concluded that the IBM works as a porous medium to provide damping inside the solid, where the damping effect can be evaluated through the eigensolution analysis. Therefore, alternative methods that provide additional damping in the solid region can be considered to develop better IBM schemes, like adding artificial viscosity or other wave absorption approaches. Note that all of these strategies must be imposed inside the solid to avoid modifying the real physics. In the present study we use a second-order term to provide additional damping to better satisfy Dirichlet boundary conditions. This strategy has also been mentioned in \citep{brown2014CBVP}, where the constitutive equations are removed from the solid region and the nonphysical diffusion term is introduced to ensure the smoothness and continuity of the solution. This is different from the present work since the constitutive equation still exists for the solid body. This results in the following equation
\begin{equation}
    \frac{\partial \velx}{\partial t} + \advcoef \frac{\partial \velx}{\partial x} = - \frac{\mask}{\IBMparam}\velx +  \mask\IBMparam_v \frac{\partial^2 \velx}{\partial x^2}, 
\end{equation}
where $\IBMparam_v$ is the diffusivity coefficient. We focus on two types of discretization with $N = 40$ and $N = 80$ elements, where the solid domain is set to be the same where $\Delta = 0.05$, indicating $\Delta = h$ and $\Delta = 2h$ for each case, respectively. The polynomial order and penalization parameter are set to $P = 3$ and $\IBMparam = 1 \times 10^{-3}$. The viscosity parameter $\IBMparam_v$ is varied to study the effect of additional dissipation. This viscous term is discretized by the Local DG (LDG) scheme \citep{cockburn1998local}, where two parameters in the scheme are set to $\beta = 0.5$ and $\tau = 0.1$. Eigensolution analysis with different $\IBMparam_v$ is firstly performed for the case with $N = 40$, as shown in Figure \ref{fig:exp2-analysis}. Firstly, it indicates that as the viscous term is added, the dispersion behavior is only affected in high wavenumber region. Secondly, from Figure \ref{fig:exp2-analysis}c, adding the second-order term will lead to smaller dissipation than the case with only volume penalization. This indicates the IBM treatment can be further enhanced by adding artificial viscosity inside the solid, making the solution further converge to the expected value. Thirdly, it should also be noted that although increasing $\IBMparam_v$ will lead to larger total dissipation (as shown by short-term diffusion in Figure \ref{fig:exp2-analysis}d), the dissipation of the physical mode does not increasing accordingly. This indicates that an optimal $\IBMparam_v$ may exist which gives the largest diffusion to the physical mode (thus the solution).

\begin{figure*}[htbp]
\begin{subfigure}{.3\textwidth}
		\includegraphics[width=120pt]{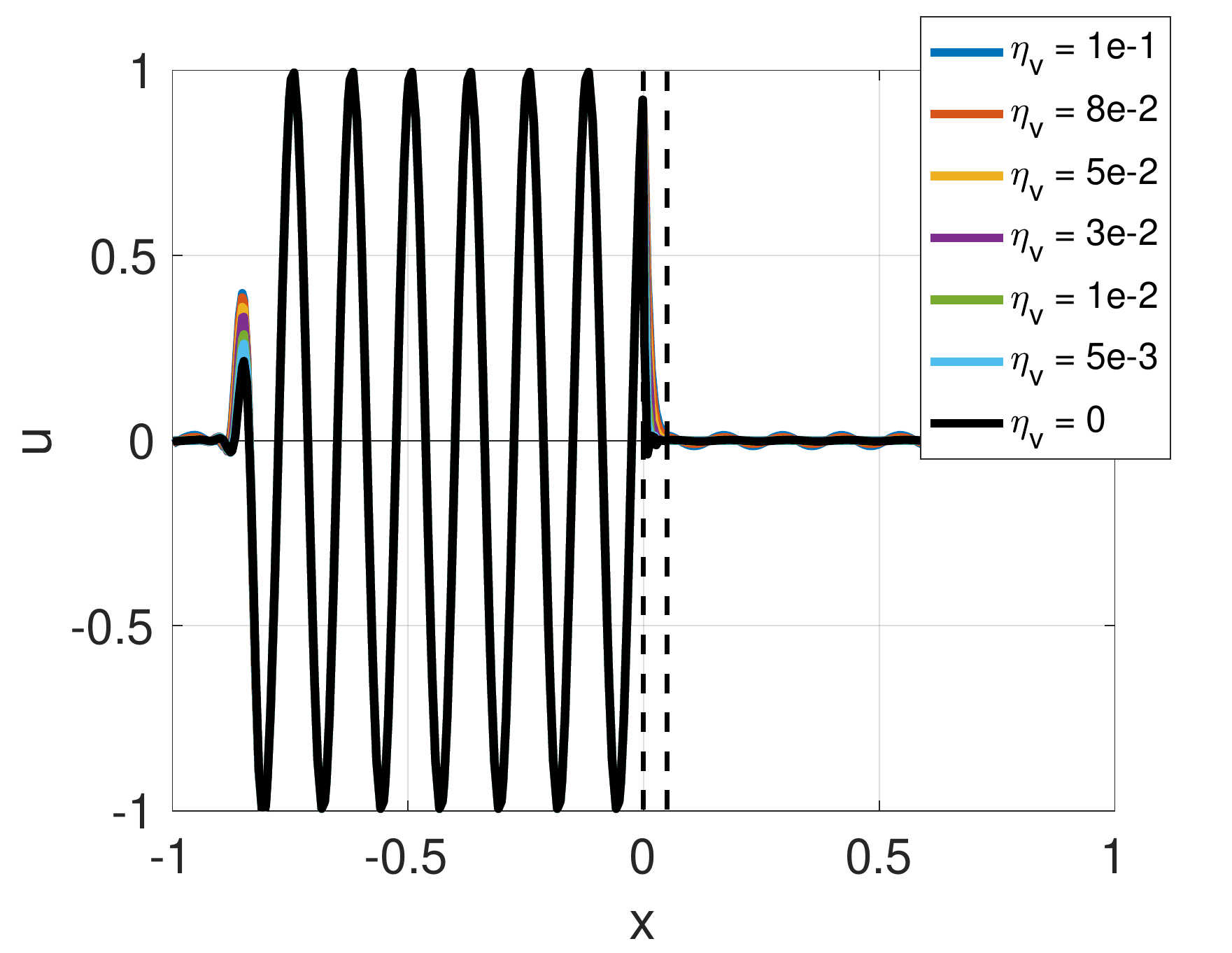}
		\caption{}
	\end{subfigure}
	\begin{subfigure}{.3\textwidth}
		\includegraphics[width=120pt]{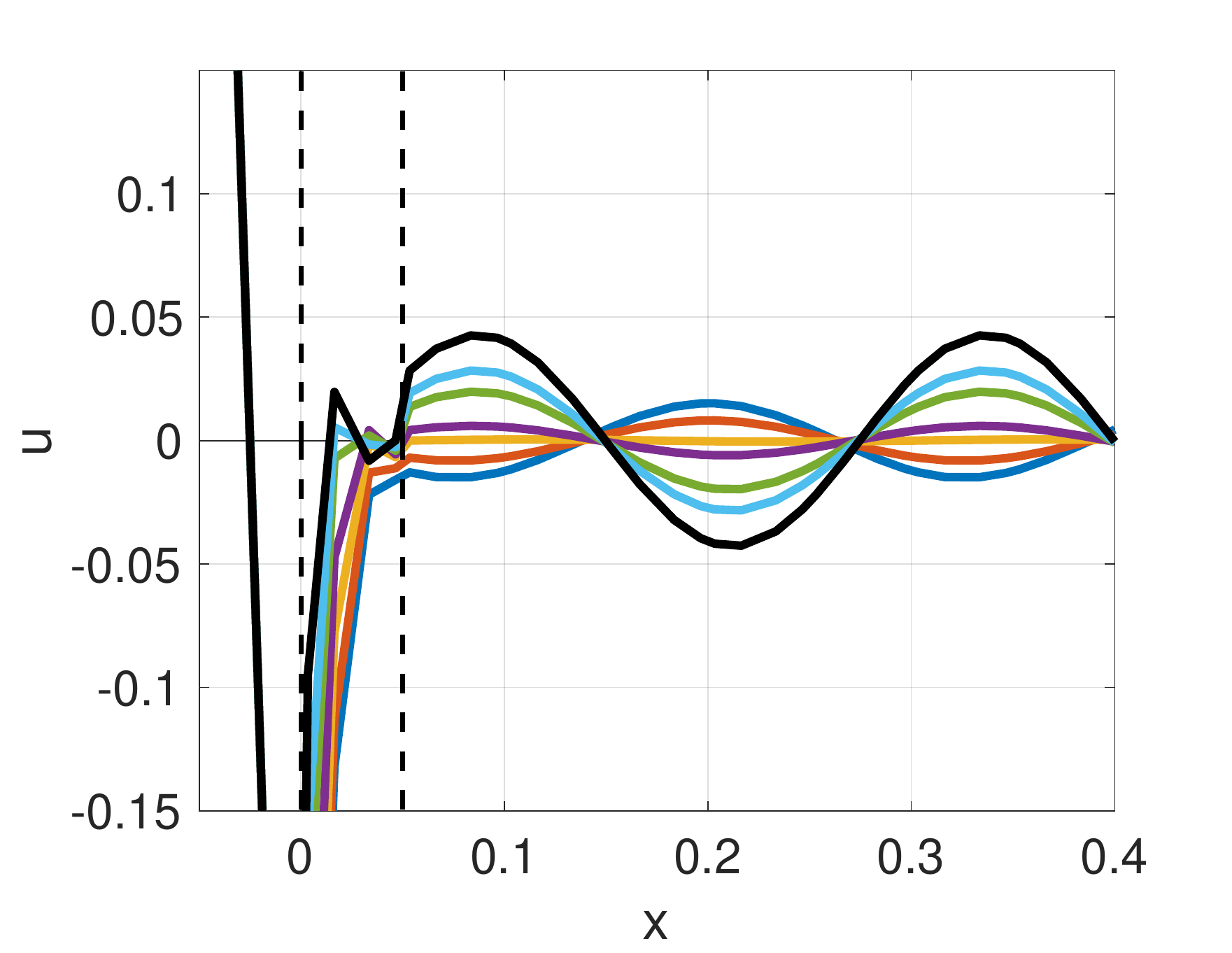}
		\caption{}
	\end{subfigure}
	\centering
	\caption{Simulation of advection equation with IBM under different diffusivity coefficients ($N = 40$, $r = 1/40$, $P = 3$, $\IBMparam = 1 \times 10^{-3}$, initial wavenumber $kh/(P+1) = 0.3142$). a) Global view. b) Zoom-in view.}
	\label{fig:exp2-sim1}
\end{figure*}

\begin{figure*}[htbp]
\begin{subfigure}{.3\textwidth}
		\includegraphics[width=120pt]{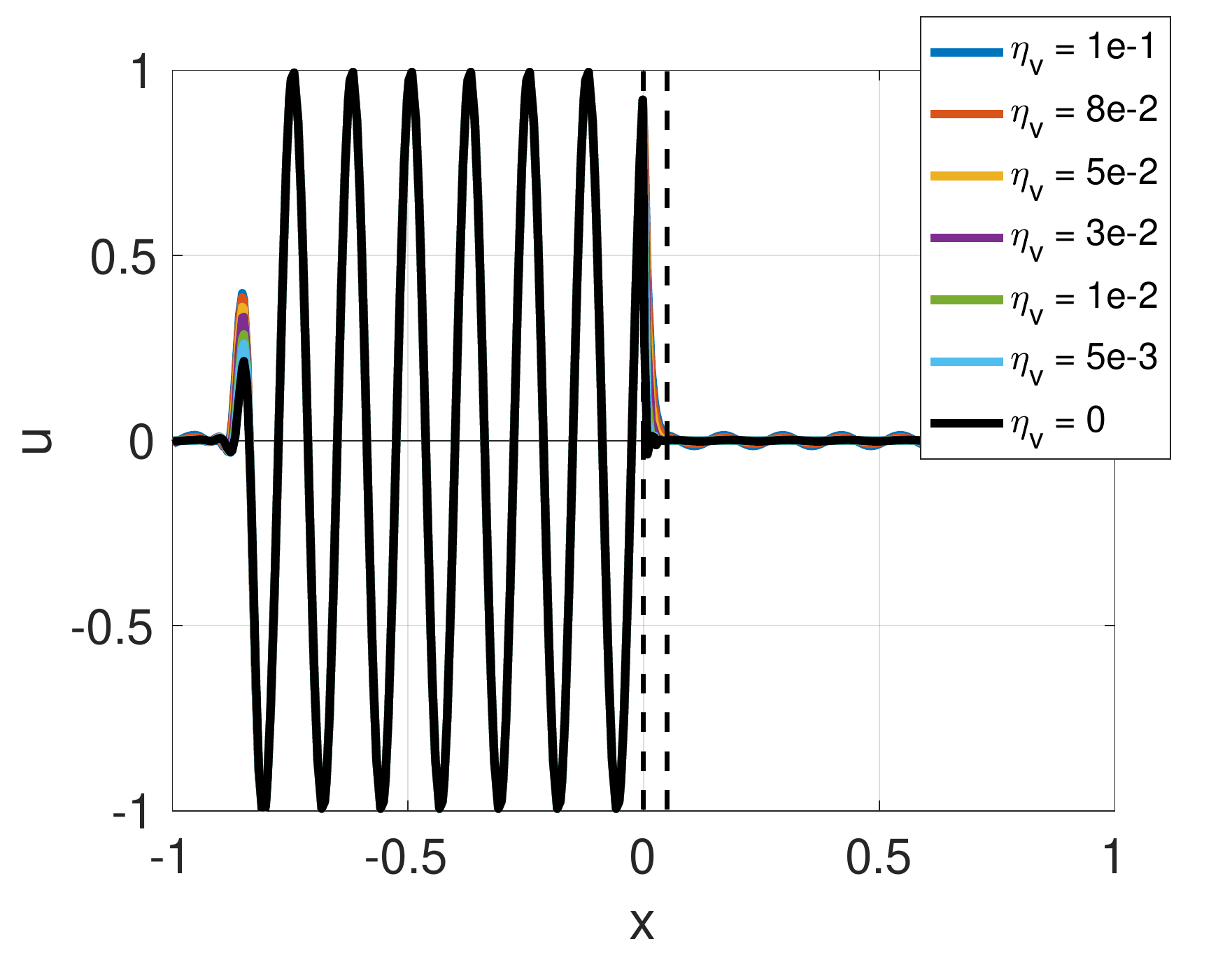}
		\caption{}
	\end{subfigure}
	\begin{subfigure}{.3\textwidth}
		\includegraphics[width=120pt]{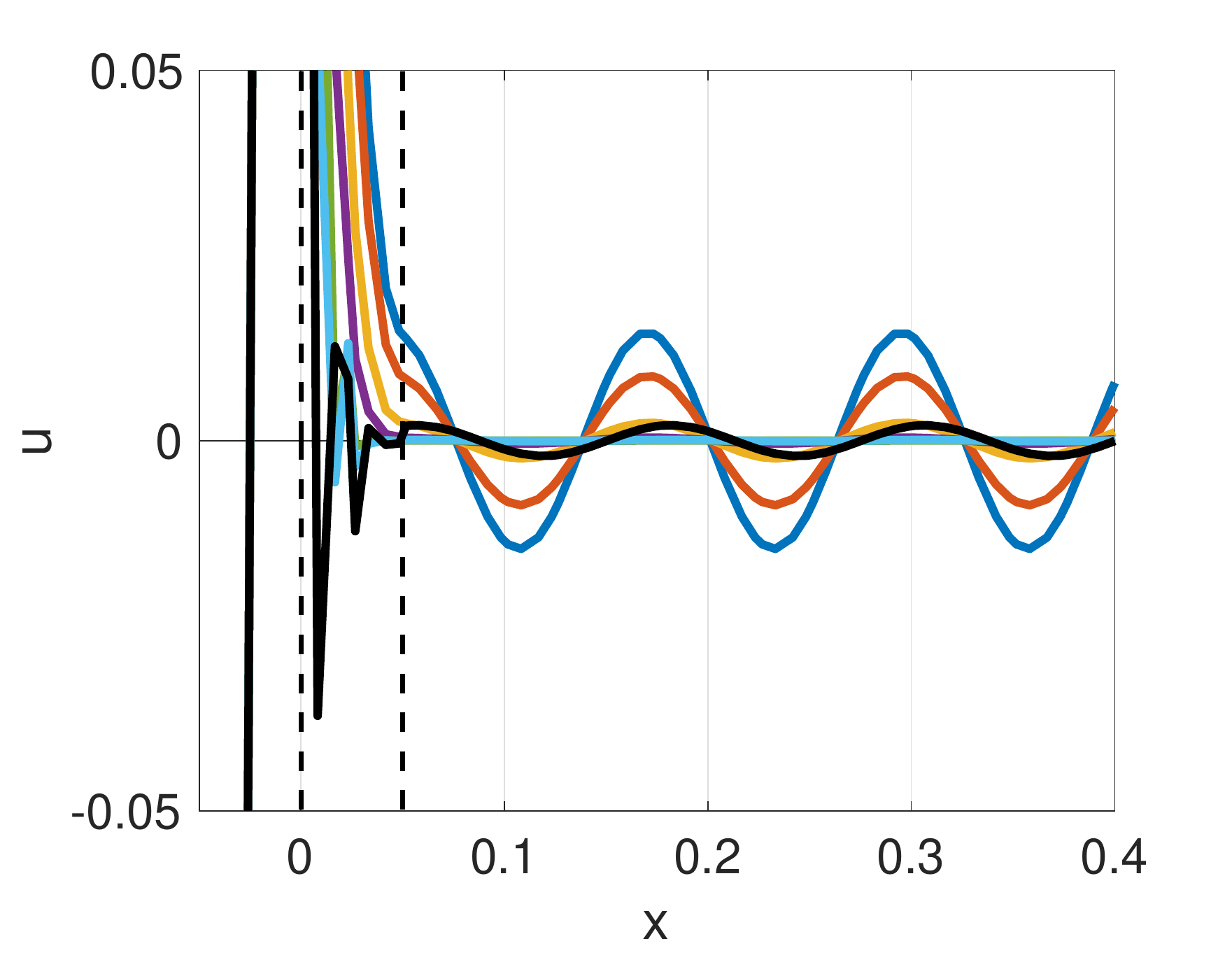}
		\caption{}
	\end{subfigure}
	\centering
	\caption{Simulation of advection equation with IBM under different diffusivity coefficients ($N = 80$, $r = 1/40$, $P = 3$, $\IBMparam = 1 \times 10^{-3}$, initial wavenumber $kh/(P+1) = 0.3142$). a) Global view. b) Zoom-in view.}
	\label{fig:exp2-sim2}
\end{figure*}

After the eigensolution analysis, we perform numerical simulation for both test cases under different $\IBMparam_v$, where the final time is set to $1.1$ and the initial condition is a sinusoidal wave with wavenumber $kh/(P+1) = 0.3142$. Results for both cases are shown in Figure \ref{fig:exp2-sim1} and Figure \ref{fig:exp2-sim2}, respectively. In Figure \ref{fig:exp2-sim1}, it is observed that for $N = 40$, when $\IBMparam_v$ is increasing, the solution is firstly damped until about $\IBMparam_v = 5 \times 10^{-2}$, where the performance is nearly optimal and the solution is quite close to zero. However, if $\IBMparam_v$ continues to increase, the solution starts to grow with another phase. This is reflected by the eigensolution as well, where the wavenumber should be shifted by $\pi/N$ before scaling. This means that for the present case, there is no need to increase $\IBMparam_v$ when $\IBMparam_v$ reach $5 \times 10^{-2}$. For the second group of results shown in Figure \ref{fig:exp2-sim2}, it can be seen that when very small $\IBMparam_v$ is set, the solution can be improved, but further increasing $\IBMparam_v$ will make the solution less accurate, e.g., when $\IBMparam_v > 3 \times 10^{-2}$. These results indicate that there exists an optimal combination of $\IBMparam$ and $\IBMparam_v$ that gives the best numerical accuracy. These optima will be explored in what follows.

\begin{figure*}[htbp]
\begin{subfigure}{.4\textwidth}
		\includegraphics[width=160pt]{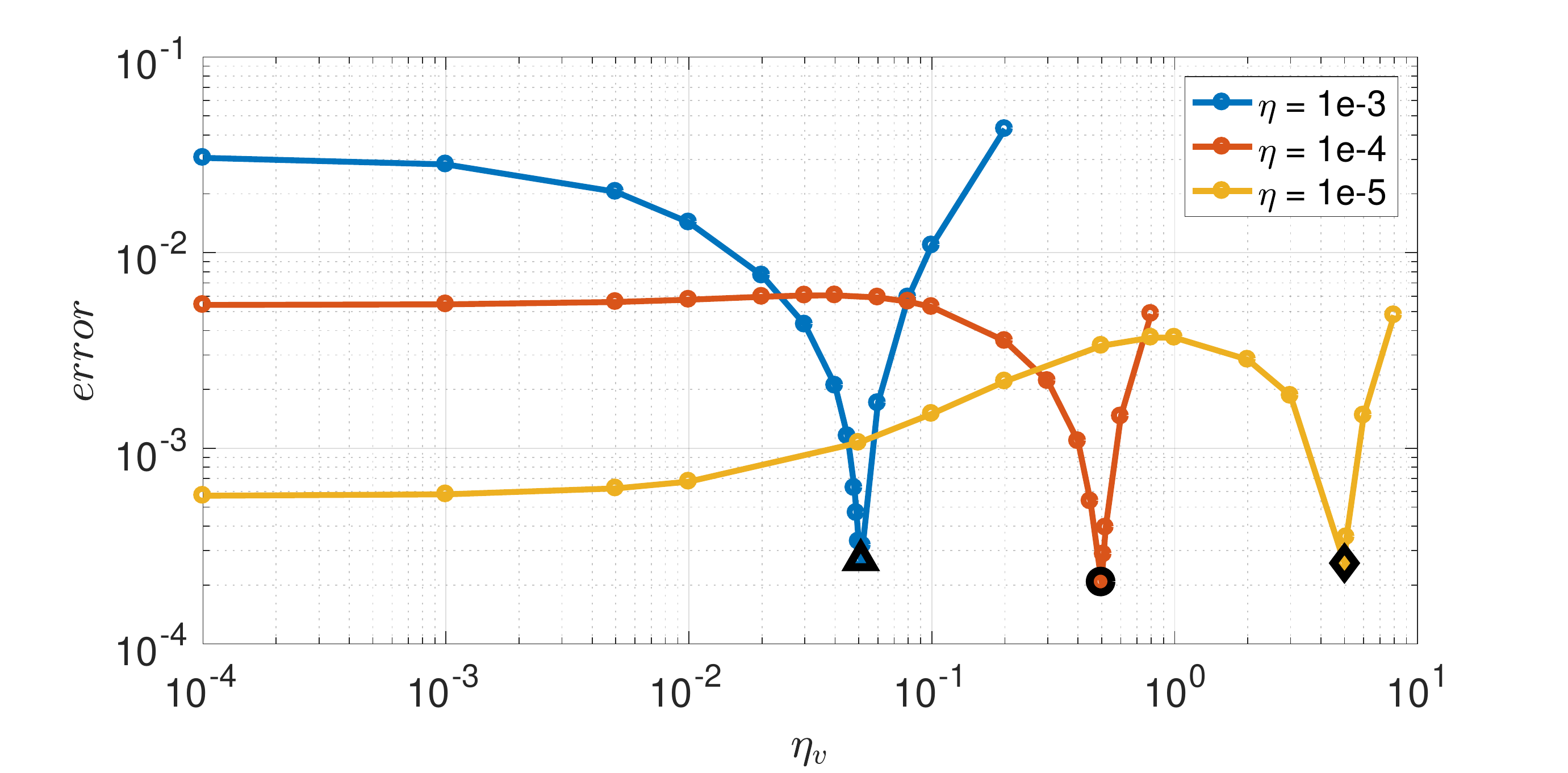}
		\caption{}
	\end{subfigure}
	\begin{subfigure}{.4\textwidth}
		\includegraphics[width=160pt]{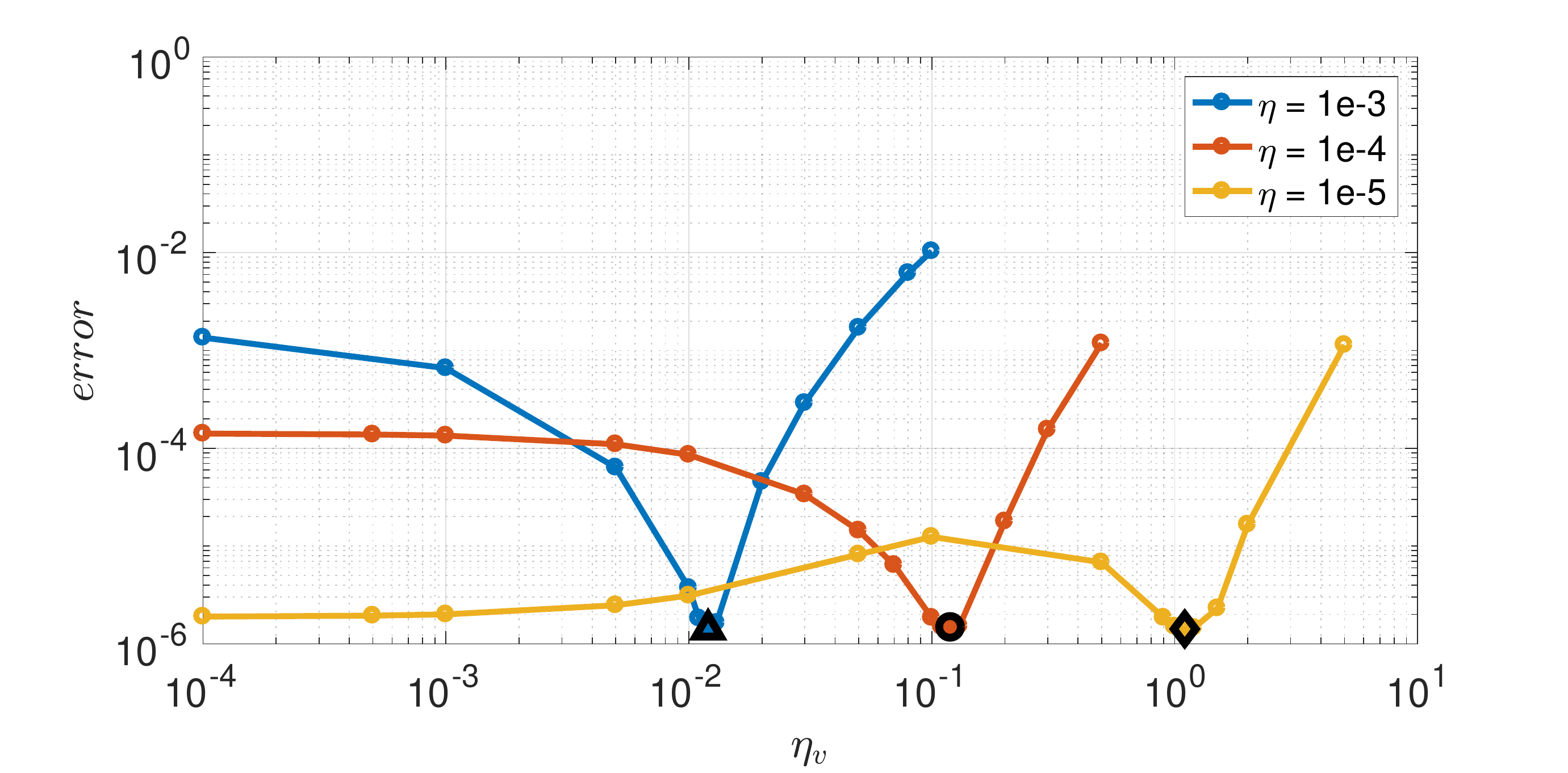}
		\caption{}
	\end{subfigure}
	\centering
	\caption{Simulation error versus diffusivity coefficient ($r = 1/40$, $P = 3$, initial wavenumber $kh/(P+1) = 0.3142$). a) $N = 40$. b) $N = 80$. Black symbols refer to the lowest error with respect to a given penalization parameter $\IBMparam$.}
	\label{fig:exp2-error}
\end{figure*}

\begin{figure*}[htbp]
\begin{subfigure}{.4\textwidth}
		\includegraphics[width=160pt]{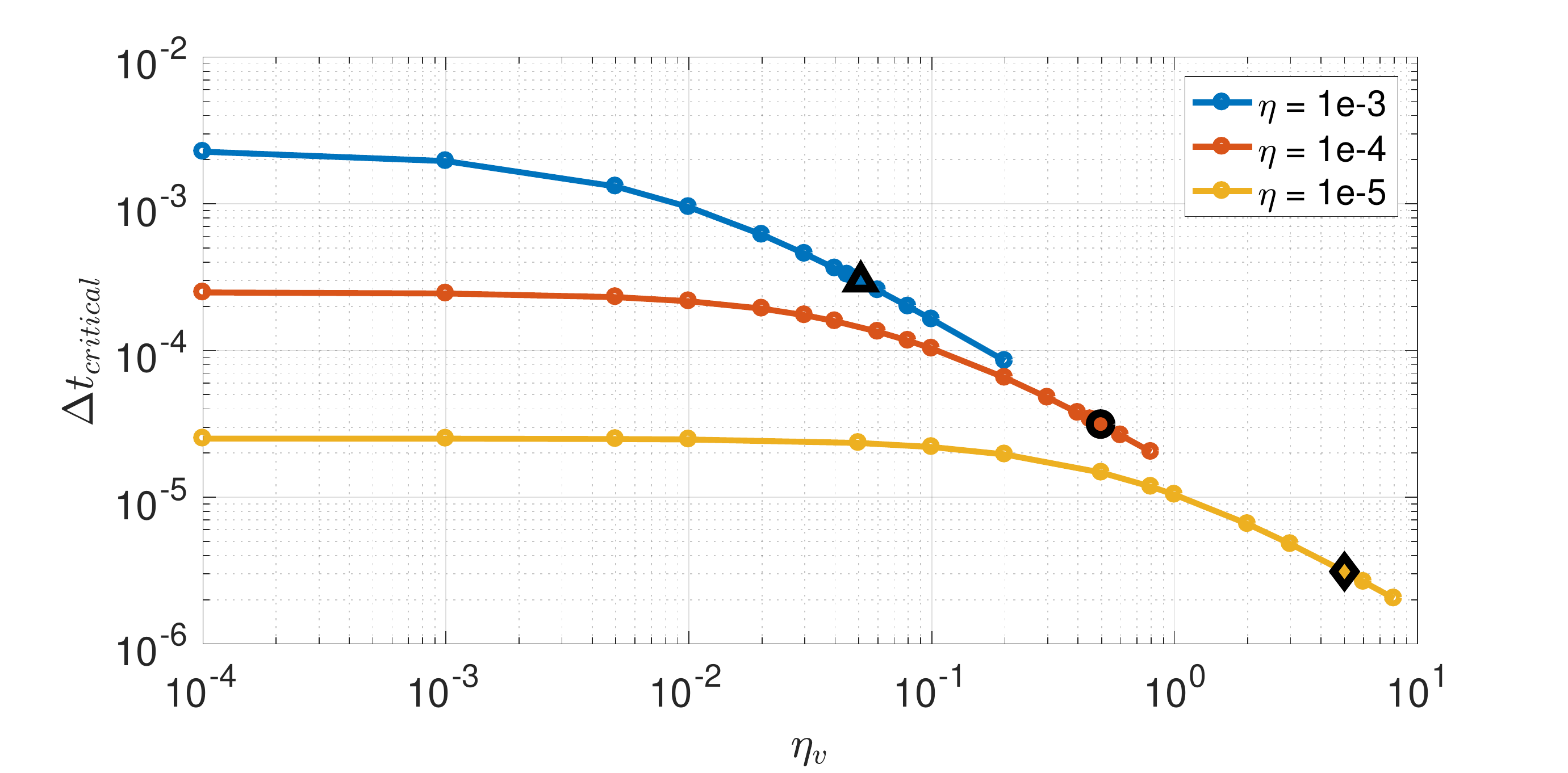}
		\caption{}
	\end{subfigure}
	\begin{subfigure}{.4\textwidth}
		\includegraphics[width=160pt]{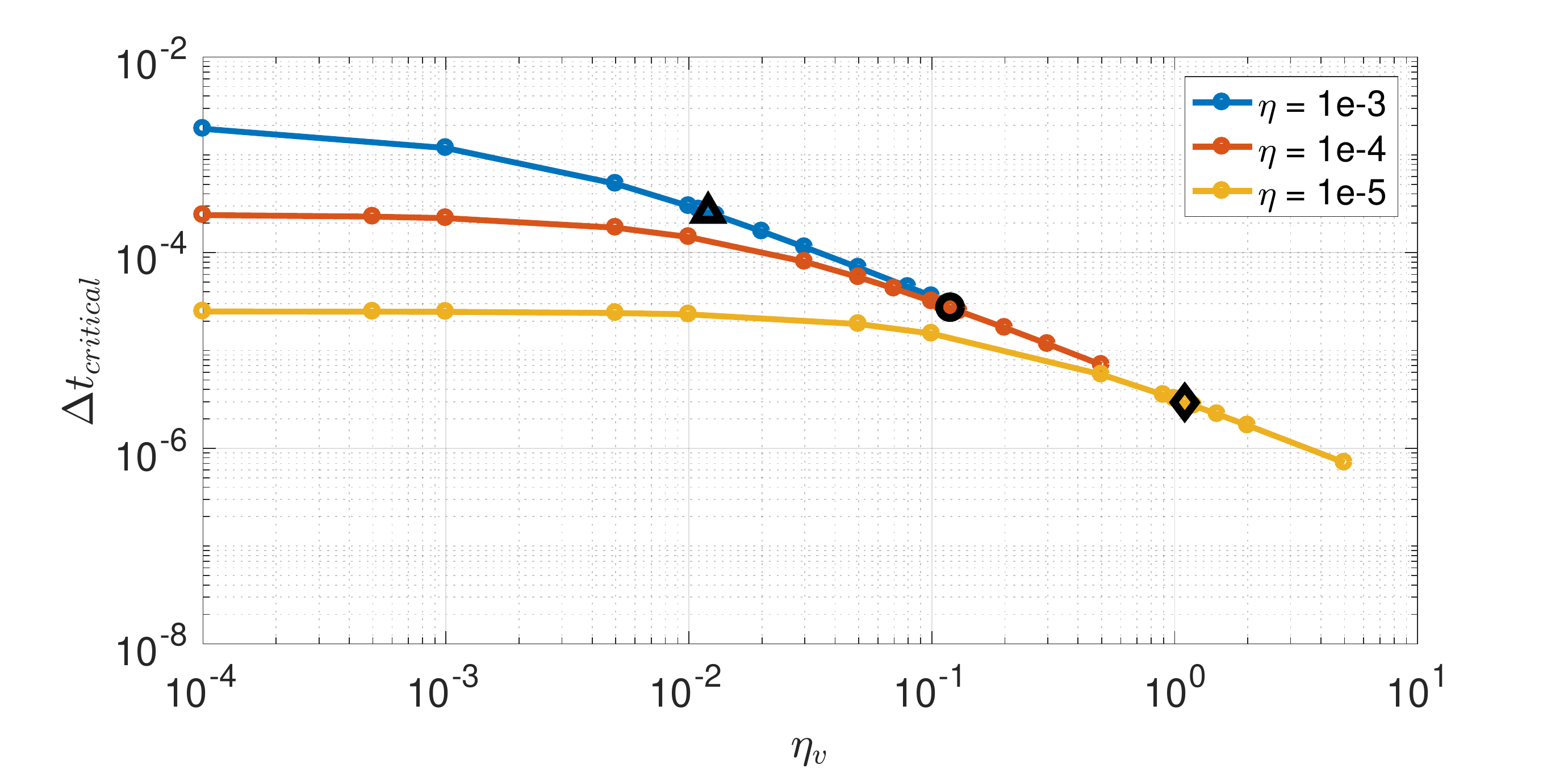}
		\caption{}
	\end{subfigure}
	\centering
	\caption{Critical time step versus diffusivity coefficient ($r = 1/40$, $P = 3$, initial wavenumber $kh/(P+1) = 0.3142$). a) $N = 40$. b) $N = 80$. Black symbols refer to the lowest error with respect to a given penalization parameter $\IBMparam$.}
	\label{fig:exp2-Dt}
\end{figure*}

To further investigate the relationship between accuracy and viscosity parameter, the errors are compared. This error is defined as the error between the solution in computational domain from the right solid boundary to the right boundary of the whole computational domain ($x \in [\Delta, 1]$) and the penalized value $\velx_s = 0$. Defining the number of solution points inside the domain of interest as $L$, we have the mean squared error
\begin{equation}
    error = \sqrt{\frac{\sum_{i=1}^{L} [\velx(x_i) - \velx^{exact}(x_i)]^2}{L}} \ , \ x_i \in [\Delta, 1], \ \velx^{exact} = 0.
\end{equation}

The simulation error at different $\IBMparam$ and $\IBMparam_v$ is compared in Figure \ref{fig:exp2-error}. It is clear that adding the second-order term will lead to lower error when the diffusivity parameter $\IBMparam_v$ is chosen properly. For a given penalization parameter and problem setting, an optimal $\IBMparam_v$ exists to provide the lowest error. However, a big gain is only seen for $\IBMparam = 1 \times 10^{-3}$ and $\IBMparam = 1 \times 10^{-4}$, where the optimal $\IBMparam_v$ will lead to a large reduction of the error. When the penalization parameter is sufficiently small, e.g., $\IBMparam = 1 \times 10^{-5}$, the optimal performance requires a relatively larger $\IBMparam_v$ and the gain becomes marginal. In addition, the critical time step $\Delta t$ for different $\IBMparam_v$ is compared in Figure \ref{fig:exp2-Dt}. This time step is the largest $\Delta t$ allowed for a stable time integration based on the present third-order Runge-Kutta time-marching scheme. It can be concluded that by adding the second-order term, less restrictions on time step are required for similar accuracy, which is a main advantage of this strategy. For example, for the first test case $N = 40$, the optimal accuracy for $\IBMparam = 1 \times 10^{-4}$ and $\IBMparam = 1 \times 10^{-5}$ is almost the same, while the time step can be about ten times larger ($3.1 \times 10^{-5}$ against $3.1 \times 10^{-6}$). In general, the strategy proposed here is useful both in terms of improving the accuracy and relaxing the time step restriction which arises from the stiff source term associated to the volume penalization.

\begin{figure*}[htbp]
\begin{subfigure}{.4\textwidth}
		\includegraphics[width=160pt]{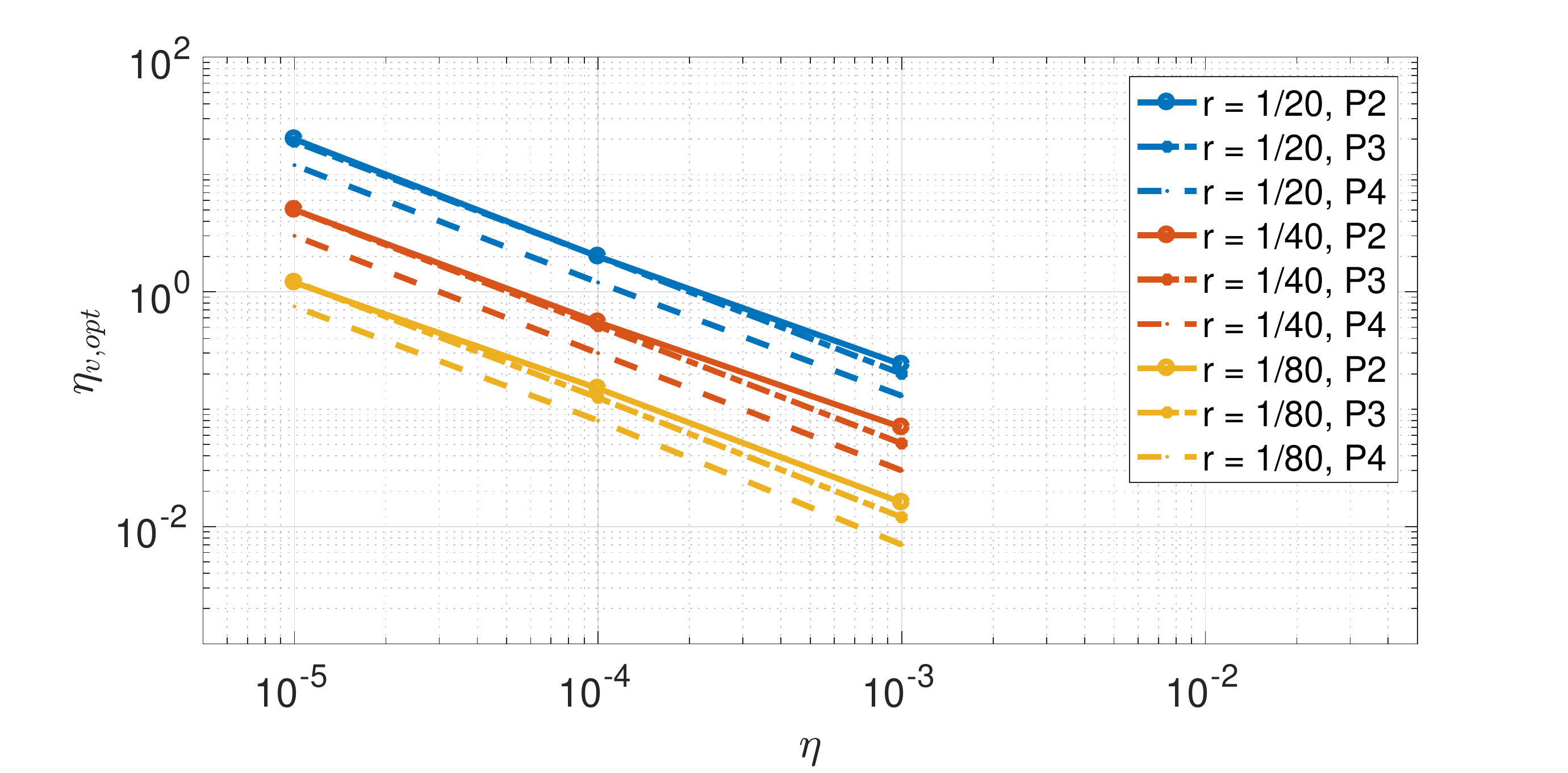}
		\caption{}
	\end{subfigure}
	\begin{subfigure}{.4\textwidth}
		\includegraphics[width=160pt]{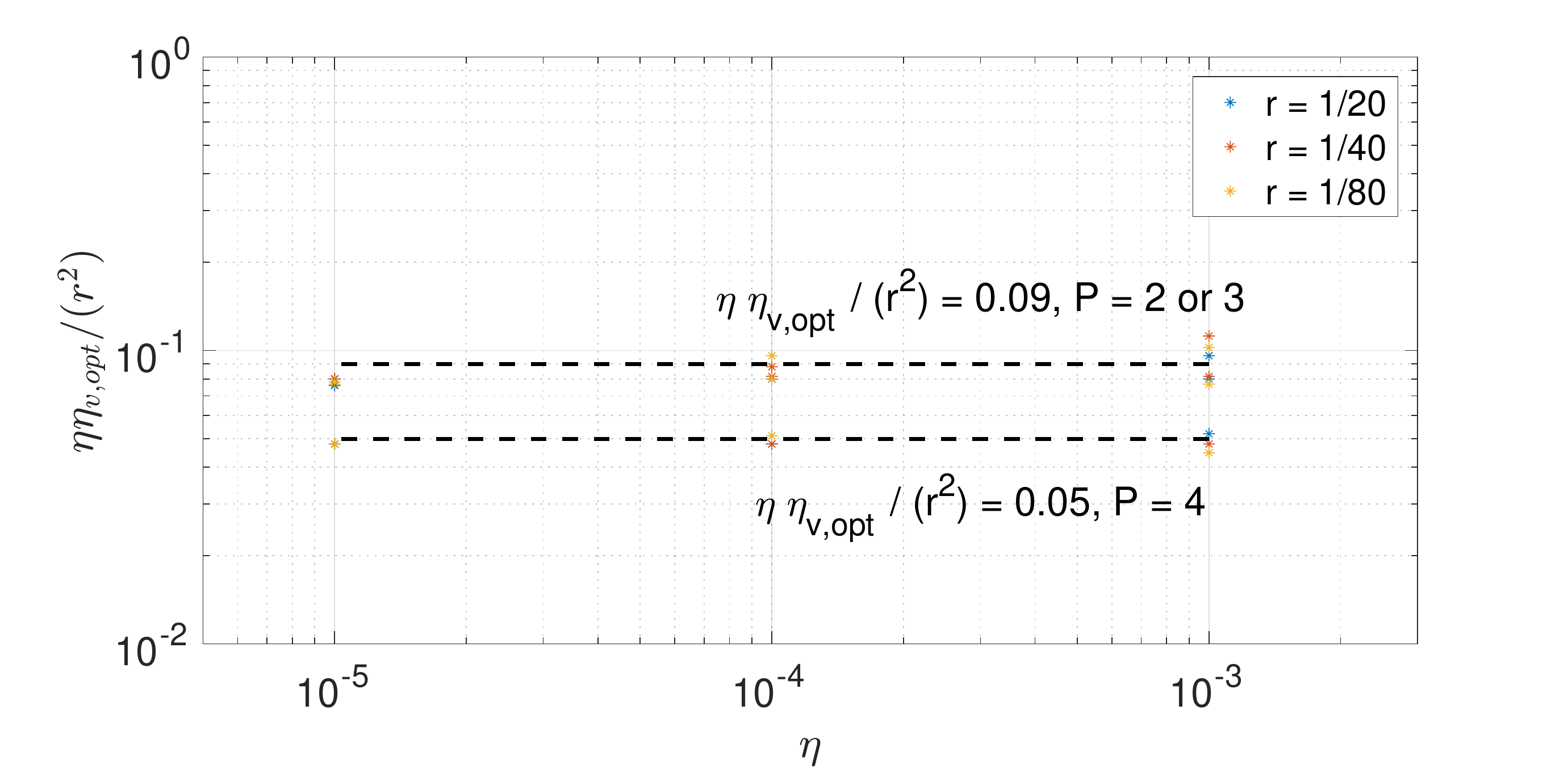}
		\caption{}
	\end{subfigure}
	\centering
	\caption{Relation between optimal $\IBMparam_{v}$ (defined as $\IBMparam_{v, opt}$) and $\IBMparam$ for various cases with different solid ratio $r$ (one solid element) and polynomial order $P$. a) $\IBMparam_{v, opt}$ versus $\IBMparam$. b) Scaling relationship between $\IBMparam_{v, opt}$ and $\IBMparam$.}
	\label{fig:eta-relation}
\end{figure*}

Finally, the relationship between optimal $\IBMparam_v$ (defined as $\IBMparam_{v,opt}$) and $\IBMparam$ is explored. From Figure \ref{fig:exp2-error}, the optimal $\IBMparam_v$ (indicated by black symbols) increases as $\IBMparam$ is decreasing, indicating that there exists a dependency between these parameters. To further study this behavior, we consider simulations for different penalization parameters, polynomial orders and solid ratios, as shown in Figure \ref{fig:eta-relation}a. We fix one solid element in the middle and choose three ratios $r$ including $1/20$, $1/40$ and $1/80$ (representing different solid widths), three $P$ ranging from $2$ to $4$ and three penalties $\IBMparam$ including $1 \times 10^{-3}$,  $1 \times 10^{-4}$ and $1 \times 10^{-5}$. For each parameter pair, an initial condition with a low wavenumber is chosen to represent the resolved wavenumber region of the scheme. Numerical simulations across a group of $\IBMparam_v$ are performed to numerically find $\IBMparam_{v,opt}$. From Figure \ref{fig:eta-relation}a, a nearly linear relation is seen between $\IBMparam_{v,opt}$ and $\IBMparam$. For a given $r$ and $\IBMparam$, $\IBMparam_{v,opt}$ remains approximately constant for polynomial order $2$ and $3$, while it reduces by half for polynomial order $4$. Therefore, scaling relationships of $\IBMparam \IBMparam_{v,opt} / (r^2) = 0.09$ and $\IBMparam \IBMparam_{v,opt} / (r^2) = 0.05$ can be fitted for $P = 2 / 3$ and $P = 4$, respectively. These fitted scaling curves are shown in Figure \ref{fig:eta-relation} b. The fitted scaling relation can provide guidelines for the selection of $\IBMparam_{v,opt}$ based on $\IBMparam$ and $P$.

\subsection{Flow past a circular cylinder}
To further investigate the observations from the proposed analyses, volume penalization for Navier-Stokes equations based on high-order FR discretization is considered \citep{kou2021IBMFR1,kou2021IBMFR2}. As a benchmark test case, steady flow past a static circular cylinder at Reynolds number 40 is chosen. For compressible Navier-Stokes equations, the following source term is added to the right-hand-side of the momentum equations \citep{abgrall2014IBM}:

\begin{equation}
\bm{S}_1 = \frac{\chi}{\eta} \times \begin{pmatrix}
\rho u_{s}-\rho u\\ 
\rho v_{s}-\rho v
\end{pmatrix} + \rho \chi \IBMparam_v  \times \begin{pmatrix}
 \frac{\partial^2 u}{\partial x^2} + \frac{\partial^2 u}{\partial x \partial y}\\ 
\frac{\partial^2 v}{\partial x^2} + \frac{\partial^2 v}{\partial y \partial x}
\end{pmatrix},
\label{eq:IBMsource-1}
\end{equation}
where $\rho$, $u$ and $v$ denote the density and velocities in $x$ and $y$ directions, respectively. $u_s$ and $v_s$ refer to the penalized velocities inside the solid, which are fixed to zero in the present study to mimic the no-slip wall boundary condition of a cylinder. An additional source term added to the energy equation is also needed following \citep{abgrall2014IBM}

\begin{equation}
\bm{S}_2 = \frac{\chi}{\eta} \times \begin{pmatrix}
\frac{\rho}{2}(u_{s}^2+v_{s}^2)-\frac{\rho}{2}(u^2+v^2)
\end{pmatrix}.
\label{eq:IBMsource-2}
\end{equation}

More details about the implementation, e.g., surface geometry representation and computation of integral aerodynamic coefficients, are given in \citep{kou2021IBMFR1,kou2021IBMFR2}. To apply the scaling law proposed in the previous section, the solid ratio $r$ should be defined for two-dimensional flows. Here we generalize the solid ratio for two-dimensional flows as 

\begin{equation}
r = \sqrt{\frac{S_{solid}}{S_{domain}}},
\label{eq:r-2d}
\end{equation}
where $S_{solid}$ and $S_{domain}$ are areas of the solid region and the whole computational domain, respectively. 

\begin{figure*}[htbp]
\begin{subfigure}{.4\textwidth}
		\includegraphics[width=160pt]{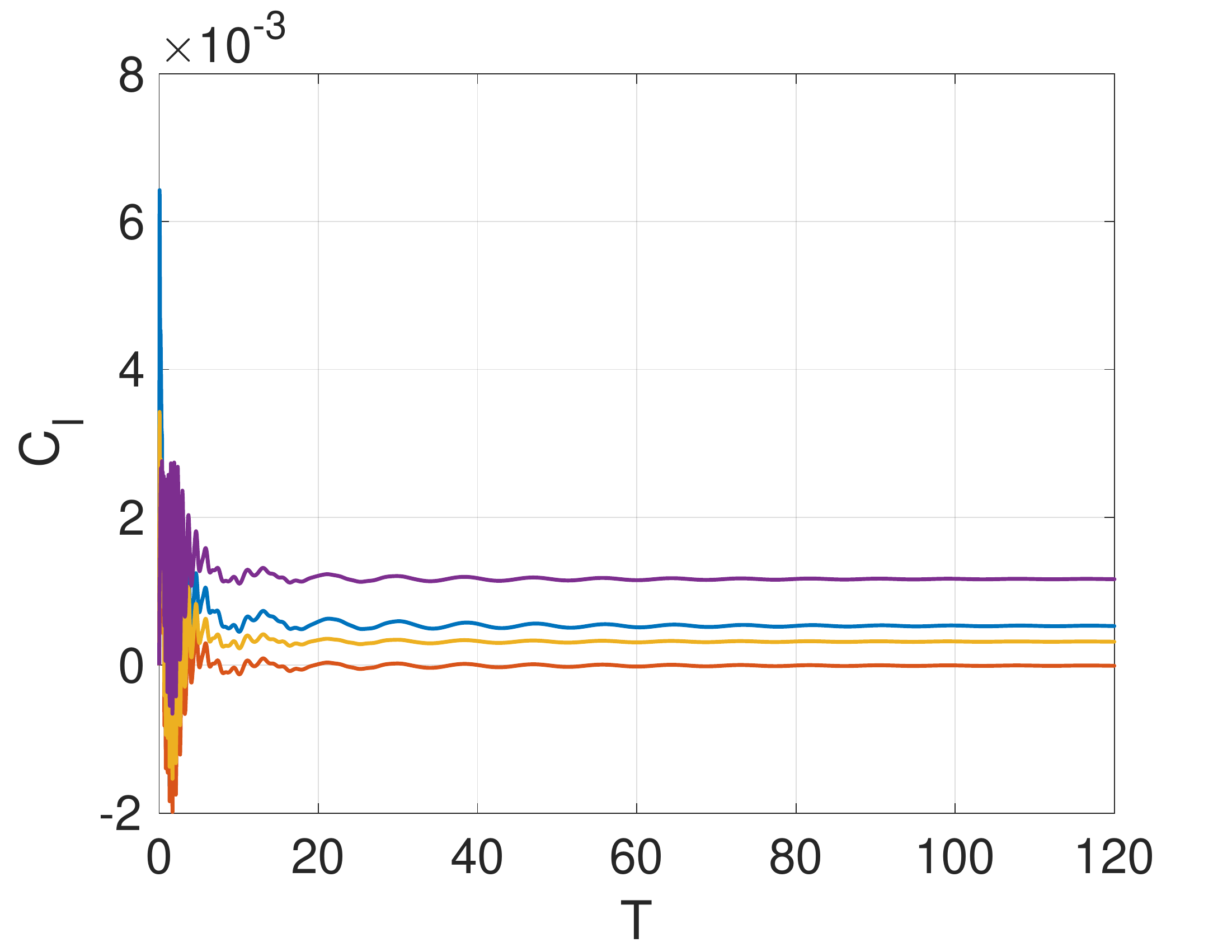}
		\caption{}
	\end{subfigure}
	\begin{subfigure}{.4\textwidth}
		\includegraphics[width=160pt]{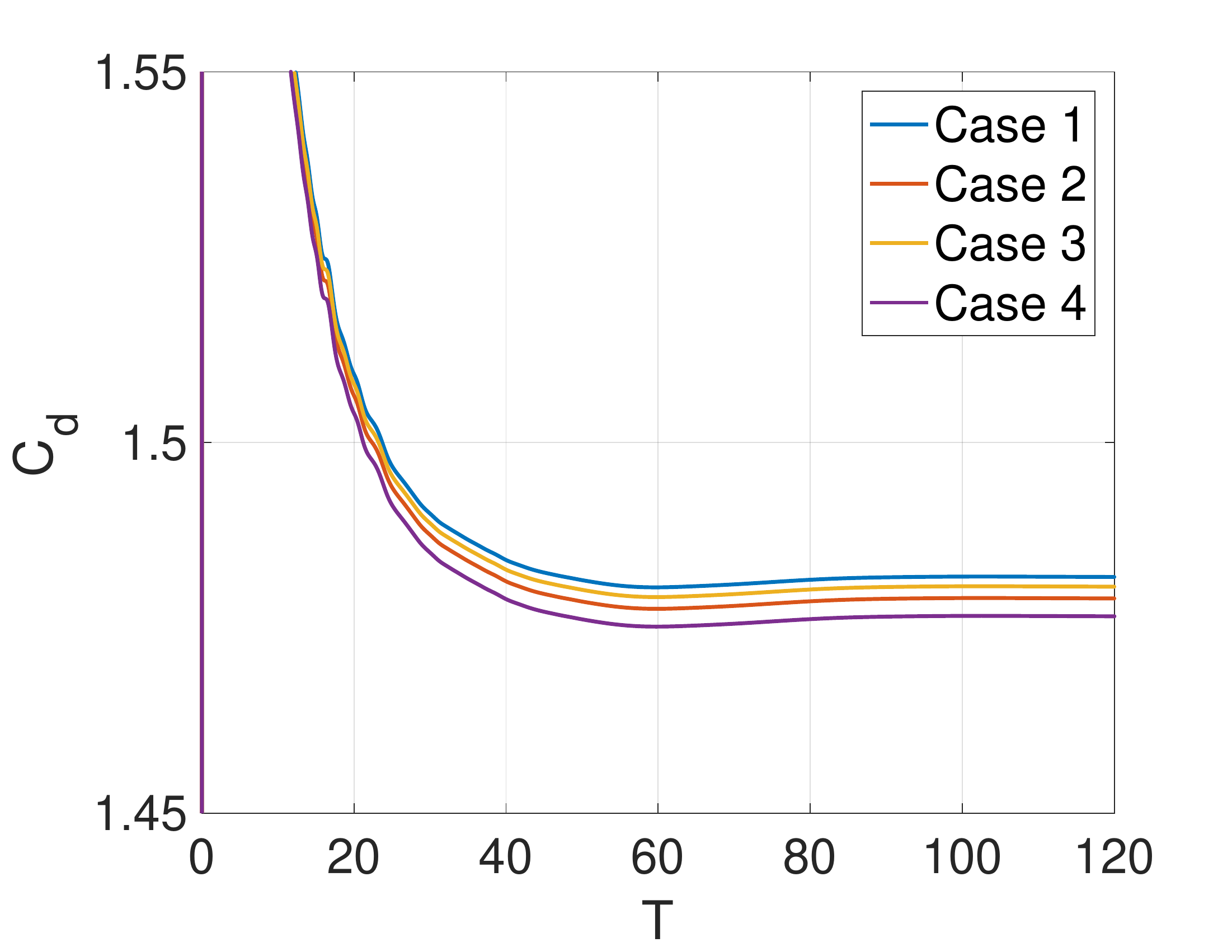}
		\caption{}
	\end{subfigure}
	\centering
	\caption{Time evolution of aerodynamic coefficients for flow past a cylinder at $Re = 40$. Case 1: $\IBMparam = 2 \times 10^{-4}$, $\IBMparam_v = 0$. Case 2: $\IBMparam = 1 \times 10^{-3}$, $\IBMparam_v = 1.5\times 10^{-2}$ (predicted by the relation $\IBMparam \IBMparam_{v,opt} / (r^2) = 0.09$). Case 3: $\IBMparam = 1 \times 10^{-3}$, $\IBMparam_v = 1 \times 10^{-2}$. Case 4: $\IBMparam = 1 \times 10^{-3}$, $\IBMparam_v = 2 \times 10^{-2}$. a) Lift coefficient. b) Drag coefficient.}
	\label{fig:cylinder}
\end{figure*}

A structured Cartesian mesh with local refinement region near the wall is considered for the present numerical experiment. The size of rectangular computational domain is $x \in [-30D, 50D]$ and $y \in [-30D, 30D]$, where $D$ is the diameter of the cylinder. In the square region $x \in [-0.6D, 0.6D]$ and $y \in [-D, D]$, uniform grid with size $0.03D$ is chosen. This results in the total number of elements $89 \times 86$. Explicit time integration based on the low-storage five-stage fourth-order explicit Runge-Kutta method (LSERK) \citep{carpenter1994fourth,hesthaven2007nodal} is used to march the solution in time. The polynomial order is set to $P = 2$ and the time step $\Delta t$ is varied according different combinations of parameters. At first, the volume penalization is applied, where we set $\eta = \Delta t$ from the traditional guideline \citep{schneider2015immersed,engels2015numerical}, which gives $\eta = \Delta t = 2 \times 10^{-4}$. Based on this case, we find another case where similar time step / computational cost is required but may lead to a smaller error, and set a smaller penalization parameter $\IBMparam = 1 \times 10^{-3}$. From the relationship obtained in the previous section for $P = 2$, i.e., $\IBMparam \IBMparam_{v,opt} / (r^2) = 0.09$, we can obtain $\IBMparam_{v,opt} \approx 1.5 \times 10^{-2}$, where $r \approx 0.0128$ is obtained from Equation \ref{eq:r-2d}. Therefore, three additional simulations with $\IBMparam = 1 \times 10^{-3}$ and $\IBMparam_v = 1.5 \times 10^{-2}$, $1.0 \times 10^{-2}$ and $2.0 \times 10^{-2}$ are included for comparison, where the same time step $\Delta t = 2 \times 10^{-4}$ is used for the first two cases, while $\Delta t = 1.5 \times 10^{-4}$ is needed for the last case to ensure stability. 

In Figure \ref{fig:cylinder}, the temporal evolution of the lift and drag coefficients is compared. Due to the symmetric geometry of the cylinder and the selected Reynolds number, the expected lift coefficient should approach zero. As shown in Figure \ref{fig:cylinder}a, the optimal prediction is given by Case 2. This is achieved by the combination of $\IBMparam$ and $\IBMparam_v$, where $\IBMparam_v$ is given by the approximated relationship. This highlights that it is possible to utilize the proposed relationship to help the selection of penalization and diffusivity parameters. In addition, different combinations of these parameters also have an impact on the predicted drag coefficient, as shown in Figure \ref{fig:cylinder}b. Note that the predicted drag coefficients (around 1.48) agree with the results reproduced from existing literature, e.g., \citep{brown2014CBVP,de2007immersed}.

\begin{figure*}[htbp]
 	\centering
	\includegraphics[width=250pt]{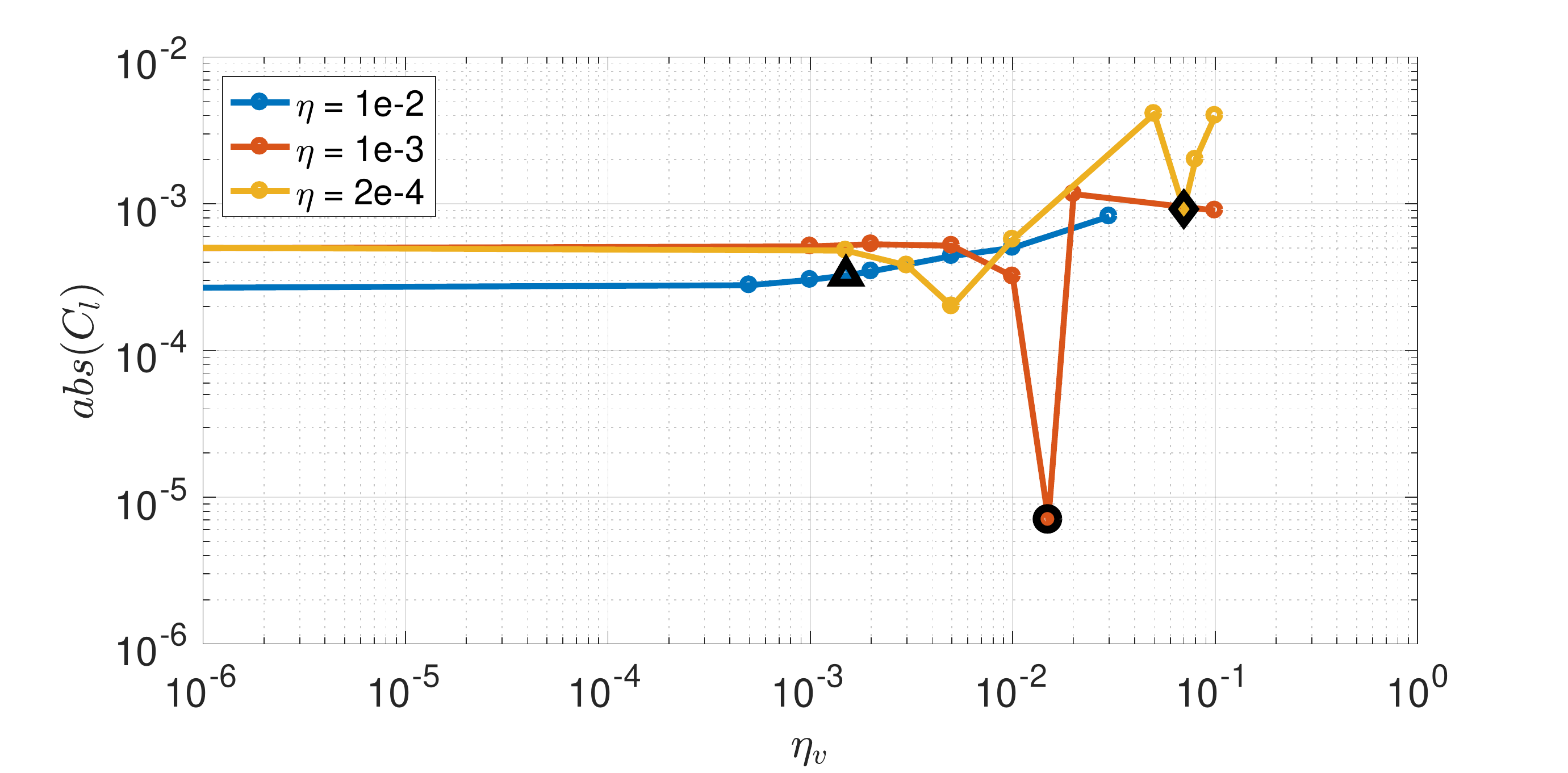}
	\caption{Predicted absolute value of lift coefficient across different diffusivity coefficients. Black symbols represent the optimal $\IBMparam_v$ obtained from the relationship $\IBMparam \IBMparam_{v,opt} / (r^2) = 0.09$.}
	\label{fig:abs-cl}
\end{figure*}

To further evaluate the proposed approach and the approximated relationship, more simulations are conducted, and the predicted lift coefficient at different combinations of penalization and diffusivity parameters are compared in Figure \ref{fig:abs-cl}. Variations of $\IBMparam_v$ at three penalization parameters, including $1 \times 10^{-2}$, $1 \times 10^{-3}$ and $2 \times 10^{-4}$ are considered. It can be seen that for small $\IBMparam$, there exists an optimal $\IBMparam_v$ that gives the best lift coefficient. When $\IBMparam = 1 \times 10^{-3}$, the proposed relationship accurately predicted the $\IBMparam_{v,opt}$, as validated in Figure \ref{fig:cylinder}. At large $\IBMparam$ where the solution is not sufficiently penalized, adding the second-order term will lead to a slightly larger error, while the predicted $\IBMparam_{v,opt}$ still gives reasonable result since a relatively small $\IBMparam_{v} = 1.5 \times 10^{-3}$ is obtained. When $\IBMparam = 2 \times 10^{-4}$, although the optimal value is not well predicted, $\IBMparam_{v,opt} = 7 \times 10^{-2}$ still gives a local minimum error. These results indicate that the proposed relationship has the potential to be used as a guideline for Navier-Stokes simulations, however, further study is still needed. In general, the results in this section show the efficacy of the proposed approach to include second-order term to improve the representation of the boundary condition.

\section{Conclusions}
This paper applies eigensolution and non-modal analyses for the immersed boundary method in the advection equation, where a high-order flux reconstruction scheme based on volume penalization is investigated. Three main conclusions are drawn:

1) The semi-discrete eigensolution analysis and non-modal analysis show that volume penalization does not incur in significant changes in the dispersion behavior of the underlying high-order scheme in the resolved wavenumber range. It imposes the boundary condition through introducing numerical dissipation inside the solid, in order to damp and drive the solution to the prescribed boundary conditions. 

2) Fully-discrete analyses show that the cause of numerical instability is associated to the secondary modes, in particular the modes coming from the IBM treatment. This analysis further reaches the stability condition of penalization parameter for explicit time marching based on third-order Runge-Kutta scheme: $0.4 \Delta t < \IBMparam_{critical} < 0.5 \Delta t$. This analysis provides guidelines to the selection of the penalization parameter, where a smaller $\IBMparam$ is preferred but should be limited to avoid increased system stiffness and associated numerical stabilities.

3) The proposed analysis shows that the IBM works as a porous medium to absorb and damp the solution inside the solid body. As a result, we propose to add a second order term inside the body. Results show that this approach, in combination with the volume penalization, achieves improved accuracy and reduces the time step restrictions due to stiff penalization terms. For a given penalization parameter, there exists an optimal viscosity coefficient that has the lowest numerical error. The relationship between volume penalization and optimal viscosity coefficient is obtained by curve-fitting.

4) We test our results on a Navier-Stokes flow and find that the proposed second-order damping inside the solid helps to achieve more accurate solutions. Furthermore, the relationship derived using linear analysis, is used as a guideline for choosing the penalty parameters. The cylinder results show improved accuracy, proving the potential of the presented methodology for Navier-Stokes simulations.

\section*{Acknowledgments}
This project has received funding from the European Union’s Horizon 2020 research and innovation programme under the Marie Skłodowska-Curie grant agreement (MSCA ITN-EID-GA ASIMIA No 813605).

\appendix
\section{Discretization of the Advection Equation}
\label{FRdist}
\subsection{Preliminaries}
As the first step of space discretization for a one-dimensional problem, the computational domain $\Domain$ is divided into $N$ nonoverlapping cells, defined as $\Domain_{n} = \left \{ x|x_{n}<x<x_{n+1} \right \}$ with element size $h_n = x_{n+1}-x_{n}$. To implement the space discretization, each cell is transformed from the physical space within the global domain to a common reference space $\Domain_{ref} = \left \{ r|-1<r<1 \right \} $. The mapping between physical space and reference space can be defined as the function $\Gamma$
\begin{equation}
    r = \Gamma_{n}(x) = 2 \left (\frac{x-x_{n}}{x_{n+1}-x_{n}} \right) - 1 ,
\end{equation}
and its inversion
\begin{equation}
    x = \Gamma_{n}^{-1}(r) = \left (\frac{1-r}{2}\right)x_{n} + \left (\frac{1+r}{2}\right)x_{n+1} .
\end{equation}

This transformation leads to the transformed conservation law:
\begin{equation}
    \frac{\partial \Tilde{u}}{\partial t} + \frac{\partial \Tilde{f}}{\partial r} = \frac{\partial u(\Gamma_{n}^{-1}(r),t)}{\partial t} + \frac{1}{J_{n}}\frac{f(\Gamma_{n}^{-1}(r),t)}{\partial r} = 0,
\end{equation}
where the Jacobian $J_{n}$ is defined as $J_{n} = (x_{n + 1} - x_{n}) / 2 $. For each cell, the solution is represented by a polynomial of degree $P$ defined at $P+1$ solution points (quadrature points). The polynomial can be represented by a vector of modal coefficients $\hat{u}_j$ of a Legendre basis, or a vector of nodal coefficients $\Tilde{u}_j$ defined at solution points. Here we retain the former set of basis.

\subsection{Spatial discretization}
The standard FR workflow for a general conservation law includes seven stages \citep{williams2013energy}, whereas for the advection equation it can be reduced to five stages. In the first and the second stage, the polynomial representation in terms of a nodal basis ($\Tilde{u}(r_{i}), i=0,1,...,P$) with order $P$ for both solution and flux is defined at $P+1$ solution points:
\begin{equation}
    \Tilde{u}_n(r) = \sum_{i=0}^{P}l_{i}(r)\Tilde{u}_n(r_{i}) = \boldsymbol{l}^T \Tilde{\boldsymbol{u}}_n,
\end{equation}
\begin{equation}
    \Tilde{f}^{D}_n(r) = \sum_{i=0}^{P}l_{i}(r) \Tilde{f}^{D}_n(r_{i}) = \sum_{i=0}^{P}l_{i}(r) \advcoef \Tilde{u}_n(r_{i}) = \advcoef \boldsymbol{l}^T \Tilde{\boldsymbol{u}}_n,
\end{equation}
where $\Tilde{f}^{D}$ denotes the transformed discontinuous flux computed from the value at solution points for the $n$th cell. $\boldsymbol{l}$ is a Lagrange interpolation operator and $l_{i}(r)$ denotes the Lagrange polynomial defined at a solution point $r_{i}$. This polynomial is defined as the interpolating polynomial through $P+1$ solution points, which takes value $1$ at $r_{i}$ and zero at all other solution points:
\begin{equation}
    l_{i}(r) = \prod^{P}_{j = 0, j \neq i}\left (\frac{r-r_{j}}{r_{i}-r_{j}}\right)
\end{equation}

In this work, Gauss quadrature points are used as the solution points, while another popular choice is based on Gauss-Lobatto quadrature rule (including element end-points). With this representation, the gradient of transformed discontinuous flux at all solution points is obtained as

\begin{equation}
\label{eq:localgrad}
    \frac{\partial \Tilde{f}^{D}_n(r)}{\partial r}
     = \sum_{i=0}^{P}\Tilde{f}^{D}_n(r_{i}) \frac{ \partial l_{i}(r) }{\partial r} = \sum_{i=0}^{P} \advcoef \Tilde{u}_n(r_{i}) \frac{ \partial l_{i}(r) }{\partial r}  = \advcoef \boldsymbol{d}^T \Tilde{\boldsymbol{u}}_n,
\end{equation}
where $\boldsymbol{d}$ is the gradient operator to obtain the discontinuous gradient at position $r$, which can be generalized to gradient matrix $\boldsymbol{D}$ at all solution points $\boldsymbol{D}[i,j]=\frac{ \partial l_{j}(r_i) }{\partial r}$. 

In the third stage of FR, the transformed numerical flux $\Tilde{f} ^{I}$ at all element interfaces (either ends of the cell $\Domain_{r}$, $r= \pm 1$) is computed, which is a common value for the two adjacent cells of each interface. To obtain this flux, the interpolated solution at either end of the interface is needed. For example, for the $n$th cell, we have the left and right interpolated solutions $\Tilde{u}^{L}_n=\boldsymbol{l}^T(-1) \Tilde{\boldsymbol{u}}_n$ and $\Tilde{u}^{R}_n=\boldsymbol{l}^T(1) \Tilde{\boldsymbol{u}}_n$, as well as the interpolated fluxes $\Tilde{f}^{D,L}_n=\advcoef \Tilde{u}^{L}_n$ and $\Tilde{f}^{D,R}_n=\advcoef \Tilde{u}^{R}_n$. The exact method for numerical flux calculation $\Tilde{f}^{I} = \Tilde{f}^{I} (\Tilde{u}_{L}, \Tilde{u}_{R})$ depends on the nature of the equation being solved, like a Roe type approximate Riemann solver \citep{roe1981approximate} or other upwind flux solvers for the Euler equations. Here we use the Lax-Friedrich type numerical flux defined as

\begin{equation}
\label{eq:intfluxl}
    \Tilde{f}^{I,L}_n = \frac{\Tilde{f}^{D,L}_{n} + \Tilde{f}^{D,R}_{n-1}}{2} - \frac{1}{2} \lambda \left | \advcoef \right |  (\Tilde{u}^{L}_{n} - \Tilde{u}^{R}_{n-1}) = \frac{1}{2}(\advcoef+\lambda \left | \advcoef \right |)\boldsymbol{l}^T(1) \Tilde{\boldsymbol{u}}_{n-1} + \frac{1}{2}(\advcoef-\lambda \left | \advcoef \right |)\boldsymbol{l}^T(-1) \Tilde{\boldsymbol{u}}_{n},
\end{equation}

\begin{equation}
\label{eq:intfluxr}
    \Tilde{f}^{I,R}_n = \frac{\Tilde{f}^{D,R}_{n} + \Tilde{f}^{D,L}_{n+1}}{2} - \frac{1}{2} \lambda \left | \advcoef \right |  (\Tilde{u}^{R}_{n} - \Tilde{u}^{L}_{n+1}) = \frac{1}{2}(\advcoef+\lambda \left | \advcoef \right |)\boldsymbol{l}^T(1) \Tilde{\boldsymbol{u}}_{n} + \frac{1}{2}(\advcoef-\lambda \left | \advcoef \right |)\boldsymbol{l}^T(-1) \Tilde{\boldsymbol{u}}_{n+1},
\end{equation}
where $\lambda$ is an upwinding parameter that controls the numerical flux. $\lambda = 1$ defines the upwind flux and $\lambda = 0$ defines the central flux. The upwind flux is considered here. The numerical interface fluxes for both boundaries of the $n$th cell are defined as $\Tilde{f}^{I,L}_{n}$ and $\Tilde{f}^{I,R}_{n}$, respectively.

The forth stage of FR process is to correct the transformed discontinuous flux to be continuous across the cell interface, while approximating the flux function as a polynomial of degree $P + 1$. The transformed correction flux $\Tilde{f}^{C}$ based on a polynomial of degree $P + 1$ is defined to make the sum between $\Tilde{f}^{C}$ and the interpolated discontinuous flux equal to the transformed numerical interface flux at $r \pm 1$: 
\begin{equation}
    \Tilde{f}^{C}_n(-1) = \Tilde{f}^{I,L}_n - \Tilde{f}^{D,L}_n\ , \ \Tilde{f}^{C}_n(1) = \Tilde{f}^{I,R}_n - \Tilde{f}^{D,R}_n.
\end{equation}

To define this correction flux, considering first the correction functions for the left and the right interface, defined as $g^{L} = g^{L}(r)$ and $g^{R} = g^{R}(r)$. This correction satisfies the following condition
\begin{equation}
    g^{L}(-1) = 1 \ , \ g^{L}(1) = 0 \ , \ g^{R}(-1) = 0 \ , \ g^{R}(1) = 1,
\end{equation}
with a symmetry consideration
\begin{equation}
    g^{L}(r) = g^{R}(-r).
\end{equation}

The specific form of correction function is determined by the above-mentioned condition and the stability criterion. Here the left and right Radau polynomials are used for the right and left boundaries, respectively. This type of polynomial belongs to a class of high-order energy stable flux reconstruction (ESFR) schemes that satisfy energy stability criterion. As mentioned in Vincent et al. \citep{vincent2011new}, the resulting scheme is equivalent to a differential formulation of a nodal DG method. The correction function is then written in terms of $g^{L}$ and $g^{R}$ as
\begin{equation}
    \Tilde{f}^{C}_n(r) = (\Tilde{f}^{I,L}_n - \Tilde{f}^{D,L}_n)g^{L}(r) + (\Tilde{f}^{I,R}_n - \Tilde{f}^{D,R}_n)g^{R}(r).
\end{equation}

Therefore the total transformed flux $\Tilde{f}$ with a polynomial of degree $P + 1$ is formed from the discontinuous and correction flux as follows
\begin{equation}
    \Tilde{f}_n = \Tilde{f}^{D}_n + \Tilde{f}^{C}_n.
\end{equation}

In the final stage of FR process, the divergence of the transformed total flux $\frac{\partial\Tilde{f}}{\partial{r}}$ at each solution point is calculated as follows

\begin{equation}
\label{eq:pointflux}
    \frac{\partial\Tilde{f}_n}{\partial{r}}(r_{i}) = \sum_{j=0}^{P}\Tilde{f}_{n}^{D}(r_j)\frac{\partial l_{j}}{\partial r}(r_{i}) + (\Tilde{f}^{I,L}_{n} - \Tilde{f}^{D,L}_{n})\frac{\partial g^{L}}{\partial r}(r_{i}) + (\Tilde{f}^{I,R}_{n} - \Tilde{f}^{D,R}_{n})\frac{\partial g^{R}}{\partial r}(r_{i}).
\end{equation}

This results in the semi-discrete formulation:
\begin{equation}
    \frac{du_n}{dt} = -\frac{1}{J_n} \frac{\partial\Tilde{f}_n}{\partial{r}}.
\end{equation}

Finally, the divergence of total flux can be used to advance the solution $u$ in time with some efficient time integration methods, like the Runge-Kutta method. For viscous fluxes, the Local Discontinuous Galerkin (LDG) \citep{cockburn1998local} formulation can be adopted. In summary, it can be concluded that the performance of the present FR scheme depends on three factors, including the position of solution points, the method for solving the transformed interface flux (the Riemann solver) and the form of correction function. 

\subsection{Local matrix formulation}
The eigensolution analysis relies on the eigendecomposition of the discretization operator in matrix form. Consider uniform element space $h$ for the computational domain, we can obtain the local matrix formulation for the $n$th element as follows, 

\begin{equation}
    \frac{d\boldsymbol{u}_n}{dt} = \boldsymbol{L} \boldsymbol{u}_{n-1} +  \boldsymbol{C} \boldsymbol{u}_{n} +  \boldsymbol{R} \boldsymbol{u}_{n+1},
\end{equation}
where the left, middle and right operators for the current element $\boldsymbol{u}_n$ are defined as $\boldsymbol{L}$, $\boldsymbol{C}$, and $\boldsymbol{R}$, respectively. These matrices can be derived based on Eq. (\ref{eq:localgrad}), Eq. (\ref{eq:intfluxl}), Eq. (\ref{eq:intfluxr}) and Eq. (\ref{eq:pointflux}):

\begin{equation}
    \boldsymbol{L} = - \frac{1}{h} \boldsymbol{g}^L_r \left [ (\advcoef+\lambda \left | \advcoef \right |)\boldsymbol{l}^T(1) \right ],
\end{equation}

\begin{equation}
   \boldsymbol{C} = - \frac{2}{h} \left [ c \boldsymbol{D} - \frac{1}{2} \boldsymbol{g}^L_r  (\advcoef+\lambda \left | \advcoef \right | ) \boldsymbol{l}^T(-1) + \frac{1}{2} \boldsymbol{g}^R_r ( \lambda \left | \advcoef \right | - \advcoef ) \boldsymbol{l}^T(1)  \right],
\end{equation}

\begin{equation}
    \boldsymbol{R} = - \frac{1}{h} \boldsymbol{g}^R_r \left [ (c-\lambda \left | \advcoef \right | )\boldsymbol{l}^T(-1) \right ],
\end{equation}
where the gradient vector of correction function at all solution points is defined as 
\begin{equation}
    \boldsymbol{g}^{dir}_r = \left [ \frac{\partial g^{dir}}{\partial r}(r_{0}), \frac{\partial g^{dir}}{\partial r}(r_{1}),..., \frac{\partial g^{dir}}{\partial r}(r_{P})  \right]^T \ , \ dir = L \ or \ R .
\end{equation}

\bibliographystyle{model1-num-names}
\bibliography{refs}

\end{document}

%% file: commands/cmdFluxReconstruction.tex
\newcommand{\Domain}{\Omega}

%% file: commands/cmdPhysicalQuantities.tex
\newcommand{\velx}{u}

%% file: commands/cmdPunctuation.tex
\newif\ifturb

\newif\ifsutherland

%% file: commands/cmdImmersedBoundary.tex
\newcommand{\mask}{\chi}
\newcommand{\IBMparam}{\eta}

\newcommand{\advcoef}{c}